\documentclass[12pt]{amsart}

\usepackage{amsmath}
\usepackage{amsfonts}
\usepackage{latexsym}
\usepackage{graphicx}
\usepackage{amssymb}
\usepackage{amsthm}
\usepackage[margin=2cm]{geometry}
\usepackage{url}
\usepackage{color}
\usepackage{enumerate}
\usepackage[shortlabels]{enumitem} %
\usepackage[small]{caption}
\usepackage{cite}       %

\usepackage{todonotes}

\usepackage{algorithm}     %
\usepackage[noend]{algpseudocode}
\algrenewcommand\algorithmicrequire{\textbf{Precondition:}}
\algrenewcommand\algorithmicensure{\textbf{Postcondition:}}

\usepackage{array}           %
\usepackage{stmaryrd}

\usepackage{hyperref}

\newcommand{\doi}[1]{
    \href{http://dx.doi.org/#1}{\nolinkurl{doi:#1}}
}

\usepackage{tikz}
\usepackage{tikz-cd}      %
\usetikzlibrary{decorations.pathreplacing} %
\usetikzlibrary{decorations.markings}      %
\usetikzlibrary{decorations.text}          %
\usetikzlibrary{arrows}
\usetikzlibrary{arrows.meta}

\usepackage[T1]{fontenc}     %
\usepackage{lmodern}         %
\usepackage[utf8]{inputenc}  %

\usepackage{sagetex}
\usepackage{accsupp}
\newcommand{\noncopynumber}[1]{
    \BeginAccSupp{method=escape,ActualText={}}
    #1
    \EndAccSupp{}
}

\lstdefinestyle{SageInput}{style=DefaultSageInput,basicstyle={\footnotesize\ttfamily}}
\newtheorem{definition}{Definition}[section]
\newtheorem{lemma}[definition]{Lemma}
\newtheorem{proposition}[definition]{Proposition}

\newtheorem{question}[definition]{Question}
\newtheorem{corollary}[definition]{Corollary}

\newtheorem{remark}[definition]{Remark}
\newtheorem{theorem}[definition]{Theorem}

\newtheorem{mainquestion}{Question}
\newtheorem{maintheorem}{Theorem}

\newtheorem*{THEOREMA}{Theorem~\ref{thm:is-the-metallic-mean-wang-shift}}
\newtheorem*{THEOREMB}{Theorem~\ref{thm:equations-satisfied-by-tiles}}
\newtheorem*{THEOREMC}{Theorem~\ref{thm:construct-valid-tilings}}
\newtheorem*{THEOREMD}{Theorem~\ref{thm:factor-map}}
\newtheorem*{THEOREME}{Theorem~\ref{thm:maximal-equicontinuous-factor}}
\newtheorem*{THEOREMF}{Theorem~\ref{thm:Markov-partition}}

\numberwithin{equation}{section}

\newcommand{\N}{\mathbb{N}}
\newcommand{\Z}{\mathbb{Z}}
\newcommand{\Q}{\mathbb{Q}}
\newcommand{\R}{\mathbb{R}}
\newcommand{\I}{\mathbb{I}} %
\newcommand{\torus}{\mathbb{T}}

\newcommand{\be}{{\boldsymbol{e}}}
\newcommand{\bk}{{\boldsymbol{k}}}

\newcommand{\bm}{{\boldsymbol{m}}}
\newcommand{\bn}{{\boldsymbol{n}}}
\newcommand{\bv}{{\boldsymbol{v}}}
\newcommand{\bx}{{\boldsymbol{x}}}
\newcommand{\zero}{{\boldsymbol{0}}}

\newcommand{\Acal}{\mathcal{A}}
\newcommand{\Bcal}{\mathcal{B}}
\newcommand{\Ccal}{\mathcal{C}}
\newcommand{\Fcal}{\mathcal{F}}
\newcommand{\Lcal}{\mathcal{L}}
\newcommand{\Pcal}{\mathcal{P}}
\newcommand{\Tcal}{\mathcal{T}}
\newcommand{\Ucal}{\mathcal{U}}
\newcommand{\Xcal}{\mathcal{X}}

\newcommand{\ibar}{\overline{i}}
\newcommand{\jbar}{\overline{j}}
\newcommand{\nbar}{\overline{n}}

\newcommand{\east}{\textsc{East}}
\newcommand{\west}{\textsc{West}}
\newcommand{\north}{\textsc{North}}
\newcommand{\south}{\textsc{South}}

\newcommand{\sctile}{\textsc{Tile}}
\newcommand{\sctop}{\textsc{Top}}
\newcommand{\scbottom}{\textsc{Bottom}}
\newcommand{\scright}{\textsc{Right}}
\newcommand{\scleft}{\textsc{Left}}

\newcommand{\dynsys}[3]{#1\overset{#2}{\curvearrowright}#3}
\newcommand{\shiftclosure}[1]{{\overline{#1}^{\sigma}}}

\newcommand{\dist}{\mathrm{dist}}
\newcommand{\card}{\mathrm{card}}

\newcommand{\Int}[1]{\operatorname{\mathrm{Interior}}\left(#1\right)}

\newcommand\defn[1]{\textbf{#1}}

\newcommand\labeloutside[6]{
\draw[draw] (#1+\size,#2) -- node[swap]{$#3$} (#1+\size,#2+\size)
                          -- node[swap]{$#4$} (#1,#2+\size)
                          -- node[swap]{$#5$} (#1,#2)
                          -- node[swap]{$#6$} (#1+\size,#2);
}

\def\dt{.25}

\newcommand\tile[7]{
\begin{scope}
\draw[draw=none,fill=#1] (#2,#3) rectangle (#2+\size,#3+\size);
\labeloutside{#2}{#3}{#4}{#5}{#6}{#7}
\end{scope}}

\newcommand\tileV[7]{
\begin{scope}[xshift=#2cm,yshift=#3cm]
\draw[draw=none,fill=#1] (#2+\dt,#3) rectangle (#2+\size-\dt,#3+\size);
\labeloutside{#2}{#3}{#4}{#5}{#6}{#7}
\end{scope}}

\newcommand\tileH[7]{
\begin{scope}[xshift=#2cm,yshift=#3cm]
\draw[draw=none,fill=#1] (#2,#3+\dt) rectangle (#2+\size,#3+\size-\dt);
\labeloutside{#2}{#3}{#4}{#5}{#6}{#7}
\end{scope}}

\def\size{1}
\def\widthYellow{.5}
\def\widthCyan{.5}
\def\greenEllipseRadius{.7}
\def\pinkEllipseRadius{.4}
\def\ourColorGreen{green!80}
\def\ourColorYellow{yellow}
\def\ourColorBlue{cyan!70}

\newcommand\tileHgreenBackground[2]{
\draw[fill=\ourColorYellow,draw=none]
({#1+\size},{#2+.5*\size+.5*\widthYellow}) arc [start angle=90, end angle=270,
x radius=\greenEllipseRadius, y radius={.5*\widthYellow}] -- cycle;
\draw[fill=\ourColorBlue,draw=none]
({#1},{#2+.5*\size+.5*\widthCyan}) arc [start angle=90, end angle=-90,
x radius=\greenEllipseRadius, y radius={.5*\widthCyan}] -- cycle;
\begin{scope}
\clip
({#1+\size},{#2+.5*\size+.5*\widthYellow}) arc [start angle=90, end angle=270,
x radius=\greenEllipseRadius, y radius={.5*\widthYellow}] -- cycle;
\draw[fill=\ourColorGreen,draw=none]
({#1},{#2+.5*\size+.5*\widthCyan}) arc [start angle=90, end angle=-90,
x radius=\greenEllipseRadius, y radius={.5*\widthCyan}] -- cycle;
\end{scope}
}

\newcommand\tileHpinkBackground[2]{
\draw[fill=\ourColorBlue,draw=none]
({#1+\size},{#2+.5*\size+.5*\widthCyan}) arc [start angle=90, end angle=270,
x radius=\pinkEllipseRadius, y radius={.5*\widthCyan}] -- cycle;
\draw[fill=\ourColorYellow,draw=none]
({#1},{#2+.5*\size+.5*\widthYellow}) arc [start angle=90, end angle=-90,
x radius=\pinkEllipseRadius, y radius={.5*\widthYellow}] -- cycle;
\begin{scope}
\clip
({#1+\size},{#2+.5*\size+.5*\widthCyan}) arc [start angle=90, end angle=270,
x radius=\pinkEllipseRadius, y radius={.5*\widthCyan}] -- cycle;
\draw[fill=pink,draw=none]
({#1},{#2+.5*\size+.5*\widthYellow}) arc [start angle=90, end angle=-90,
x radius=\pinkEllipseRadius, y radius={.5*\widthYellow}] -- cycle;
\end{scope}
}

\newcommand\tileVgreenBackground[2]{
\draw[fill=\ourColorYellow,draw=none]
({#1+.5*\size+.5*\widthYellow},{#2+\size}) arc [start angle=0, end angle=-180,
x radius={.5*\widthYellow}, y radius=\greenEllipseRadius] -- cycle;
\draw[fill=\ourColorBlue,draw=none]
({#1+.5*\size+.5*\widthCyan},{#2}) arc [start angle=0, end angle=180,
x radius={.5*\widthCyan}, y radius=\greenEllipseRadius] -- cycle;
\begin{scope}
\clip
({#1+.5*\size+.5*\widthYellow},{#2+\size}) arc [start angle=0, end angle=-180,
x radius={.5*\widthYellow}, y radius=\greenEllipseRadius] -- cycle;
\draw[fill=\ourColorGreen,draw=none]
({#1+.5*\size+.5*\widthCyan},{#2}) arc [start angle=0, end angle=180,
x radius={.5*\widthCyan}, y radius=\greenEllipseRadius] -- cycle;
\end{scope}
}

\newcommand\tileVpinkBackground[2]{
\draw[fill=\ourColorBlue,draw=none]
({#1+.5*\size+.5*\widthCyan},{#2+\size}) arc [start angle=0, end angle=-180,
x radius={.5*\widthCyan}, y radius=\pinkEllipseRadius] -- cycle;
\draw[fill=\ourColorYellow,draw=none]
({#1+.5*\size+.5*\widthYellow},{#2}) arc [start angle=0, end angle=180,
x radius={.5*\widthYellow}, y radius=\pinkEllipseRadius] -- cycle;
\begin{scope}
\clip
({#1+.5*\size+.5*\widthCyan},{#2+\size}) arc [start angle=0, end angle=-180,
x radius={.5*\widthCyan}, y radius=\pinkEllipseRadius] -- cycle;
\draw[fill=pink,draw=none]
({#1+.5*\size+.5*\widthYellow},{#2}) arc [start angle=0, end angle=180,
x radius={.5*\widthYellow}, y radius=\pinkEllipseRadius] -- cycle;
\end{scope}
}

\newcommand\tileJunctionBackground[6]{
\def\widthR{0.5}
\def\widthT{0.5}
\def\widthL{0.5}
\def\widthB{0.5}
\coordinate (R) at (#1+\size,{#2+(\size+\widthR)/2}) {};
\coordinate (T) at ({#1+(\size+\widthT)/2},#2+\size) {};
\coordinate (L) at (#1,{#2+(\size-\widthL)/2}) {};
\coordinate (B) at ({#1+(\size-\widthB)/2},#2) {};
\def\du{.02}
\draw [fill=#3,draw=none] (R) to [bend right=20] (#1+\size/2+\du,#2+\size/2-\du) -- (#1+\size,#2) -- cycle;
\draw [fill=#4,draw=none] (T) to [bend left=20]  (#1+\size/2-\du,#2+\size/2+\du) -- (#1,#2+\size) -- cycle;
\draw [fill=#5,draw=none] (L) to [bend right=20] (#1+\size/2-\du,#2+\size/2+\du) -- (#1,#2+\size) -- cycle;
\draw [fill=#6,draw=none] (B) to [bend left=20]  (#1+\size/2+\du,#2+\size/2-\du) -- (#1+\size,#2) -- cycle;
\draw[fill=white,draw=none] ({#1},{#2+\size})
-- ++ (0,{-(\size-\widthL)/2}) arc
[start angle=-90, end angle=0,
x radius={(\size-\widthT)/2}, y radius={(\size-\widthL)/2}];
\draw[fill=white,draw=none] ({#1+\size},{#2})
-- ++ (0,{(\size-\widthR)/2}) arc
[start angle=90, end angle=180,
x radius={(\size-\widthB)/2}, y radius={(\size-\widthR)/2}];
}

\newcommand\crossX[2]{
\draw[very thick] ({#1+\size*.4},{#2+\size*.4}) --
                    ({#1+\size*.6},{#2+\size*.6});
\draw[very thick] ({#1+\size*.4},{#2+\size*.6}) --
                    ({#1+\size*.6},{#2+\size*.4});
}

\newcommand\tileJunctionOOOO[6]{
\begin{scope}
\tileJunctionBackground{#1}{#2}{\ourColorBlue}{\ourColorBlue}{\ourColorBlue}{\ourColorBlue}
\labeloutside{#1}{#2}{#3}{#4}{#5}{#6}
\end{scope}}
\newcommand\tileJunctionOOOI[6]{
\begin{scope}
\tileJunctionBackground{#1}{#2}{\ourColorBlue}{\ourColorBlue}{\ourColorBlue}{\ourColorYellow}
\labeloutside{#1}{#2}{#3}{#4}{#5}{#6}
\end{scope}}
\newcommand\tileJunctionOOIO[6]{
\begin{scope}
\tileJunctionBackground{#1}{#2}{\ourColorBlue}{\ourColorBlue}{\ourColorYellow}{\ourColorBlue}
\labeloutside{#1}{#2}{#3}{#4}{#5}{#6}
\end{scope}}
\newcommand\tileJunctionOOII[6]{
\begin{scope}
\tileJunctionBackground{#1}{#2}{\ourColorBlue}{\ourColorBlue}{\ourColorYellow}{\ourColorYellow}
\labeloutside{#1}{#2}{#3}{#4}{#5}{#6}
\end{scope}}

\newcommand\tileJunctionOIIO[6]{
\begin{scope}
\tileJunctionBackground{#1}{#2}{\ourColorBlue}{\ourColorYellow}{\ourColorYellow}{\ourColorBlue}
\labeloutside{#1}{#2}{#3}{#4}{#5}{#6}
\crossX{#1}{#2}
\end{scope}}
\newcommand\tileJunctionOIII[6]{
\begin{scope}
\tileJunctionBackground{#1}{#2}{\ourColorBlue}{\ourColorYellow}{\ourColorYellow}{\ourColorYellow}
\labeloutside{#1}{#2}{#3}{#4}{#5}{#6}
\end{scope}}

\newcommand\tileJunctionIOOI[6]{
\begin{scope}
\tileJunctionBackground{#1}{#2}{\ourColorYellow}{\ourColorBlue}{\ourColorBlue}{\ourColorYellow}
\labeloutside{#1}{#2}{#3}{#4}{#5}{#6}
\crossX{#1}{#2}
\end{scope}}

\newcommand\tileJunctionIOII[6]{
\begin{scope}
\tileJunctionBackground{#1}{#2}{\ourColorYellow}{\ourColorBlue}{\ourColorYellow}{\ourColorYellow}
\labeloutside{#1}{#2}{#3}{#4}{#5}{#6}
\end{scope}}

\newcommand\tileJunctionIIII[6]{
\begin{scope}
\tileJunctionBackground{#1}{#2}{\ourColorYellow}{\ourColorYellow}{\ourColorYellow}{\ourColorYellow}
\labeloutside{#1}{#2}{#3}{#4}{#5}{#6}
\end{scope}}

\newcommand\tileJunctionGRAY[7]{
\begin{scope}
\tileJunctionBackground{#1}{#2}{#7}{#7}{#7}{#7}
\labeloutside{#1}{#2}{#3}{#4}{#5}{#6}
\end{scope}}

\newcommand\tileHgreen[6]{
\begin{scope}
\tileHgreenBackground{#1}{#2}
\labeloutside{#1}{#2}{#3}{#4}{#5}{#6}
\end{scope}}

\newcommand\tileVgreen[6]{
\begin{scope}
\tileVgreenBackground{#1}{#2}
\labeloutside{#1}{#2}{#3}{#4}{#5}{#6}
\end{scope}}

\newcommand\tileHantigreen[6]{
\begin{scope}
\tileHpinkBackground{#1}{#2}
\labeloutside{#1}{#2}{#3}{#4}{#5}{#6}
\end{scope}}

\newcommand\tileVantigreen[6]{
\begin{scope}
\tileVpinkBackground{#1}{#2}
\labeloutside{#1}{#2}{#3}{#4}{#5}{#6}
\end{scope}}

\begin{document}

\title
[Metallic mean Wang tiles II: the dynamics of an aperiodic computer chip]
{Metallic mean Wang tiles II:\\the dynamics of an aperiodic computer chip}

\author[S.~Labb\'e]{S\'ebastien Labb\'e}
\address[S.~Labb\'e]{CNRS, LaBRI, UMR 5800, F-33400 Talence, France}
\email{sebastien.labbe@labri.fr}
\urladdr{http://www.slabbe.org/}

\makeatletter
\@namedef{subjclassname@2020}{\textup{2020} Mathematics Subject Classification}
\makeatother


\date{\today}

\begin{abstract}
We consider a new family $(\mathcal{T}_n)_{n\geq1}$ of aperiodic sets of Wang tiles and we describe the dynamical properties of the set $\Omega_n$ of valid configurations $\mathbb{Z}^2\to\mathcal{T}_n$. The tiles can be defined as the different instances of a square-shaped computer chip whose inputs and outputs are 3-dimensional integer vectors. The family include the Ammann aperiodic set of 16 Wang tiles and gathers the hallmarks of other small aperiodic sets of Wang tiles. Notably, the tiles satisfy additive versions of equations verified by the Kari--Culik aperiodic sets of 14 and 13 Wang tiles. Also configurations in $\Omega_n$ are the codings of a $\mathbb{Z}^2$-action on a 2-dimensional torus like the Jeandel--Rao aperiodic set of 11 Wang tiles. The family broadens the relation between quadratic integers and aperiodic tilings beyond the omnipresent golden ratio as the dynamics of $\Omega_n$ involves the positive root $\beta$ of the polynomial $x^2-nx-1$, also known as the $n$-th metallic mean. We show the existence of an almost one-to-one factor map $\Omega_n\to\mathbb{T}^2$ which commutes with the shift action on $\Omega_n$ with horizontal and vertical translations by $\beta$ on $\mathbb{T}^2$. The factor map can be explicitly defined by the average of the top labels from the same row of tiles as in Kari and Culik examples. The proofs are based on the minimality of $\Omega_n$ (proved in a previous article) and a polygonal partition of $\mathbb{T}^2$ which we show is a Markov partition for the toral $\mathbb{Z}^2$-action. The partition and the sets of Wang tiles are symmetric which makes them, like Penrose tilings, worthy of investigation.
\end{abstract}

\maketitle

\setcounter{tocdepth}{1}
\tableofcontents

\section{Introduction}

Turing machines can be encoded into a finite set of Wang
tiles (unit squares with labeled edges) in such a way that the Turing machine
does not halt if and only if there exists a tiling of the plane by translated
copies of the tiles respecting the condition that the common edge of
adjacent tiles have the same label \cite{MR0216954},
see also \cite{MR0297572,MR2507046,MR4200110}.
As a consequence, the existence of a valid tiling of the plane with a given
finite set of Wang tiles (called the domino problem) can not be decided by an algorithm. 
Indeed, if the domino problem were decidable, we could use the algorithm
solving the domino problem to solve the halting problem, which is a
contradiction \cite{MR1577030}.

Therefore, we can think of Wang tiles as if their tilings are computing something.
As observed by Wang, the undecidability of the domino problem implies the
existence of aperiodic sets of Wang tiles \cite{wang_proving_1961}.
Shortly after, Berger proved the undecidability of the domino problem 
and constructed the first known aperiodic set of Wang tiles \cite{MR0216954}.
Since then, aperiodic tilings has developed into an active subject of study 
with applications to the theory of quasicrystals
\cite{MR857454,zbMATH00768067,MR3136260,zbMATH06785522}.
Thus, sets of Wang tiles (and their computations) can be classified into three cases:
\begin{itemize}
    \item \textbf{Finite:}    the Wang tiles do not tile the plane,
    \item \textbf{Periodic:}  the Wang tiles tile the plane and one of the valid tiling is periodic,
    \item \textbf{Aperiodic:} the Wang tiles tile the plane and none of the valid tilings are periodic.
\end{itemize}
The finite cases can be associated with computations that halt.
The periodic cases can be associated with computations that do not halt and
fall into an infinite loop.
The aperiodic cases can be associated with computations that do not halt and
never repeat.

For applications, computations that halt are usually preferred over
computations that loop forever.
Among computations that halt, the description of those ``busy beavers'' 
\cite{zbMATH04123308,10.1145/3427361.3427369}
running for the maximum number of steps before halting 
is an open question even for Turing machines made of only 6 rules \cite{OEISA060843}
(it was recently solved for 5 rules\footnote{\url{https://github.com/ccz181078/Coq-BB5}}).
In this article, we are interested in the description of computations that do
not halt and never repeat. We focus on those that happen to be performed by
small aperiodic sets of Wang tiles. 
We aim to reveal their links with dynamical systems and the coding of their orbits.

\subsection*{The Kari--Culik outliers}
\begin{figure}
    \begin{minipage}{.35\textwidth}
        \centering
    \[
        g(x)=
        \begin{cases}
            2x           & \hspace{-3mm}\text{ if } x\leq1,\\
            \frac{2}{3}x & \hspace{-3mm}\text{ if } x>1.
        \end{cases}
    \]
    \includegraphics{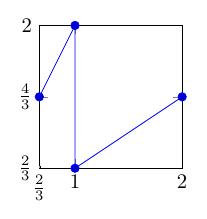}
    \end{minipage}
    \begin{minipage}{.64\textwidth}
\begin{center}
    \resizebox{.97\linewidth}{!}{
        \includegraphics{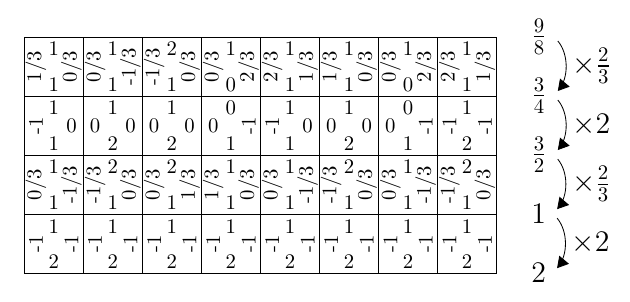}
    }
\end{center}
\end{minipage}
    \caption{Averages of horizontal labels in a tiling with Kari's 14 tiles are
    orbits under the map $g$ on the interval $[\frac{2}{3},2]$; see
    \cite{MR3606059,zbMATH06754662}.}
    \label{fig:Kari-tiling}
\end{figure}

The smallest sets of aperiodic Wang tiles until 2015
were discovered by Kari and Culik in 1996.
Kari \cite{MR1417578} proved that a well-chosen set of 14 Wang tiles admits
tilings of the plane, and that none of them is periodic.
The proof that they are not periodic is cleverly short. It is based on an
arithmetic interpretation of the edge labels of the Wang tiles.
The tiles have labels $r,t,\ell,b\in\Q$ 
satisfying an equation
\begin{equation}\label{eq:wang-kari}
    \raisebox{-9.5mm}{
    \begin{tikzpicture}[auto]
    \tile{white}{0}{0}{r}{t}{\ell}{b}
    \end{tikzpicture}
    }
    \qquad
    \qquad
    qt+\ell = b + r
\end{equation}
for some $q\in\Q$.
We may interpret the Wang tile as a computation (the multiplication by $q$) with
value $t$ as an input and $b$ as an output. 
The value $\ell$ is a carry input on the left and $r$ is a carry output on the right.
Kari \cite{MR1417578} proposed a set of four tiles
satisfying \eqref{eq:wang-kari} with $q=2$ and ten tiles with $q=\frac{2}{3}$. 
The proof of the non-existence of a periodic tiling with those 14 tiles
uses the fact that the equation $2^m 3^n = 1$
has only one solution over the integers ($m=n=0$), see Figure~\ref{fig:Kari-tiling}.
Based on the same idea, Culik \cite{MR1417576} proposed a smaller aperiodic set
of 13 tiles (four tiles satisfying \eqref{eq:wang-kari} with $q=3$ and nine
tiles with $q=\frac{1}{2}$). 
Note that generalizations of Kari--Culik tilings exist \cite{MR2369448} 
and that further results were obtained about 
their entropy \cite{MR3606059}
and on a minimal subsystem \cite{MR3668002}.

Among aperiodic tilings of the plane by Wang tiles, Kari and Culik sets seem like outliers.
The aperiodicity of Penrose tiles \cite{zbMATH03663938},
Berger tiles \cite{MR0216954},
Robinson tiles \cite{MR0297572},
Knuth tiles \cite{MR0286317},
Ammann tiles \cite{MR857454,MR1156132}
can be explained by the hierarchical decomposition of their tilings.
Often, aperiodic tilings have a self-similar structure
\cite{MR1452190,MR1637896,MR1854103,MR1615950,zbMATH07258478}
and this is the case for recently discovered aperiodic geometrical tiles
\cite{MR2834173,MR4770585,MR4807152}.
However, Kari and Culik tilings have positive entropy. Thus, they are not
self-similar and do not possess a hierarchical decomposition 
\cite{MR3606059}. Note that the absence of hierarchical decomposition also 
follows from a cylindricity argument
proposed by Thierry Monteil and explained in \cite[\S 4.2]{MR3606059}.
Moreover, except some extensions of Kari and Culik sets \cite[\S 6]{MR2369448},
no other known aperiodic sets of tiles satisfy equations explaining their non-periodicity.

\subsection*{The metallic mean family of aperiodic Wang tiles}

The current article is the second article about a new family of aperiodic
Wang tiles related to the metallic mean. 
Recall that the metallic mean $\beta$
is the positive root 
of the polynomial $x^2-nx-1$ where $n\geq1$ is an integer \cite{MR1729917}, 
that is,
\[
    \beta = [n; n, n, \cdots] =
n + \frac{1}{
n + \frac{1}{
n + \frac{1}{
n + \cdots}}}
= n + \frac{1}{\beta}.
\]
Metallic means were also called \emph{silver means} in \cite{zbMATH00051561}
and \emph{noble means} in \cite{MR3136260}.

\begin{figure}[h]
\begin{center}
    \includegraphics{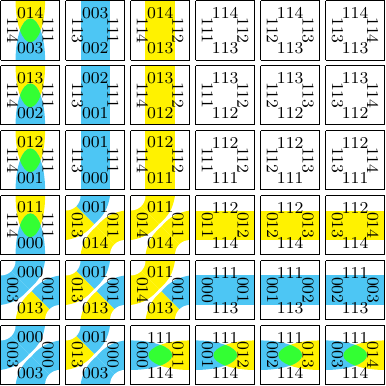}
\end{center}
    \caption{The metallic mean Wang tile set $\Tcal_n$ for $n=3$.}
    \label{fig:T3}
\end{figure}

\begin{figure}[h]
\begin{center}
    \includegraphics[width=.8\linewidth]{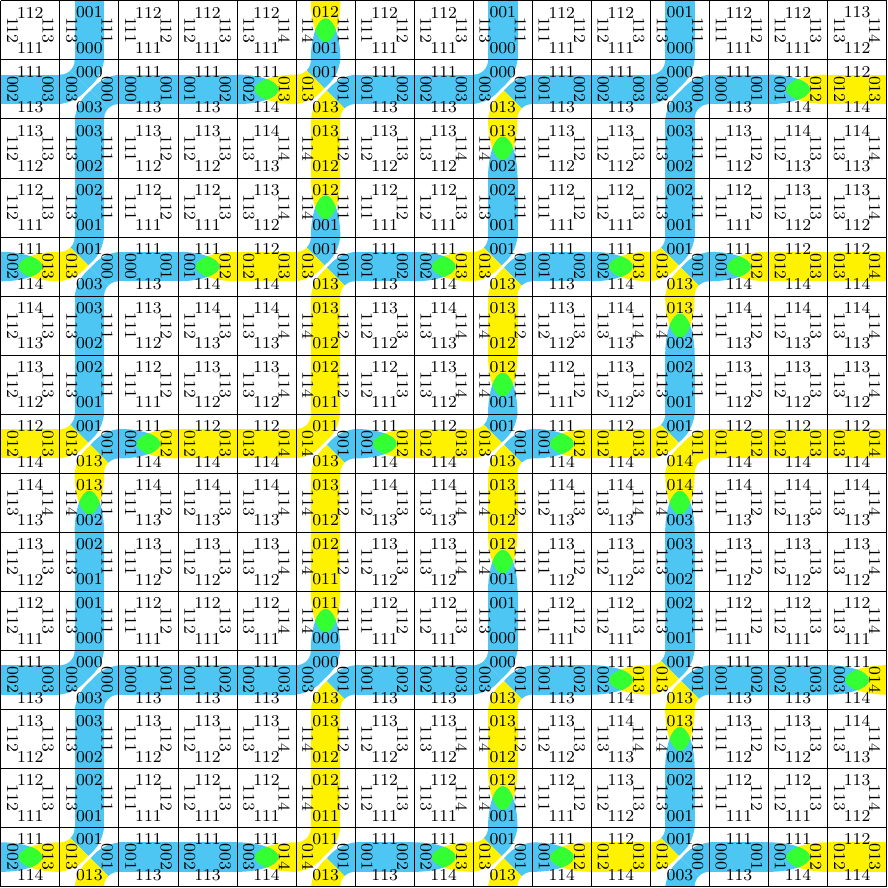}
\end{center}
    \caption{A valid $15\times 15$ pattern with Wang tile set $\Tcal_3$.}
    \label{fig:T3-15x15}
\end{figure}

Let us recall the main results proved in the first article of the series.
For every integer $n\geq1$,
the $n^{th}$ metallic mean Wang shift $\Omega_n$ is defined from a set $\Tcal_n$
of $(n+3)^2$ Wang tiles. 
An illustration of the set $\Tcal_3$ is shown in Figure~\ref{fig:T3}.
The labels of the Wang tiles are vectors in $\N^3$.
In Figure~\ref{fig:T3}, we represent vectors as words
for economy of space reasons. 
For instance, the vector $(1,1,4)$ is represented as $114$.
A finite rectangular valid tiling is shown in
Figure~\ref{fig:T3-15x15} for the set $\Tcal_3$.
More images of valid tilings with metallic mean Wang tiles are available in
\cite{labbe_metallic_I_2025}.

It was shown in the previous article that the metallic mean Wang shift
$\Omega_n$ is self-similar, aperiodic and minimal.
We gather in the next theorem the main results already proved about $\Omega_n$.

\begin{theorem}[\cite{labbe_metallic_I_2025}]
    For every integer $n\geq1$,
    \begin{enumerate}[(i)]
        \item the metallic mean Wang shift $\Omega_n$ is self-similar,
            aperiodic and minimal,
        \item the inflation factor of the self-similarity of $\Omega_n$ is
            the $n$-th metallic mean,
            that is, the positive root of $x^2-nx-1$.
    \end{enumerate}
    Also, when $n=1$, $\Omega_1$ is equivalent to the Wang shift defined
    from the 16 Ammann Wang tiles \cite[p.595, Figure 11.1.13]{MR857454}.
\end{theorem}

In order to describe the substitutive structure of the Wang shift $\Omega_n$ generated from
the set $\Tcal_n$, it was needed in \cite{labbe_metallic_I_2025}
to introduce a larger set $\Tcal_n'$ satisfying $\Tcal_n\subseteq\Tcal_n'$.
It was shown that the set $\Tcal_n'$ 
is in bijection with the set of possible return blocks allowing to
decompose uniquely the configurations of $\Omega_n$.
The return blocks are rectangular blocks of
tiles with a unique junction tile (a tile where horizontal and vertical color
stripes intersect) at the lower left corner.
Also, it was proved in \cite{labbe_metallic_I_2025} that in a valid
configuration of $\Omega_n'$, only the tiles from $\Tcal_n$ appear.
From this observation follows the self-similarity of $\Omega_n$.

\subsection*{This article}

\begin{figure}
\begin{center}
    \includegraphics{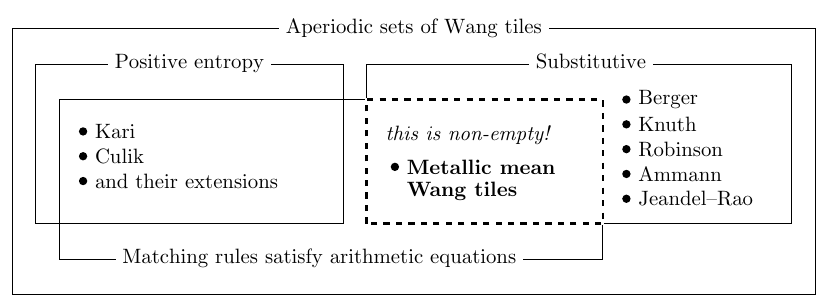}
\end{center}
    \caption{A Venn diagram of aperiodic sets of Wang tiles.
    Aperiodicity of Kari \cite{MR1417578} and Culik \cite{MR1417576} sets of
    tiles and their extensions \cite{MR2369448}
    follows from the arithmetic equations satisfied by their matching rules.
    In this article, we show that the dashed region in the Venn diagram
    is non-empty, that is, there exists a family of
    substitutive (self-similar) aperiodic sets of Wang tiles 
    whose matching rules satisfy arithmetic equations.}
    \label{fig:venn-diagram}
\end{figure}

In this article, we demonstrate that Kari and Culik tilings are not a complete
oddity within aperiodic sets of tiles. In particular, we show for the first
time that substitutive aperiodic sets of Wang tiles can also satisfy equations
and even be defined by them, see Figure~\ref{fig:venn-diagram}.
This article is devoted to a family of aperiodic Wang tiles associated with the
metallic mean numbers, the positive roots of the polynomials $x^2-nx-1$ where
$n\geq1$ is a positive integer. When $n=1$, the family recovers
the Ammann set of 16 Wang tiles \cite{MR857454}.

The labels of the Wang tiles are not numbers like in Kari and Culik sets, but
rather integer vectors.
Note that integers vectors were already used as labels of Wang tiles
in \cite{zbMATH05523458,kari_undecidability_2008}, see also \cite{zbMATH06754662}.
The equations satisfied by the tiles are derived from a function that expresses a
relation between the labels of the Wang tiles. 
The function
provides an independent definition of the family of metallic mean Wang tiles
as the instances of an aperiodic computer chip.
The family $(\Omega_n)_{n\geq1}$ of metallic mean Wang shifts
was introduced separately in \cite{labbe_metallic_I_2025} where
it was shown to be aperiodic as a consequence of its self-similarity.

Here, in this second article on the metallic mean Wang tiles, we prove
that $\Omega_n$ is aperiodic for another reason.
Namely, we show that the $\Z^2$ shift action on $\Omega_n$ is an almost 1-to-1
extension of a minimal $\Z^2$-action by rotations on $\torus^2$.
This reminds of a result proved for Penrose tilings \cite{MR1355301}
and the two reasons for them to be aperiodic.
Aperiodicity of Penrose tilings follows from its self-similarity \cite{zbMATH03663938}
and from the fact of being a cut-and-project scheme \cite{MR609465,MR3136260}.

For every integer $n\geq1$,
we show that valid configurations in $\Omega_n$
are computing the orbits of a dynamical system defined by a $\Z^2$-action $R_n$
on the 2-dimensional torus $\torus^2$.
The dynamical system $\dynsys{\Z^2}{R_n}{\torus^2}$
is defined by horizontal and vertical translation on
$\torus^2$ by the $n$-th metallic mean modulo 1.
As for the Jeandel--Rao Wang shift \cite{MR4213162}, the proof is based 
on a polygonal partition of $\torus^2$ which we prove is a Markov partition for
the toral $\Z^2$-action. We also prove the existence of an almost one-to-one factor map
$\Omega_n\to\torus^2$ commuting the shift
    $\dynsys{\Z^2}{\sigma}{\Omega_n}$
    with the toral $\Z^2$-rotation 
    $\dynsys{\Z^2}{R_n}{\torus^2}$.
    Since $R_n$ is a free action, 
    this provides a second reason for the Wang shift $\Omega_n$ to be aperiodic.

The factor map can be defined by taking averages of 
the dot product involving the top labels of the Wang tiles in the biinfinite
row of tiles passing through the origin in a configuration.
The existence of the factor map proves that the average changes from row to row
by an irrational rotation by the $n$-th metallic mean number.
This can be seen as an additive version of a multiplicative phenomenon known
for Kari--Culik tilings.
Recall that the average of top label values along a row is at the heart of Kari
and Culik's construction of aperiodic tilings where the average change
by a rational multiplication from row to row \cite[Theorem 6]{MR3606059}.

The polygonal partition used to encode the toral $\Z^2$-action is symmetric and
is much more simple to define compared to the Markov partition associated with the
Jeandel--Rao Wang shift.
Moreover, the label of the polygonal atoms of the partition have a meaning in
the sense that they define the linear inequalities describing their boundaries.
The symmetry and simplicity of the partition was helpful to extend the family
beyond the golden ratio. The results proved here for the metallic mean Wang
tiles should serve as an inspiration to replace the labels of the Jeandel--Rao
tiles by integer vectors satisfying equations. 
Understanding the matching rules of Jeandel--Rao tiles by mean
of arithmetic would open the door for discovering a vast family of
aperiodic sets of Wang tiles beyond the family of metallic mean Wang tiles.
See Section~\ref{sec:open-questions} for more open questions.

\subsection*{Structure of the article}
In Section~\ref{sec:statements-main-results}, we state the main results proved
in this article.
In Section~\ref{sec:prelim-wang-shifts},
we present preliminary notions on dynamical systems, subshifts and Wang shifts.
In Section~\ref{sec:family-metallic}, we recall the definition
of the family of metallic mean Wang tiles.
In Section~\ref{sec:thetan-chip},
we show that instances of the $\theta_n$-chip are the metallic mean Wang tiles.
This proves Theorem~\ref{thm:is-the-metallic-mean-wang-shift}.
In Section~\ref{sec:equations},
we prove Theorem~\ref{thm:equations-satisfied-by-tiles}
and we present more equations satisfied by the metallic mean tiles and their tilings.
In Section~\ref{sec:valid-tilngs-as-codings},
we use the floor function on linear forms to construct valid tilings
with the metallic mean Wang tiles and we prove
Theorem~\ref{thm:construct-valid-tilings}.
In Section~\ref{sec:explicit-factor-map},
we define an explicit factor map $\Omega_n\to\torus^2$
and we prove Theorem~\ref{thm:factor-map}.
In Section~\ref{sec:isomorphism},
we define the partition $\Pcal_n$ for every integer $n\geq1$ and
we show that the metallic mean Wang shift is equal to the symbolic dynamical system
defined by the coding of a toral $\Z^2$-action by this partition.
This shows that $\Omega_n$ is isomorphic
as measure-preserving dynamical systems to a toral $\Z^2$-action.
We prove
Theorem~\ref{thm:Markov-partition} and
Theorem~\ref{thm:maximal-equicontinuous-factor}
in this section.
In Section~\ref{sec:renormalization-rauzy},
we compute the renormalization of the partition $\Pcal_n$ and $\Z^2$-action $R_n$
using computations performed in SageMath when $n=3$.
We illustrate how the Rauzy induction of $\Z^2$-actions
and of polygonal partitions can be used to show the self-similarity
of the symbolic dynamical system $\Xcal_{\Pcal_n,R_n}$.
In Section~\ref{sec:open-questions}, we discuss some open questions raised by
the current work.

\section{Statements of the main results}\label{sec:statements-main-results}

\subsection*{An aperiodic computer chip}

\begin{figure}[h]
\begin{center}
    \includegraphics{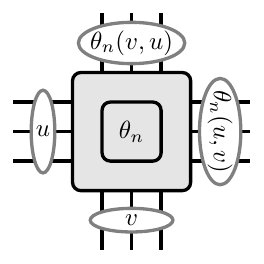}
\end{center}
    \caption{The $\theta_n$-chip is a computer chip computing $\theta_n(u,v)$ and $\theta_n(v,u)$ 
             from the left input $u$ and bottom input $v$.}
    \label{fig:chip}
\end{figure}

For every integer $n\geq1$, we define a finite subset $V_n\subset\N^3$ of vectors
\[
    V_n = \{(v_0,v_1,v_2)\in\N^3\colon 
                       0\leq v_0\leq v_1\leq 1 
                       \text{ and } 
                       v_1\leq v_2\leq n+1\}
\]
with nondecreasing entries where the middle entry is at most 1.
We introduce a function 
\[
\begin{array}{rccl}
    \theta_n:&V_n\times V_n & \to & \Z^3\\
    &(u_0,u_1,u_2), (v_0,v_1,v_2) 
    & \mapsto & (r_0,r_1,r_2),
\end{array}
\]
taking two vectors as input and returning one vector.
Its image is defined by the rule
\begin{equation}\label{eq:defintion-map-theta}
\left\{
\begin{array}{ll}
    &r_0=u_0,\\
    &r_1=\begin{cases}
            v_2-n       & \text{ if } u_0 = 0,\\
            1  & \text{ if } u_0 = 1,
         \end{cases}\\
    &r_2=\begin{cases}
            v_1+u_0             & \text{ if } v_0 = 0,\\
            u_2+1  & \text{ if } v_0 = 1.
         \end{cases}
\end{array}
\right.
\end{equation}
Notice that $(r_0,r_1,r_2)$ does not depend on $u_1$.
For every integer $n\geq1$, we construct a symmetric $\theta_n$-chip, that is, a computer
chip taking as inputs $u\in V_n$ on the left and $v\in V_n$ on the bottom and
producing as outputs $\theta_n(u,v)$ on the right and $\theta_n(v,u)$ on the
top (see Figure~\ref{fig:chip}).

\begin{figure}
\begin{center}
    \includegraphics{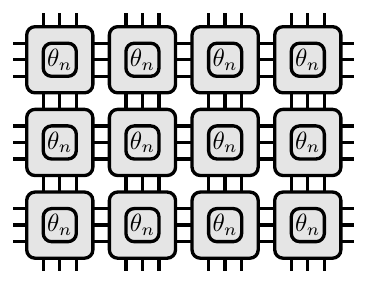}
\end{center}
    \caption{A rectangular cluster of copies of the $\theta_n$-chip.
    }
    \label{fig:chip-cluster}
\end{figure}

If $\theta_n(u,v)$ and $\theta_n(v,u)$ are in $V_n$, then
one can use multiple copies of the $\theta_n$-chip and connect them to each
other horizontally and vertically into an arbitrarily large rectangular cluster
of $\theta_n$-chips (see Figure~\ref{fig:chip-cluster}).

We prove in this work the existence of arbitrarily large rectangular clusters of the $\theta_n$-chip
all of them performing correct computations.
Also we show that no rectangular cluster of the $\theta_n$-chip perform
a periodic computation. 
Thus, we say that the $\theta_n$-chip is an \textbf{aperiodic computer chip}.
Perhaps we can say it is an \emph{aperiodic monochip}, but we can not say it is
an \emph{aperiodic monotile} as in
\cite{MR4770585,MR4807152} because the same chip with
different inputs
has to be considered a distinct Wang tile.

\subsection*{Instances of the chip are metallic mean Wang tiles}
If we consider all possible values of inputs $u$ and $v$ in $V_n$
and if we restrict the outputs to be in the set $V_n$,
then we obtain a finite set of Wang tiles
\begin{equation}\label{eq:Ccal}
    \Ccal_n=
    \left\{
    \raisebox{-11.0mm}{
    \begin{tikzpicture}[auto,scale=.3]
    \def\t{.8}
    \def\dx{5}
    \def\dy{5}
    \def\x{0}
    \def\y{0}
            \draw[ultra thick,rounded corners,fill=gray!20] (\dx*\x,\dy*\y+0) rectangle (\dx*\x+4,\dy*\y+4);
            \draw[ultra thick,rounded corners] (\dx*\x+1,\dy*\y+1) rectangle (\dx*\x+3,\dy*\y+3);
    \node at (-1,2)  {$u$};
    \node at (2,-1)  {$v$};
    \node[right] at (4,2)   {$\theta_n(u,v)$};
    \node at (2,5)   {$\theta_n(v,u)$};
    \end{tikzpicture}
    }
    \middle|\,
        u,v\in V_n
        \text{ such that }
        \theta_n(u,v),\theta_n(v,u)\in V_n
    \right\}
\end{equation}
which is the finite set of all possible instances of the $\theta_n$-chip.

\newcommand\MainTheoremA{
    For every integer $n\geq1$,
    the Wang shift $\Omega_{\Ccal_n}$ defined by the $\theta_n$-chip is the
    $n^{th}$ metallic mean Wang shift $\Omega_n$.
}

\begin{maintheorem}\label{thm:is-the-metallic-mean-wang-shift}
    \MainTheoremA
\end{maintheorem}

Something unexpected and surprising happens in the proof of
Theorem~\ref{thm:is-the-metallic-mean-wang-shift}. The set $\Ccal_n$ of
instances of the $\theta_n$-chip is exactly equal to the extended set
$\Tcal_n'$ of metallic mean Wang tiles introduced in
\cite{labbe_metallic_I_2025} in order to 
prove the self-similarity of $\Omega_n$, see Proposition~\ref{prop:Ccal-equals-Tcaln'}.

\subsection*{Tile labels satisfy Equations}
The next result states that every tile in $\Ccal_n$ satisfy a system of equations.
While the equations associated with Kari's \cite{MR1417578} and Culik's
\cite{MR1417576} aperiodic set of Wang tiles are multiplicative, the ones
associated with $\Ccal_n$ are additive.

\newcommand\MainTheoremB{
    Let $n\geq1$ be an integer, $d=(0,-1,1)$ and $e=(1,0,0)$.
    The set of Wang tiles defined by the $\theta_n$-chip satisfy
    the following system of equations:
    \[
    \Ccal_n\subset
    \left\{
    \raisebox{-11.0mm}{
    \begin{tikzpicture}[auto,scale=.3]
    \def\t{.8}
    \def\dx{5}
    \def\dy{5}
    \def\x{0}
    \def\y{0}
            \draw[ultra thick,rounded corners,fill=gray!20] (\dx*\x,\dy*\y+0) rectangle (\dx*\x+4,\dy*\y+4);
            \draw[ultra thick,rounded corners] (\dx*\x+1,\dy*\y+1) rectangle (\dx*\x+3,\dy*\y+3);
    \node at (-1,2)  {$\ell$};
    \node at (2,-1)  {$b$};
    \node[right] at (4,2)   {$r$};
    \node at (2,5)   {$t$};
    \end{tikzpicture}
    }
    \in V_n\times V_n\times V_n\times V_n\,\;
    \middle|\,
    \arraycolsep=2pt
    \begin{array}{rcl}
        \langle \frac{1}{n}d, t+\ell\rangle - \langle e, \ell\rangle
             &=&\langle \frac{1}{n}d, b+r\rangle - \langle e, b\rangle\\
        \langle e, \ell  \rangle&=&\langle e, r  \rangle\\
        \langle e, b     \rangle&=&\langle e, t  \rangle
    \end{array}
    \right\}
\]
    where $\langle\_,\_\rangle$ denotes the canonical inner product of $\Z^3$.
}

\begin{maintheorem}\label{thm:equations-satisfied-by-tiles}
    \MainTheoremB
\end{maintheorem}

Equivalently, 
if we let
$\ell=(\ell_0,\ell_1,\ell_2)$,
$b=(b_0,b_1,b_2)$,
$r=(r_0,r_1,r_2)$ and
$t=(t_0,t_1,t_2)$,
the equations in the theorem say that
tiles in $\Ccal_n$ satisfy $\ell_0=r_0$, $b_0=t_0$ and
\begin{equation}
    \frac{t_2-t_1+\ell_2-\ell_1}{n}-\ell_0
    =
    \frac{b_2-b_1+r_2-r_1}{n}-b_0
\end{equation}
which reminds of Equation~\eqref{eq:wang-kari}.

Like Kari's and Culik's tiles,
these equations behave well with tilings and more equations
can be deduced for valid tilings of a rectangle, see Section~\ref{sec:equations}.
In particular, Equation~\eqref{eq:equations-for-a-rectangle-with-equal-sides}
says that in a tiling of a cylinder of height $k$,
    the average of the inner product with $\textstyle\frac{1}{n}d$ of the top
    labels of the cylinder is obtained 
    from the average of the inner product with $\textstyle\frac{1}{n}d$ of the bottom
    labels of the cylinder by $k$ rotations on the unit circle by a fixed angle.
    The angle is equal to the frequency of columns in the cylinder containing
    junction tiles and vertical strip colored tiles, which is a rational number.
    Therefore, the existence of a cyclic rectangle is not directly forbidden from 
    these equations. Note that we know from the self-similarity
    of $\Omega_n$ that the frequency of columns containing junction tile
    in every valid configuration in $\Omega_n$ is equal to $\beta^{-1}$,
    which is a irrational number \cite{labbe_metallic_I_2025}.

    It remains an open problem to deduce the aperiodicity of the Wang shift $\Omega_n$ 
    from the equations satisfied by the labels of $\theta_n$-chip
    as this is nicely done for Kari and Culik sets of tiles.
    See Section~\ref{sec:open-questions} for related open questions.

\subsection*{Existence of valid tilings}
Valid configurations in $\Omega_n$ can be constructed using the floor function on
linear forms.
Let $\Lambda_n:[0,1)^2\to\Z^3$ be defined as
\[
    \Lambda_n(x,y)
    =
    \left(
    \begin{array}{r}
        \lfloor               y-\beta^{-1}+1\rfloor\\
        \lfloor \beta^{-1}x + y-\beta^{-1}+1\rfloor\\
        \lfloor \beta     x + y-\beta^{-1}+1\rfloor
    \end{array}
    \right).
\]
where $\beta$ is the $n^{th}$ metallic mean, that is, the positive root of the
polynomial $x^2-nx-1$.
For every $(x,y)\in\R^2$, let
    \[
        \sctile_n(x,y)=
    \raisebox{-11mm}{
    \begin{tikzpicture}[auto,scale=1]
        \tile{white}{0}{0}{\Lambda_n(\{x\},\{y\})        }
                          {\Lambda_n(\{y\},\{x\})        }
                          {\Lambda_n(\{x{-}\beta^{-1}\},\{y\})}
                          {\Lambda_n(\{y{-}\beta^{-1}\},\{x\})}
    \end{tikzpicture}}
    \]
    be a Wang tile
    where $\{x\}=x-\lfloor x\rfloor$ is the fractional part of a number $x\in\R$.

\newcommand\MainTheoremC{
    For every integer $n\geq1$
        and every $(x,y)\in[0,1)^2$, the configuration
    \[
        \begin{array}{rccl}
            c_{(x,y)}:&\Z^2&\to&\Tcal_n\\
            &(i,j)&\mapsto &\sctile_n\left(x{+}i\beta^{-1},y{+}j\beta^{-1}\right)
        \end{array}
    \]
    is a valid tiling of the plane by the set of metallic mean Wang tiles $\Tcal_n$.
}

\begin{maintheorem}\label{thm:construct-valid-tilings}
    \MainTheoremC
\end{maintheorem}

This construction reminds of the proof of existence of tilings with
Kari and Culik tiles based on the balanced representation of real numbers
and first difference of Beatty sequences \cite{MR1417578,MR1417576}, see also
\cite{MR2369448,MR3668002}.

\subsection*{A factor map defined from averages of tile labels}
\begin{figure}[h]
\begin{center}
    \includegraphics{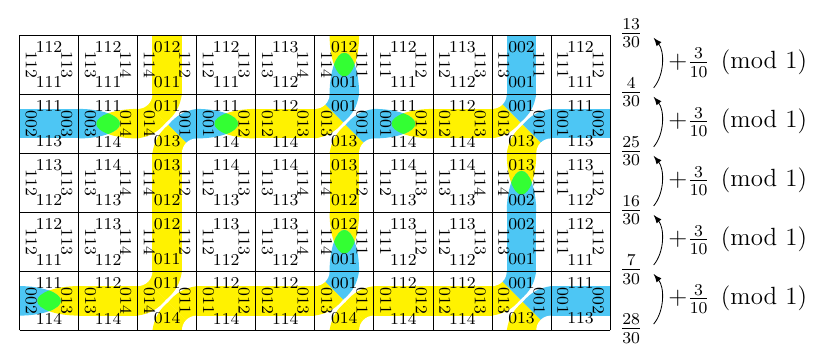}
\end{center}
    \caption{A $10\times5$ valid rectangular tiling with the set $\Tcal_n$ with $n=3$.
    The numbers indicated in the right margin are the average of the 
    inner products $\langle\frac{1}{n}d,v\rangle$ over the vectors $v$
    appearing as top (or bottom) labels of a horizontal row of tiles and where $d=(0,-1,1)$.
    We observe that these numbers increase by $\frac{3}{10}\pmod1$ from row to row.
    The number $\frac{3}{10}$ is equal to the frequency of columns containing
    junction tiles (a junction tile is a tile whose labels all start with 0).
    Observe that this is a cylindrical tiling (left and right outer labels of
    the rectangle match) which simplifies the equations involved because the
    left and right carries cancel.
    }
    \label{fig:T3-tiling-10x5-with-values}
\end{figure}

In Kari--Culik tilings \cite{MR1417578,MR1417576}, there is a well-defined notion
of average \cite{MR3606059} of the top tile labels along a bi-infinite horizontal row.
The change of value from one row to the next row is described by a piecewise
rationally multiplicative map.
In this context, metallic mean Wang shifts also behave like
Kari--Culik tilings.
It involves the consideration of the average of specific inner products and
irrational rotations instead of multiplications, see
Figure~\ref{fig:T3-tiling-10x5-with-values}
which can be compared with Figure~\ref{fig:Kari-tiling}.

We show that the average of the dot products of the vector
$\frac{1}{n}d=\frac{1}{n}(0,-1,1)$ with the top labels of a given row in a valid
configuration $\Z^2\to\Tcal_n$ in $\Omega_n$ is well-defined and takes a value in the interval
$[0,1]$ (see Equation~\eqref{eq:phi_n}).
By symmetry of the set $\Tcal_n$, the same holds
for the right labels of a given column. By considering the row and column going
through the origin of a configuration, the two averages define a map 
$\Phi_n:\Omega_n\to\torus^2$ (see Equation~\eqref{eq:Phi}).
We prove that this map is a factor map from the Wang shift to the 2-torus.

\newcommand\MainTheoremD{
    Let $d=(0,-1,1)$,
    $n\geq1$ be an integer
    and $\Omega_n$ be the $n^{th}$ metallic mean Wang shift.
    The map 
    \begin{equation}\label{eq:Phi_n}
    \begin{array}{rccl}
        \Phi_n:&\Omega_n & \to & \torus^2\\
        &w & \mapsto &
    \displaystyle
        \lim_{k\to\infty}\frac{1}{2k+1}\sum_{i=-k}^k
        \left(\begin{array}{c}
            \langle \textstyle\frac{1}{n} d,\scright(w_{0,i})\rangle\\[1mm]
        \langle \textstyle\frac{1}{n} d,\sctop(w_{i,0})\rangle
        \end{array}\right)
    \end{array}
    \end{equation}
    is a factor map, that is, it is continuous, onto and 
    commutes the shift $\dynsys{\Z^2}{\sigma}{\Omega_n}$
    with the toral $\Z^2$-rotation $\dynsys{\Z^2}{R_n}{\torus^2}$
    by the equation
    $
        \Phi_n\circ\sigma^k
        = R_n^k\circ \Phi_n 
    $
    for every $k\in\Z^2$
    where
    \[
    \begin{array}{rccl}
        R_n&:\Z^2\times\torus^2 & \to & \torus^2\\
        &(k,x) & \mapsto & R_n^k(x):=x + \beta k
    \end{array}
    \]
    and $\beta=\frac{n+\sqrt{n^2+4}}{2}$
    is the $n^{th}$ metallic mean, that is, the positive root of the polynomial $x^2-nx-1$.
}

\begin{maintheorem}\label{thm:factor-map}
    \MainTheoremD
\end{maintheorem}

As a consequence of Theorem~\ref{thm:factor-map}, we deduce that $\Omega_n$ is
aperiodic because $\beta$ is irrational and $R_n$ is a free $\Z^2$-action, see
Corollary~\ref{cor:Omegan_is_aperiodic}.
Note that since $\beta-\beta^{-1}=n$, we have $\beta=\beta^{-1}\pmod 1$.

Theorem~\ref{thm:factor-map} is an analogue of a result known for Kari and
Culik aperiodic Wang tilings which satisfy equations involving balanced
representations of real numbers and orbits of piecewise rationally
multiplicative maps, see also Theorem 16 in \cite{MR2369448} and Proposition 3
in \cite{MR3668002}. Here the result applies to all of the configurations in
the Wang shift $\Omega_n$.

\subsection*{A symbolic dynamical system and a Markov partition}
The Wang shift $\Omega_n$ can be independently described as a symbolic representation
of the dynamical system $\dynsys{\Z^2}{R_n}{\torus^2}$ by encoding its orbits with an
appropriate topological partition of $\torus^2$. The partition of
$\torus^2$ naturally emerges from the set of preimages of the map $\sctile_n$ and from
Theorem~\ref{thm:construct-valid-tilings}.

Since $\Lambda_n$ is defined as the floor of linear forms,
for every tile $t\in\Tcal_n$, the set
\[
    P_t=\Int{\sctile_n^{-1}(t)}
\]
is a polygonal open region in the unit square.
It satisfies that $\Pcal_n=\{P_t\mid t\in\Tcal_n\}$
is a topological partition of $\torus^2$ made of $(n+3)^2$ atoms.
The polygonal partition $\Pcal_n$ is the refinement
of two polygonal partitions 
$\east_n = \{\Lambda_n^{-1}(v)\colon v\in V_n\}$ and $\north_n$,
the second one being the image of the first under a
symmetry by the positive diagonal.
The partition $\east_n$
can be constructed by drawing
the following geodesics on the torus $\torus^2$:
\begin{itemize}
    \item two closed geodesics of slope $0$ and $\infty$ going through the origin $(0,0)$,
    \item a closed geodesic of slope 0 going through the point $(0,\beta^{-1})$,
    \item a geodesic of slope $-\beta^{-1}$ from $(0,\beta^{-1})$ to $(1,0)$,
    \item a geodesic of slope $-\beta$ from $(0,\beta^{-1})$ to $(1,0)$
        wrapping around the unit square fundamental domain $n$ times.
\end{itemize}
\begin{figure}[h]
\begin{center}
    \begin{tabular}{ccc}
        $\east_3$ & $\north_3$ & $\Pcal_3=\east_3\wedge\north_3$\\
 \includegraphics[scale=.85]{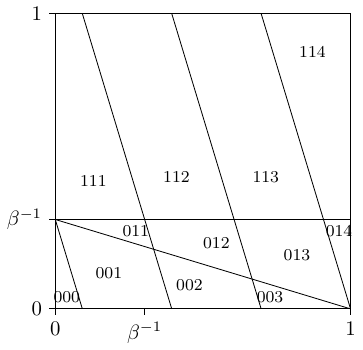}
&\includegraphics[scale=.85]{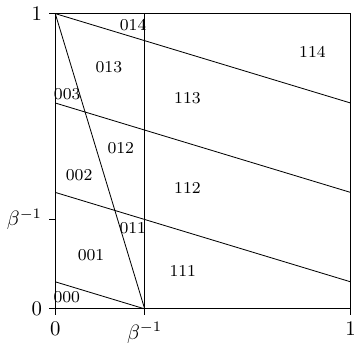}
&\includegraphics[scale=.85]{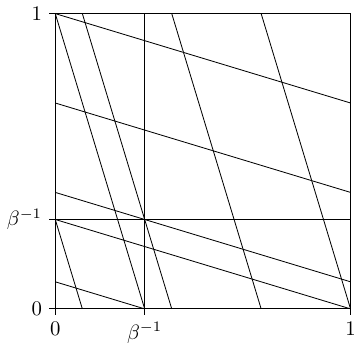}\\
    \end{tabular}
\caption{The partition $\east_3$ and its image $\north_3$ under a symmetry with
    the positive diagonal.  Their refinement is $\Pcal_3$ which is 
    a partition of the unit square into 36 polygonal atoms.
    Here $\beta$ is the third metallic mean, that is, the positive root of $x^2-3x-1$.}
\label{fig:partition-n3-intro}
\end{center}
\end{figure}
See an illustration of $\Pcal_n$ when $n=3$ in
Figure~\ref{fig:partition-n3-intro}.
Every open region defined by the complement of the geodesics can be identified with
a pair of vectors in $V_n$ and a unique tile in $\Tcal_n$ with such top and right labels.
As opposed to the four topological polygonal partitions associated with Jeandel-Rao
tilings \cite{MR4213162}, $\Pcal_n$ can be computed only from $\east_n$ and $\north_n$ 
without considering the $\south_n$ and $\west_n$ partitions. This is because
the set $\Tcal_n$ of tiles is NE-deterministic, see Theorem~\ref{thm:NE-SW-deterministic}.

The encoding of $\Z^2$-orbits under $R_n$ by the topological partition $\Pcal_n$
are 2-dimensional configurations whose
topological closure is the symbolic dynamical system $\Xcal_{\Pcal_n,R_n}$.
We prove that $\Xcal_{\Pcal_n,R_n}=\Omega_n$, and since $\Omega_n$ is a subshift of finite
type by definition, we have the following theorem.

\newcommand\MainTheoremE{
For every integer $n\geq1$, 
the symbolic dynamical system $\Xcal_{\Pcal_n,R_n}$
corresponding to $\Pcal_n,R_n$ is equal to the metallic mean Wang shift $\Omega_n$:
\[
    \Omega_n = \Xcal_{\Pcal_n,R_n}.
\]
In particular, $\Pcal_n$ is a Markov partition for the dynamical system
$\dynsys{\Z^2}{R_n}{\torus^2}$.
}

\begin{maintheorem}\label{thm:Markov-partition}
    \MainTheoremE
\end{maintheorem}

Markov partitions were originally defined for one-dimensional dynamical
systems $\dynsys{\Z}{T}{\torus^2}$
and were extended to $\Z^d$-actions by automorphisms of
compact Abelian group in \cite{MR1632169}.
Following \cite{MR4213162,MR4347332},
we use the same terminology 
and extend the definition proposed in \cite[\S 6.5]{MR1369092}
for dynamical systems defined by higher-dimensional actions by rotations,
see Definition~\ref{def:Markov}.

\subsection*{The maximal equicontinuous factor and an isomorphism}
Using Theorem~\ref{thm:Markov-partition}
and applying the results already proved for Jeandel--Rao Wang shift
\cite{MR4213162}, we have the following additional topological and measurable 
properties for the factor map.
We refer the reader to the preliminary Section~\ref{sec:prelim-wang-shifts}
for the notions and vocabulary on topological and measure-preserving dynamical systems
that are used in the statement.
A similar result holds for Penrose tilings \cite{MR1355301}.

\newcommand\MainTheoremF{
    The Wang shift $\Omega_n$ and the $\Z^2$-action $R_n$ have the following properties:
\begin{enumerate}[\rm (i)]
\item $\dynsys{\Z^2}{R_n}{\torus^2}$ is the maximal equicontinuous
    factor of $\dynsys{\Z^2}{\sigma}{\Omega_n}$,
\item the factor map $\Phi_n:\Omega_n\to\torus^2$ is almost one-to-one
    and its set of fiber cardinalities is $\{1,2,8\}$,
\item the shift-action $\dynsys{\Z^2}{\sigma}{\Omega_n}$ 
    on the metallic mean Wang shift is uniquely ergodic,
\item the measure-preserving dynamical system $(\Omega_n,\Z^2,\sigma,\nu)$
    is isomorphic to $(\torus^2,\Z^2,R_n,\lambda)$ where $\nu$ is the unique
        shift-invariant probability measure on $\Omega_n$ and $\lambda$ is the
        Haar measure on $\torus^2$.
\end{enumerate}
}

\begin{maintheorem}\label{thm:maximal-equicontinuous-factor}
    \MainTheoremF
\end{maintheorem}

\section{Preliminaries on dynamical systems, subshifts and Wang shifts}
\label{sec:prelim-wang-shifts}

This section follows the preliminary section of the chapter
\cite{labbe_three_2020} and article \cite{MR4213162}.

\subsection{Topological dynamical systems}

Most of the notions introduced here can be found in \cite{MR648108}.
A \defn{dynamical system} is
a triple $(X,G,T)$, where $X$ is a topological space, $G$ is a topological
group and $T$ is a continuous function $G\times X\to X$ defining a left action
of $G$ on $X$:
if $x\in X$, $e$ is the identity element of $G$ and $g,h\in G$, then using
additive notation for the operation in $G$ we have $T(e,x)=x$
and $T(g+h,x)=T(g,T(h,x))$.
In other words, if one denotes the transformation $x\mapsto T(g,x)$
by $T^g$, then $T^{g+h}=T^g T^h$.
In this work, we consider the Abelian group $G=\Z\times\Z$.

If $Y\subset X$, let $\overline{Y}$ denote the topological closure of $Y$ and
let $\overline{Y}^T:=\cup_{g\in G}T^g(Y)$ denote the $T$-closure of $Y$.
A subset $Y\subset X$ is \defn{$T$-invariant} if $\overline{Y}^T=Y$.
A dynamical system $(X,G,T)$ is called \defn{minimal} if $X$ does
not contain any nonempty, proper, closed $T$-invariant subset.
The left action of $G$ on $X$ is \defn{free}
if $g=e$ whenever there exists $x\in X$ such that $T^g(x)=x$.

Let $(X,G,T)$ and $(Y,G,S)$ be two dynamical systems with
the same topological group $G$.
A \defn{homomorphism} $\theta:(X,G,T)\to(Y,G,S)$ is a continuous
function $\theta:X\to Y$ satisfying the commuting property
that $S^g\circ\theta=\theta\circ T^g$ for every $g\in G$.
A homomorphism $\theta:(X,G,T)\to(Y,G,S)$ is called an \defn{embedding}
if it is one-to-one, a \defn{factor map} if it is onto, and a \defn{topological
conjugacy} if it is both one-to-one and onto and its inverse map is continuous.
If $\theta:(X,G,T)\to(Y,G,S)$ is a factor map,
then $(Y,G,S)$ is called a \defn{factor} of $(X,G,T)$
and $(X,G,T)$ is called an \defn{extension} of $(Y,G,S)$.
Two dynamical systems are \defn{topologically conjugate} if there is a
topological conjugacy between them.
A \defn{measure-preserving dynamical system} is defined as a system
$(X,G,T,\mu,\Bcal)$, where $\mu$ is a probability measure defined on
the Borel $\sigma$-algebra $\Bcal$ of subsets of $X$,
and $T^g:X\to X$ is a measurable map
which preserves the measure $\mu$ for all $g\in G$, that is,
$\mu(T^g(B))=\mu(B)$ for all $B\in\Bcal$. The measure $\mu$ is said to be
\defn{$T$-invariant}.
In what follows, 
when it is clear from the context,
we omit the Borel $\sigma$-algebra $\Bcal$ of subsets of $X$ and write $(X,G,T,\mu)$ 
to denote a measure-preserving dynamical system.

The set of all $T$-invariant probability measures of a dynamical
system $(X,G,T)$ is denoted by $\mathcal{M}^T(X)$.
A $T$-invariant probability measure on $X$ is called \defn{ergodic} if for every set
$B\in\Bcal$ such that $T^{g}(B)=B$ for all $g\in G$, we have that $B$ has either
zero or full measure. A
dynamical system $(X,G,T)$ is \defn{uniquely ergodic}
if it has only one invariant probability measure, i.e., $|\mathcal{M}^T(X)|=1$.
One can prove that a uniquely ergodic dynamical system is ergodic.
A dynamical system $(X,G,T)$ is said \defn{strictly ergodic}
if it is uniquely ergodic and minimal.

Let $(X,G,S,\mu,\Acal)$
and $(Y,G,T,\nu,\Bcal)$ be two measure-preserving dynamical systems.
We say that the two systems are
\defn{isomorphic} (mod 0) %
if there exist measurable sets $X_0\subset X$ and $Y_0\subset Y$
of full measure (i.e., $\mu(X_0)=1$ and $\nu(Y_0)=1$) with
$S^g(X_0)\subset X_0$, $T^g(Y_0)\subset Y_0$ for all $g\in G$
and there exists a bi-measurable bijection $\phi_0:X_0\to Y_0$,
\begin{itemize}
    \item which is measure-preserving, that is, $\mu(\phi_0^{-1}(B))=\nu(B)$ for all measurable sets 
        $B\subset Y_0$,
    \item satisfying $\phi_0\circ S^g(x)=T^g\circ\phi_0(x)$ for all $x\in X_0$ and $g\in G$.
\end{itemize}
The role of the set $X_0$ is to make precise the fact that the properties of
the isomorphism need to hold only on a set of full measure.
In this case, we call $\phi_0$ an \defn{isomorphism} (mod 0) with respect to $\mu$ and $\nu$. 
We also refer to an everywhere defined measurable map
$\phi:X\to Y$ as an \defn{isomorphism} (mod 0) with respect to $\mu$ and $\nu$ if 
$\phi(x)=\phi_0(x)$ with $x\in X$ for some $\phi_0$ and $X_0$ as above.
When $\phi$ is also a factor map, some authors say that
$\phi$ is a \defn{topo-isomorphism} in order to express both its topological
and measurable nature \cite{MR4425329}.

\subsection{Maximal equicontinuous factor}

A metrizable dynamical system $(X,G,T)$ is called \defn{equicontinuous} if
the family of homeomorphisms $\{T^g\}_{g\in G}$ is equicontinuous, i.e., if for
all $\varepsilon>0$ there exists $\delta>0$ such that
\[
    \dist(T^g(x), T^g(y)) < \varepsilon
\]
for all $g\in G$ and all $x,y\in X$ with $\dist(x,y)<\delta$.
According to a well-known theorem~\cite[Theorem 3.2]{MR3381481},
equicontinuous
minimal systems defined by the action of an Abelian group
are rotations on groups.

We say that $\theta:(X,G,T)\to(Y,G,S)$ is an \defn{equicontinuous factor} if
$\theta$ is a factor map and $(Y,G,S)$ is equicontinuous.
We say that $(X_{\rm max}, G, T_{\rm max})$ is the \defn{maximal equicontinuous
factor} of $(X,G,T)$ if 
there exists an equicontinuous factor
$\pi_{\rm max}:(X,G,T)\to(X_{\rm max}, G, T_{\rm max})$,
such that for any
equicontinuous factor $\theta:(X,G,T)\to(Y,G,S)$,
there exists a unique factor map $\psi:(X_{\rm max}, G, T_{\rm max})\to(Y,G,S)$
with $\psi\circ\pi_{\rm max}=\theta$.
The maximal equicontinuous factor exists and is unique (up to topological
conjugacy), see \cite[Theorem 3.8]{MR3381481} and \cite[Theorem 2.44]{MR2041676}.

Let $\theta:(X,G,T)\to(Y,G,S)$ be a factor map.
We call the preimage set $\theta^{-1}(y)$ of a point $y\in Y$ the \defn{fiber} of $\theta$ over $y$.
The cardinality of the fiber $\theta^{-1}(y)$ for some $y\in Y$ has an important role and is related to the definition of other notions, see \cite{MR3381481}.
In particular, the factor map $\theta$ is \defn{almost one-to-one} if 
$\{y\in Y:\card(\theta^{-1}(y))=1\}$
is a $G_\delta$-dense set in $Y$ (that is a countable intersection of open sets
which is dense in $Y$).
In that case, $(X,G,T)$ is an \defn{almost one-to-one extension} of $(Y,G,S)$.
The \defn{set of fiber cardinalities}
of a factor map $\theta:(X,G,T)\to(Y,G,S)$
is the set $\{\card(\theta^{-1}(y)) : y \in Y\}\subset\N\cup\{\infty\}$, see \cite{MR1877329}.
The set of fiber cardinalities of the maximal equicontinuous factor of a
minimal dynamical system is invariant under topological conjugacy,
see for instance \cite[Lemma~2.2]{MR4213162}.

\subsection{Subshifts and shifts of finite type}

In this section, we introduce multidimensional subshifts,
a particular type of dynamical systems 
\cite[\S 13.10]{MR1369092},
\cite{MR1861953,MR2078846,MR3525488}.
Let $\Acal$ be a finite set, $d\geq 1$, and let $\Acal^{\Z^d}$ be the set of all maps
$x:\Z^d\to\Acal$, equipped with the compact product topology. 
An element $x\in\Acal^{\Z^d}$ is called \defn{configuration}
and we write it as $x=(x_\bm)=(x_\bm:\bm\in\Z^d)$,
where $x_\bm\in\Acal$ denotes the value of $x$ at $\bm$. 
The topology on $\Acal^{\Z^d}$ is compatible with the metric defined for all
configurations $x,x'\in\Acal^{\Z^d}$ by $\dist(x,x')=2^{-\min\left\{\Vert\bn\Vert\,:\,
x_\bn\neq x'_\bn\right\}}$
where $\Vert\bn\Vert = |n_1| + \dots + |n_d|$.
The \defn{shift action} $\sigma:\bn\mapsto
\sigma^\bn$ of the additive group $\Z^d$ on $\Acal^{\Z^d}$ is defined by
\begin{equation}\label{eq:shift-action}
    (\sigma^\bn(x))_\bm = x_{\bm+\bn}
\end{equation}
for every $x=(x_\bm)\in\Acal^{\Z^d}$ and $\bn\in\Z^d$. 
If $X\subset \Acal^{\Z^d}$,
let $\overline{X}$ denote the topological closure of $X$
and let $\shiftclosure{X}:=\{\sigma^\bn(x)\mid x\in X, \bn\in\Z^d\}$
denote the shift-closure of $X$.
A subset $X\subset
\Acal^{\Z^d}$ is \defn{shift-invariant} if 
$\shiftclosure{X}=X$. A closed, shift-invariant subset
$X\subset\Acal^{\Z^d}$ is a \defn{subshift}. 
If $X\subset\Acal^{\Z^d}$ is a subshift we write
$\sigma=\sigma^X$ for the restriction of the shift action
\eqref{eq:shift-action} to $X$. 
When $X$ is a subshift,
the triple $(X,\Z^d,\sigma)$ is a dynamical system
and the notions presented in the previous section hold.

A configuration $x\in X$ is \defn{periodic} if there is a nonzero vector
$\bn\in\Z^d\setminus\{\zero\}$ such that $x=\sigma^\bn(x)$
and otherwise it is \defn{nonperiodic}.
We say that a nonempty subshift $X$ is \defn{aperiodic}
if the shift action $\sigma$ on $X$ is free.

For any subset $S\subset\Z^d$ let $\pi_S:\Acal^{\Z^d}\to\Acal^S$ denote the
projection map which restricts every $x\in\Acal^{\Z^d}$ to $S$. 
A \defn{pattern} is a function $p\in\Acal^S$ for some finite subset
$S\subset\Z^d$.
To every pattern $p\in\Acal^S$ corresponds
a subset $\pi_S^{-1}(p)\subset\Acal^{\Z^d}$ called \defn{cylinder}.
A nonempty set $X\subset\Acal^{\Z^d}$ is a
\defn{subshift} if and only if there exists a set $\Fcal$
of \defn{forbidden} patterns such that
\begin{equation}\label{eq:SFT}
    X = \{x\in\Acal^{\Z^d} \mid \pi_S\circ\sigma^\bn(x)\notin\Fcal
    \text{ for every } \bn\in\Z^d \text{ and } S\subset\Z^d\},
\end{equation}
see \cite[Prop.~9.2.4]{MR3525488}.
A subshift $X\subset\Acal^{\Z^d}$ is a 
\defn{subshift of finite type} (SFT) if there exists a finite set $\Fcal$ such that \eqref{eq:SFT} holds.
In this article, we consider shifts of finite type on $\Z\times\Z$, that is, the case
$d=2$.

\subsection{Wang shifts}

A \defn{Wang tile} 
is a tuple of four colors $(a,b,c,d)\in I\times J\times
I\times J$
where $I$
is a finite set of vertical colors
and $J$
is a finite set of horizontal colors, see
\cite{wang_proving_1961,MR0297572}.
A Wang tile is represented as a unit square with colored edges:
\begin{center}
    \raisebox{-9.5mm}{
    \begin{tikzpicture}[auto]
    \tile{white}{0}{0}{a}{b}{c}{d}
    \end{tikzpicture}}
\end{center}
For each Wang tile $\tau=(a,b,c,d)$, let
$\scright(\tau)=a$,
$\sctop(\tau)=b$,
$\scleft(\tau)=c$,
$\scbottom(\tau)=d$
denote respectively the colors of the right, top, left and bottom edges of $\tau$.

\begin{figure}[h]
\begin{center}
    \includegraphics{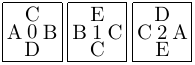}
\end{center}
    \caption{The set of 3 Wang tiles introduced
    in \cite{wang_proving_1961} using letters $\{A,B,C,D,E\}$ instead of
    numbers from the set $\{1,2,3,4,5\}$ for labeling the edges.
    Each tile is identified uniquely by an index from the
    set $\{0,1,2\}$ written at the center each tile.}
    \label{fig:wang-three-tiles}
\end{figure}

\begin{figure}[h]
\begin{center}
\begin{tikzpicture}
    \node (A) at (0,0) {$\arraycolsep=1.8pt
                        \input{SAGEOUTPUT/Wang_1961_three_tiles_tiling_config.tex}
                        $};
    \node (B) at (3,0) {\includegraphics{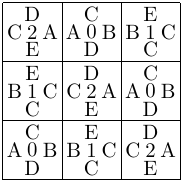}};
    \draw[|->] (A) to (B);
\end{tikzpicture}
\end{center}
    \caption{A finite $3\times 3$ pattern on the left is valid with respect to the Wang tiles
    since it respects Equations~\eqref{eq:validwangtiling1}
    and~\eqref{eq:validwangtiling2}. Validity can be verified on the tiling shown on
    the right.}
    \label{fig:wang-three-tiles-3x3-tiling}
\end{figure}

Let $\Tcal=\{t_0,\dots,t_{m-1}\}$ be a set of Wang tiles as the one shown in Figure~\ref{fig:wang-three-tiles}.
A configuration $x:\Z^2\to\{0,\dots,m-1\}$ is \defn{valid} with respect to $\Tcal$ if
it assigns a tile in $\Tcal$ to each position of $\Z^2$ so that contiguous edges
of adjacent tiles have the same color, that is,
\begin{align}
    \scright(t_{x(\bn)})&=\scleft(t_{x(\bn+\be_1)})\label{eq:validwangtiling1}\\
    \sctop(t_{x(\bn)})&=\scbottom(t_{x(\bn+\be_2)})\label{eq:validwangtiling2}
\end{align}
for every $\bn\in\Z^2$ where $\be_1=(1,0)$ and $\be_2=(0,1)$.
A finite pattern which is valid with respect to $\Ucal$ is shown in 
Figure~\ref{fig:wang-three-tiles-3x3-tiling}.

Let $\Omega_\Tcal\subset\{0,\dots,m-1\}^{\Z^2}$ denote the set of all valid 
configurations with respect to $\Tcal$.
Together with the shift action $\sigma$ of $\Z^2$,
$\Omega_\Tcal$ is a subshift that we call a \defn{Wang shift}.
Furthermore, $\Omega_\Tcal$ is a subshift of finite type (SFT) of the form \eqref{eq:SFT}
since $\Omega_\Tcal$ is the subshift defined from the finite set of
forbidden patterns made of all horizontal and vertical dominoes of two tiles
that do not share an edge of the same color.
Reciprocally, every subshift of finite type can be encoded into a Wang shift
following a well-known construction (see \cite[p.~141-142]{MR1014984}).

To a configuration $x\in\Omega_\Tcal$ corresponds a tiling of the plane $\R^2$ by
the tiles $\Tcal$ where the unit square Wang tile $t_{x(\bn)}$ is placed at position $\bn$ for every
$\bn\in\Z^2$, as in Figure~\ref{fig:wang-three-tiles-3x3-tiling}.
In this article, we consider tilings from the symbolic point of view.
In particular, we represent tilings of plane by Wang tiles symbolically 
by configurations $\Z^2\to\Tcal$.

A configuration $x\in\Omega_\Tcal$ is \defn{periodic} if there exists
$\bn\in\Z^2\setminus\{0\}$ such that $x=\sigma^\bn(x)$.
A set of Wang tiles $\Tcal$ is \defn{periodic} if there exists a periodic configuration
$x\in\Omega_\Tcal$. 
Originally, Wang thought that every set of Wang tiles $\Tcal$ is periodic 
as soon as $\Omega_\Tcal$ is nonempty \cite{wang_proving_1961}.
This statement is equivalent to the existence of an algorithm 
solving the \emph{domino problem}, that is, taking as input a set of Wang tiles
and returning \textit{yes} or \textit{no} whether there exists a valid
configuration with these tiles. 
Berger, a student of Wang, later proved that the domino problem is undecidable
and he also provided a first example of an aperiodic set of Wang tiles
\cite{MR0216954}.
A set of Wang tiles $\Tcal$ is \defn{aperiodic} if
the Wang shift $\Omega_\Tcal$ is a nonempty aperiodic subshift.
This means that in general one can not decide the emptiness of a Wang shift
$\Omega_\Tcal$. 

\section{The family of metallic mean Wang tiles}
\label{sec:family-metallic}

In this section, we recall from \cite{labbe_metallic_I_2025}
the definition of 
the set $\Tcal_n$ of metallic mean Wang tiles
and the extended set $\Tcal_n'$ which satisfies $\Tcal_n\subset\Tcal_n'$.
The extended set $\Tcal_n'$ was used to prove the self-similarity of the Wang
shift $\Omega_n$ defined over $\Tcal_n$.

For every integer $n\in\Z$, 
we write $\overline{n}$ to denote $n+1$
and $\underline{n}$ to denote $n-1$:
\begin{align*}
    \overline{n} &:= n+1,\\
    \underline{n} &:= n-1.
\end{align*}

For every Wang tile $\tau=(a,b,c,d)$, we define its
symmetric image under a symmetry by the positive diagonal as
 $\widehat{\tau}=(b,a,d,c)$:
\[
    \text{ if }
    \tau=
    \raisebox{-9.8mm}{
    \begin{tikzpicture}[auto]
    \tile{white}{0}{0}{a}{b}{c}{d}
    \end{tikzpicture}},
    \qquad
    \text{ then }
    \qquad
    \widehat{\tau}=
    \raisebox{-8.0mm}{
    \begin{tikzpicture}[auto]
    \tile{white}{0}{0}{b}{a}{d}{c}
    \end{tikzpicture}}.
\]

\subsection{The tiles}

For every integer $n\geq 1$, let
\[
    V_n = \{(v_0,v_1,v_2)\in\Z^3\colon 0\leq v_0\leq v_1\leq 1\text{ and }v_1\leq v_2\leq n+1\}.
\]
be a set of vectors having non-decreasing entries with an upper bound of 1 on the middle entry
and an upper bound of $n+1$ on the last entry.
The label of the edges of the Wang tiles considered in this article are vectors
in $V_n$. To lighten the figures and the presentation of the Wang tiles, it is convenient to denote
the vector $(v_0,v_1,v_2)\in V_n$ more compactly as a word $v_0v_1v_2$.
For instance the vector $(1,1,1)$ is represented as $111$.

To help the reading of the tiles and tilings, we assign a color to the vectors
according to the following rule:
a vector $v\in 00\N$ is drawn in blue,
a vector $v\in 01\N$ is drawn in yellow and
a vector $v\in 11\N$ is drawn in white.
Overlap between blue and yellow region will be shown in green.

For every integer $n\geq1$ and
for every $i,j\in\N$ such that $0\leq i\leq n$ and $0\leq j\leq n$,
we have the following white tiles:
\[
\begin{array}{|c|}
    \hline
    \text{\textbf{white} tiles}\\
    \hline
    w_n^{i,j} = 
    \raisebox{-10mm}{
    \begin{tikzpicture}[auto]
        \tikzstyle{every node}=[font=\footnotesize]
        \tile{white}{0}{0}{11\ibar}{11\jbar}{11i}{11j}
    \end{tikzpicture}}\\
    \hline
\end{array}
\]
For every $i,n\in\N$ such that $0\leq i\leq n$, 
we have the following blue, yellow, green and antigreen tiles:
\[
\begingroup
\allowdisplaybreaks
\arraycolsep=5pt
\begin{array}{|c|c|c|}
    \hline
    &
    \text{horizontal tiles}&
    \text{vertical tiles}\\
    \hline
    \text{\textbf{blue} tiles}&
b_n^i = 
    \raisebox{-9mm}{
    \begin{tikzpicture}[auto]
        \tikzstyle{every node}=[font=\footnotesize]
    \tileH{\ourColorBlue}{0}{0}{00\ibar}{111}{00i}{11n}
    \end{tikzpicture}}\qquad
&
\widehat{b_n^i} = 
    \raisebox{-9mm}{
    \begin{tikzpicture}[auto]
        \tikzstyle{every node}=[font=\footnotesize]
    \tileV{\ourColorBlue}{5}{0}{111}{00\ibar}{11n}{00i}
    \end{tikzpicture}}\\
    \hline
    \text{\textbf{yellow} tiles}&
y_n^i = 
    \raisebox{-9mm}{
    \begin{tikzpicture}[auto]
        \tikzstyle{every node}=[font=\footnotesize]
    \tileH{\ourColorYellow}{0}{0}{01\ibar}{112}{01i}{11\nbar}
    \end{tikzpicture}}
&
\widehat{y_n^i} = 
    \raisebox{-9mm}{
    \begin{tikzpicture}[auto]
        \tikzstyle{every node}=[font=\footnotesize]
    \tileV{\ourColorYellow}{5}{0}{112}{01\ibar}{11\nbar}{01i}
    \end{tikzpicture}}\\
\hline
    \text{\textbf{green} overlap tiles}&
g_n^i = 
    \raisebox{-9mm}{
    \begin{tikzpicture}[auto]
        \tikzstyle{every node}=[font=\footnotesize]
        \tileHgreen{0}{0}{01\ibar}{111}{00i}{11\nbar}
    \end{tikzpicture}}
&
\widehat{g_n^i} = 
    \raisebox{-9mm}{
    \begin{tikzpicture}[auto]
        \tikzstyle{every node}=[font=\footnotesize]
    \tileVgreen{5}{0}{111}{01\ibar}{11\nbar}{00i}
    \end{tikzpicture}}\\
\hline
    \text{\textbf{antigreen} no overlap tiles}&
a_n^i = 
    \raisebox{-9mm}{
    \begin{tikzpicture}[auto]
        \tikzstyle{every node}=[font=\footnotesize]
    \tileHantigreen{0}{0}{00\ibar}{112}{01i}{11n}
    \end{tikzpicture}}
&
\widehat{a_n^i} = 
    \raisebox{-9mm}{
    \begin{tikzpicture}[auto]
        \tikzstyle{every node}=[font=\footnotesize]
    \tileVantigreen{5}{0}{112}{00\ibar}{11n}{01i}
    \end{tikzpicture}}\\
\hline
\end{array}
\endgroup
\]

For every $n\in\N$ and $k,\ell,r,s\in\{0,1\}$ such that $k\leq\ell$ and $r\leq s$,
we have the following junction tiles (the gray region will be drawn in a blue
or yellow color depending on the specific values of $k,\ell,r,s$ according
to the same rule as above):
\[
\begingroup
\allowdisplaybreaks
\arraycolsep=20pt
\begin{array}{|c|}
    \hline
    \text{\textbf{junction} tiles}\\
    \hline
j_n^{k,\ell,r,s}=
\raisebox{-10mm}{
\begin{tikzpicture}[auto]
\tikzstyle{every node}=[font=\footnotesize]
    \tileJunctionGRAY{0}{2}{(0,k,\ell)}{(0,r,s)}{(0,s,r+n)}{(0,\ell,k+n)}{gray!50}
\end{tikzpicture}}\\
    \hline
\end{array}
\endgroup
\]
Junction tiles play a similar role as junction tiles in \cite{MR1014984}.

\subsection{The extended set $\Tcal_n'$ of metallic mean Wang tiles}
In this section, we give the
definition of the family of extended sets of Wang tiles 
$(\Tcal_n')_{n\geq1}$.

From the above, we define the following sets of tiles:
\begin{align*}
W_n &= \left\{ w_n^{i,j} \mid 1\leq i\leq n, 1\leq j\leq n \right\}
    &&(n^2\text{ white tiles}),\\
B'_n &= \left\{ b_n^i \mid 0\leq i \leq n \right\}
    &&(n+1\text{ horizontal blue tiles}),\\
Y_n &= \left\{ y_n^i \mid 1\leq i\leq n \right\}
    &&(n\text{ horizontal yellow tiles}),\\
G_n &= \left\{ g_n^i\mid 0\leq i\leq n \right\}
    &&(n+1\text{ horizontal green tiles}),\\
A_n &= \left\{ a_n^i \mid 1\leq i \leq n \right\}
    &&(n \text{ horizontal antigreen tiles}).
\end{align*}

Finally, we have a set of 9 junction tiles:
\begingroup
\allowdisplaybreaks
\begin{align*}
J_n' &=\left\{ 
    j_n^{0,0,0,0},
    j_n^{0,0,0,1},
    j_n^{0,0,1,1},
    j_n^{0,1,0,0},
    j_n^{0,1,0,1},
    j_n^{0,1,1,1},
    j_n^{1,1,0,0}, 
    j_n^{1,1,0,1},
    j_n^{1,1,1,1}
    \right\}\\
     &= 
    \left\{
    \raisebox{-9mm}{
    \begin{tikzpicture}[auto]
        \def\dt{2.2}
        \tikzstyle{every node}=[font=\footnotesize]
        \tileJunctionIIII{2*\dt}{0}  {}{011}{01\nbar}{}
        \draw[fill=white,draw=white] (2*\dt,0) -- ++ (1,1) -- ++ (0,-1) -- cycle;
        \draw[dashed,thick]  (2*\dt,0) -- ++ (1,1);
        \tileJunctionOOII{\dt}{0}  {}{001}{01n}{}
        \draw[fill=white,draw=white] (\dt,0) -- ++ (1,1) -- ++ (0,-1) -- cycle;
        \draw[dashed,thick]  (\dt,0) -- ++ (1,1);
        \tileJunctionOOOO{0}{0}  {}{000}{00n}{}        
        \draw[fill=white,draw=white] (0,0) -- ++ (1,1) -- ++ (0,-1) -- cycle;
        \draw[dashed,thick]  (0,0) -- ++ (1,1);
    \end{tikzpicture}}
    \right\}
    \times
    \left\{
    \raisebox{-9mm}{
    \begin{tikzpicture}[auto]
        \def\dt{2.2}
        \tikzstyle{every node}=[font=\footnotesize]
        \tileJunctionIIII{2*\dt}{0}  {011}{}{}{01\nbar}
        \draw[fill=white,draw=white] (2*\dt,0) -- ++ (1,1) -- ++ (-1,0) -- cycle;
        \draw[dashed,thick]  (2*\dt,0) -- ++ (1,1);
        \tileJunctionOOII{\dt}{0}{001}{}{}{01n}
        \draw[fill=white,draw=white] (\dt,0) -- ++ (1,1) -- ++ (-1,0) -- cycle;
        \draw[dashed,thick]  (\dt,0) -- ++ (1,1);
        \tileJunctionOOOO{0}{0}  {000}{}{}{00n}        
        \draw[fill=white,draw=white] (0,0) -- ++ (1,1) -- ++ (-1,0) -- cycle;
        \draw[dashed,thick]  (0,0) -- ++ (1,1);
    \end{tikzpicture}}
    \right\}\\
    &= 
    \left\{
    \raisebox{-31mm}{
    \begin{tikzpicture}[auto]
        \tikzstyle{every node}=[font=\footnotesize]
        \tileJunctionIOOI{0}{4.4}{011}{000}{00n}{01\nbar}
        \tileJunctionIOII{3}{4.4}{011}{001}{01n}{01\nbar}
        \tileJunctionIIII{6}{4.4}{011}{011}{01\nbar}{01\nbar}
        \tileJunctionOOOI{0}{2.2}{001}{000}{00n}{01n}        
        \tileJunctionOOII{3}{2.2}{001}{001}{01n}{01n}
        \tileJunctionOIII{6}{2.2}{001}{011}{01\nbar}{01n}
        \tileJunctionOOOO{0}{0}  {000}{000}{00n}{00n}        
        \tileJunctionOOIO{3}{0}  {000}{001}{01n}{00n}
        \tileJunctionOIIO{6}{0}  {000}{011}{01\nbar}{00n}
    \end{tikzpicture}}
    \right\}
    \qquad\qquad(9 \text{ junction tiles}).
\end{align*}
\endgroup

We may observe that $\widehat{W_n}=W_n$
and $\widehat{J'_n}=J'_n$ are closed under reflection.
Also, 
$\widehat{B'_n}$ are $n+1$ vertical blue tiles,
$\widehat{Y_n}$ are $n$ vertical yellow tiles,
$\widehat{G_n}$ are $n+1$ vertical green tiles and
$\widehat{A_n}$ are $n$ vertical antigreen tiles.

The extended set of metallic mean Wang tiles $\Tcal_n'$ can be described in
terms of the white, yellow, green, blue, antigreen and junction tiles seen
before.

\begin{figure}[h]
\begin{center}
\includegraphics{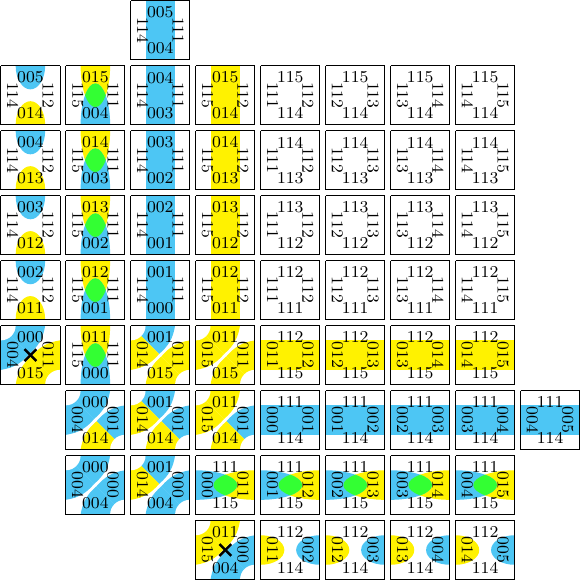}
\end{center}
    \caption{Extended metallic mean Wang tile sets $\Tcal_n'$ for $n=4$.
    The junction tiles 
    $j_n^{0,0,1,1}$ and $j_n^{1,1,0,0}$
    are shown with a $\times$-mark
    in their center.
    }
    \label{fig:T'n-for-1-a-5}
\end{figure}

\begin{definition}[Extended set of metallic mean Wang tiles
    \cite{labbe_metallic_I_2025}]\label{def:Tcaln'} 
Let
\[
\Tcal_n'= W_n\cup Y_n \cup \widehat{Y_n} 
             \cup G_n \cup \widehat{G_n} 
             \cup B'_n \cup \widehat{B'_n} 
             \cup A_n \cup \widehat{A_n} 
             \cup J'_n.
\]
The set $\Tcal_n'$ defines the \defn{extended metallic mean Wang shift}
$\Omega'_n=\Omega_{\Tcal_n'}$.
\end{definition}

The set $\Tcal_n'$
contains $n^2+2(n+1+n+n+1+n)+9=n^2+8n+13$ Wang tiles.
The set of Wang tiles $\Tcal_n'$ for $n=4$ is shown in
Figure~\ref{fig:T'n-for-1-a-5}.

\subsection{The family $\Tcal_n$ of $(n+3)^2$ Wang tiles}

\begin{figure}[h]
\begin{center}
    \begin{tabular}{ccc}
         \includegraphics{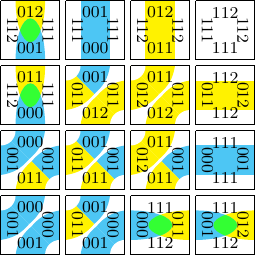}
         &\includegraphics{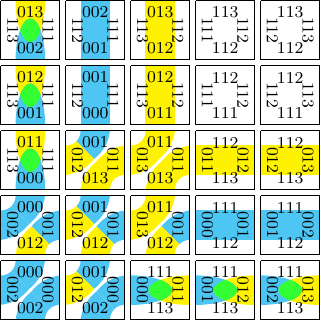}
         &\includegraphics{SAGEOUTPUT/W3_tiles.pdf}\\
        $\Tcal_1$ & $\Tcal_2$ & $\Tcal_3$\\
    \end{tabular}
    \begin{tabular}{ccc}
         \includegraphics{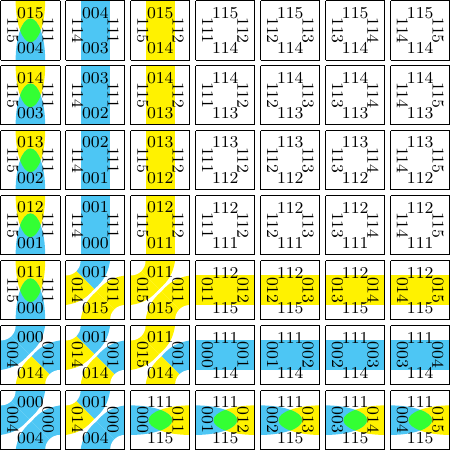}
         &\includegraphics{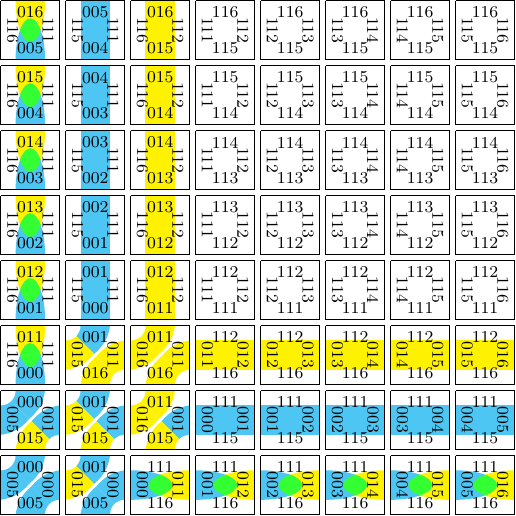}\\
        $\Tcal_4$ & $\Tcal_5$
    \end{tabular}
\end{center}
    \caption{Metallic mean Wang tile sets $\Tcal_n$ for $n=1,2,3,4,5$.}
    \label{fig:Tn-for-1-a-5}
\end{figure}

In this section, we give the
definition of the family of sets of Wang tiles 
$(\Tcal_n)_{n\geq1}$.
The set $\Tcal_n$ is a subset of $\Tcal_n'$ defined as follows.
Let
\begin{align*}
    B_n &= B_n' \setminus \left\{ b_n^n\right\}
    &&(\text{subset of $n$ horizontal blue tiles}),\\
    J_n &= J_n' \setminus \left\{ j_n^{1,1,0,0}, j_n^{0,0,1,1} \right\}
    &&(\text{subset of 7 junction tiles}).
\end{align*}

\begin{definition}[Metallic mean Wang tiles\cite{labbe_metallic_I_2025}]
For every positive integer $n$,
we construct the set of Wang tiles
\[
    \Tcal_n= W_n\cup Y_n \cup \widehat{Y_n} 
                \cup G_n \cup \widehat{G_n} 
                \cup B_n \cup \widehat{B_n} 
                \cup J_n.
\]
The set of tiles defines the \defn{Metallic mean Wang shift}
$\Omega_n=\Omega_{\Tcal_n}$.
\end{definition}
The subset $\Tcal_n$ contains $n^2+2(n+n+1+n)+7=(n+3)^2$ Wang tiles.
They are shown in Figure~\ref{fig:Tn-for-1-a-5} for $n=1,2,3,4,5$.

\section{The $\theta_n$-chip and metallic mean Wang tiles}
\label{sec:thetan-chip}

In this section, we relate the
$\theta_n$-chip with metallic mean Wang tiles.
The proposition below provides an independent
characterization of the extended set $\Tcal_n'$ of metallic-mean Wang tiles 
as instances of the $\theta_n$-chip, see Equation~\ref{eq:Ccal}.

\begin{proposition}\label{prop:Ccal-equals-Tcaln'}
    For every $n\geq1$,
    the set of instances of the computer chip
    is equal to the extended set of metallic mean Wang tiles, that is,
    $\Ccal_n=\Tcal_n'$.
\end{proposition}

\begin{proof}
    ($\subseteq$)
    Let
    $\tau=\raisebox{-9.5mm}{
    \begin{tikzpicture}[auto]
        \tile{white}{0}{0}{\theta_n(u,v)}{\theta_n(v,u)}{u}{v}
    \end{tikzpicture}
    }$ be a Wang tile such that $u=(u_0,u_1,u_2)\in V_n$,
    $v=(v_0,v_1,v_2)\in V_n$, $\theta_n(u,v)\in V_n$ and $\theta_n(v,u)\in V_n$.
    We proceed case by case:
    \begin{itemize}
        \item If $u_0=1$ and $v_0=1$, then $1=u_1\leq u_2$, $1=v_1\leq v_2$ and
            \begin{align*}
                \theta_n(u,v)&=(u_0,1,u_2+1)=(1,1,u_2+1)\in V_n,\\
                \theta_n(v,u)&=(v_0,1,v_2+1)=(1,1,v_2+1)\in V_n.
            \end{align*}
            Thus, $0\leq u_2\leq n$ and $0\leq v_2\leq n$ and $\tau\in W_n$
            is a white tile.
        \item If $u_0=0$ and $v_0=1$, then 
            \begin{align*}
                \theta_n(u,v)&=(u_0,v_2-n,u_2+1)=(0,v_2-n,u_2+1)\in V_n,\\
                \theta_n(v,u)&=(v_0,1,u_1+v_0)=(1,1,u_1+1)\in V_n,
            \end{align*}
            where $0\leq u_2\leq n$, $n\leq v_2\leq n+1$ and $0\leq u_1\leq 1$.
            There are four possibilities according to the values
            of $v_2\in\{n,n+1\}$ and $u_1\in\{0,1\}$ that we consider case by case:
            \begin{itemize}
                \item If $v_2=n$ and $u_1=0$, then
                $\tau
                =\raisebox{-9.5mm}{
                \begin{tikzpicture}[auto]
                    \tile{white}{0}{0}{(0,0,u_2+1)}{(1,1,1)}{(0,0,u_2)}{(1,1,n)}
                \end{tikzpicture}}
                = b_n^{u_2}\in B_n \cup \{b_n^n\}
                $ is a blue horizontal stripe tile with $0\leq u_2\leq n$.
                \item If $v_2=n$ and $u_1=1$, then
                $\tau
                =\raisebox{-9.5mm}{
                \begin{tikzpicture}[auto]
                    \tile{white}{0}{0}{(0,0,u_2+1)}{(1,1,2)}{(0,1,u_2)}{(1,1,n)}
                \end{tikzpicture}}
                = a_n^{u_2}\in A_n
                $ is an antigreen horizontal tile with $1\leq u_2\leq n$.
                \item If $v_2=n+1$ and $u_1=0$, then
                $\tau
                =\raisebox{-9.5mm}{
                \begin{tikzpicture}[auto]
                    \tile{white}{0}{0}{(0,1,u_2+1)}{(1,1,1)}{(0,0,u_2)}{(1,1,n+1)}
                \end{tikzpicture}}
                = g_n^{u_2}\in G_n
                $ is a green horizontal overlap tile with $0\leq u_2\leq n$.
                \item If $v_2=n+1$ and $u_1=1$, then
                $\tau
                =\raisebox{-9.5mm}{
                \begin{tikzpicture}[auto]
                    \tile{white}{0}{0}{(0,1,u_2+1)}{(1,1,2)}{(0,1,u_2)}{(1,1,n+1)}
                \end{tikzpicture}}
                = y_n^{u_2}\in Y_n
                $ is a yellow horizontal stripe tile with $1\leq u_2\leq n$.
            \end{itemize}
        \item If $u_0=1$ and $v_0=0$, the possibilities are the symmetric image
            of the previous case. Thus, 
            $\tau\in \widehat{B_n}\cup \{\widehat{b_n^n}\} \cup
                     \widehat{A_n}\cup 
                     \widehat{G_n}\cup 
                     \widehat{Y_n}$ is a blue, antigreen, green or yellow
                     vertical tile.
        \item If $u_0=0$ and $v_0=0$, then
            \begin{align*}
                \theta_n(u,v)&=(u_0,v_2-n,v_1+u_0)=(0,v_2-n,v_1)\in V_n,\\
                \theta_n(v,u)&=(v_0,u_2-n,u_1+v_0)=(0,u_2-n,u_1)\in V_n,
            \end{align*}
            where
            $0\leq u_2-n\leq u_1\leq 1$ and
                  $0\leq v_2-n\leq v_1\leq 1$.
            In particular, $(v_2-n,v_1),(u_2-n,u_1)\in\{(0,0),(0,1),(1,1)\}$.
            In all cases, we have
            $\tau
            =\raisebox{-9.5mm}{
            \begin{tikzpicture}[auto]
                \tile{white}{0}{0}{(0,v_2-n,v_1)}{(0,u_2-n,u_1)}{(0,u_1,u_2)}{(0,v_1,v_2)}
            \end{tikzpicture}}
            \in J_n \cup \{j_n^{0,0,1,1}, j_n^{1,1,0,0}\}$
            is a junction tile.
    \end{itemize}

    ($\supseteq$)
    Proving $\Ccal_n\supseteq\Tcal_n'$ is not necessary to conclude the proof,
    since $\Ccal_n\subseteq\Tcal_n'$ and $\Tcal_n'$ is a finite set.
    Indeed, the set $\Tcal_n'$ contains $\#\Tcal_n'=n^2+8n+13$ elements.
    Also, in the proof that $\Ccal_n\subseteq\Tcal_n'$ made above,
    we exhibited 
    $n^2$ white tiles,
    $2(n+1)$ blue tiles,
    $2n$ antigreen tiles,
    $2(n+1)$ green tiles,
    $2n$ yellow tiles and
    $9$ junction tiles
    in $\Ccal_n$.
    Therefore, $\Ccal_n$ contains $n^2+2(n+1+n+n+1+n)+9=n^2+8n+13$ elements.
    We conclude that $\Ccal_n=\Tcal_n'$.

    Alternatively, $\Ccal_n\supseteq\Tcal_n'$ can be proved directly.
    One may check that for every 
    $\tau=
    \raisebox{-9.5mm}{
    \begin{tikzpicture}[auto]
        \tile{white}{0}{0}{r}{t}{\ell}{b}
    \end{tikzpicture}
    }
    \in
    \Tcal_n'$,
    we have
    $\{r,t,\ell,b\}\subset V_n$, $r=\theta_n(\ell,b)$ and $t=\theta_n(b,\ell)$.
    Thus, $\tau\in\Ccal_n$.
\end{proof}

We may now prove the first main result.

\begin{THEOREMA}
    \MainTheoremA
\end{THEOREMA}

\begin{proof}
    From Proposition~\ref{prop:Ccal-equals-Tcaln'}, we have
    $\Ccal_n=\Tcal_n'$.
    It was shown in \cite{labbe_metallic_I_2025}
    that the tiles in the difference set
    $\Tcal_n'\setminus\Tcal_n$ do not appear in valid configurations
    of $\Omega_{\Tcal_n'}$, so that $\Omega_{\Tcal_n'}=\Omega_{\Tcal_n}$.
    Thus, we conclude the equalities
    \[
        \Omega_{\Ccal_n}
        = \Omega_{\Tcal_n'}
        = \Omega_{\Tcal_n}
        = \Omega_n.\qedhere
    \]
\end{proof}

Now, we show that the computation performed by $\theta_n$ is invertible.
Let 
\[
\begin{array}{rccl}
    \psi_n:&V_n\times V_n & \to & \Z^3\\
    &(r_0,r_1,r_2), (t_0,t_1,t_2) 
    & \mapsto & (\ell_0,\ell_1,\ell_2),
\end{array}
\]
be the function defined by
\begin{equation}\label{eq:defintion-map-psi}
\left\{
\begin{array}{ll}
    &\ell_0=r_0,\\
    &\ell_1=\begin{cases}
            t_2-t_0  & \text{ if } r_0 = 0,\\
            1        & \text{ if } r_0 = 1,
         \end{cases}\\
    &\ell_2=\begin{cases}
            t_1+n  & \text{ if } t_0 = 0,\\
            r_2-1  & \text{ if } t_0 = 1.
         \end{cases}
\end{array}
\right.
\end{equation}
The following proposition states that the south and west colors 
of tiles in $\Ccal_n$
can be deduced from the right and top colors
using the map $\psi_n$.

\begin{proposition}\label{prop:inverse-of-theta}
    We have
\begin{equation}
    \Ccal_n=
    \left\{
    \raisebox{-9.5mm}{
    \begin{tikzpicture}[auto]
        \tile{white}{0}{0}{r}{t}{\psi_n(r,t)}{\psi_n(t,r)}
    \end{tikzpicture}
    }
    \middle|\,
        r,t\in V_n
        \text{ such that }
        \psi_n(r,t),\psi_n(t,r)\in V_n
    \right\}.
\end{equation}
\end{proposition}

\begin{proof}
    Let $\ell,b\in V_n$ and
    suppose that $r=(r_0,r_1,r_2)=\theta_n(\ell,b)$ and
                 $t=(t_0,t_1,t_2)=\theta_n(b,\ell)$.
    From Equation~\eqref{eq:defintion-map-theta}, we have
    \begin{equation}\label{eq:r-and-t-from-l-and-b}
    \left\{
    \begin{array}{ll}
        &r_0=\ell_0,\\
        &r_1=\begin{cases}
                b_2-n     & \text{ if } \ell_0 = 0,\\
                1         & \text{ if } \ell_0 = 1,
            \end{cases}\\
        &r_2=\begin{cases}
                b_1+\ell_0   & \text{ if } b_0 = 0,\\
                \ell_2+1     & \text{ if } b_0 = 1,
            \end{cases}
    \end{array}
    \right.
    \quad
    \text{ and }
    \quad
    \left\{
    \begin{array}{ll}
        &t_0=b_0,\\
        &t_1=\begin{cases}
                \ell_2-n     & \text{ if } b_0 = 0,\\
                1         & \text{ if } b_0 = 1,
            \end{cases}\\
        &t_2=\begin{cases}
                \ell_1+b_0   & \text{ if } \ell_0 = 0,\\
                b_2+1     & \text{ if } \ell_0 = 1.
            \end{cases}
    \end{array}
    \right.
    \end{equation}
    The above holds if and only if
    \[
    \left\{
    \begin{array}{ll}
        &\ell_0=r_0,\\
        &\ell_1=\begin{cases}
                t_2-t_0  & \text{ if } r_0 = 0,\\
                1        & \text{ if } r_0 = 1,
            \end{cases}\\
        &\ell_2=\begin{cases}
                t_1+n  & \text{ if } t_0 = 0,\\
                r_2-1  & \text{ if } t_0 = 1,
            \end{cases}
    \end{array}
    \right.
    \quad
    \text{ and }
    \quad
    \left\{
    \begin{array}{ll}
        &b_0=t_0,\\
        &b_1=\begin{cases}
                r_2-r_0  & \text{ if } t_0 = 0,\\
                1        & \text{ if } t_0 = 1,
            \end{cases}\\
        &b_2=\begin{cases}
                r_1+n  & \text{ if } r_0 = 0,\\
                t_2-1  & \text{ if } r_0 = 1.
            \end{cases}
    \end{array}
    \right.
    \]
    if and only if $\ell=(\ell_0,\ell_1,\ell_2)=\psi_n(r,t)$ and
                   $b=(b_0,b_1,b_2)=\psi_n(t,r)$.
    Thus, from Equation~\eqref{eq:Ccal}, we have
    \begin{align*}
        \Ccal_n&=
        \left\{
        \raisebox{-9.5mm}{
        \begin{tikzpicture}[auto]
            \tile{white}{0}{0}{\theta_n(\ell,b)}{\theta_n(b,\ell)}{\ell}{b}
        \end{tikzpicture}
        }
        \middle|\,
            \ell,b\in V_n
            \text{ such that }
            \theta_n(\ell,b),\theta_n(b,\ell)\in V_n
        \right\}\\
        &=
        \left\{
        \raisebox{-9.5mm}{
        \begin{tikzpicture}[auto]
            \tile{white}{0}{0}{r}{t}{\psi_n(r,t)}{\psi_n(t,r)}
        \end{tikzpicture}
        }
        \middle|\,
            r,t\in V_n
            \text{ such that }
            \psi_n(r,t),\psi_n(t,r)\in V_n
        \right\}.\qedhere
    \end{align*}
\end{proof}

As a consequence of
Proposition~\ref{prop:inverse-of-theta},
there is a bijection between the south-west and the north-east colors
for the tiles in $\Ccal_n$.
Using the vocabulary of \cite{MR1692474}, 
we may state the following result.
A set $\Tcal$ of Wang tiles is called \defn{SW-deterministic} if there do not exist two
different tiles in $\Tcal$ that would have same colors on their bottom and left edges,
respectively. In other words, for all colors $C_1$ and $C_2$
there exists at most one tile in $\Tcal$ whose bottom and left edges have
colors $C_1$ and $C_2$, respectively.
\defn{NW-}, \defn{NE-} and \defn{SE-deterministic} sets of Wang tiles are
defined analogously.
Thus, we obtain a conceptual proof for a result already obtained in
\cite{labbe_metallic_I_2025}.

\begin{theorem}
    [{\cite[Lemma~4.3]{labbe_metallic_I_2025}}]
    \label{thm:NE-SW-deterministic}
    For every integer $n\geq1$, the set of Wang tiles $\Ccal_n$ is
    NE-deterministic and SW-deterministic.
\end{theorem}

\begin{proof}
The set of Wang tile $\Ccal_n$ is SW-deterministic by definition
and NE-deterministic from Proposition~\ref{prop:inverse-of-theta}.
\end{proof}

\section{Equations satisfied by the Wang tiles and their tilings}\label{sec:equations}

In this section, we show that the set $\Ccal_n$ of Wang tiles satisfy a system
of equations. Moreover, we show that the rectangular 
tilings (of sizes $h\times 1$, $\infty\times 1$ and $h\times k$) generated by
them satisfy equations.
While the equations associated with Kari's \cite{MR1417578} and Culik's
\cite{MR1417576} aperiodic sets of Wang tiles are multiplicative, the ones
associated with $\Ccal_n$ are additive.

In the next theorem, we show that tiles in $\Ccal_n$ satisfy 
$\ell_0=r_0$, $b_0=t_0$ and the equation
\[
    \frac{t_2-t_1+\ell_2-\ell_1}{n}-\ell_0
    =
    \frac{b_2-b_1+r_2-r_1}{n}-b_0
\]
which reminds of Equation~\eqref{eq:wang-kari}.

\begin{THEOREMB}
    \MainTheoremB
\end{THEOREMB}

\begin{proof}
    Let
    $\ell=(\ell_0,\ell_1,\ell_2)$,
    $b=(b_0,b_1,b_2)$,
    $r=(r_0,r_1,r_2)$ and
    $t=(t_0,t_1,t_2)$.
    We always have $r_0=\ell_0$ and $t_0=b_0$.
    Thus, $\langle e, \ell\rangle=\ell_0=r_0=\langle e, r\rangle$ and 
         $\langle e, b\rangle=b_0   =t_0=\langle e, t\rangle$.
    Moreover,
    \begin{align*}
        \langle d, b\rangle &= b_2-b_1,\\
        \langle d, \ell\rangle &= \ell_2-\ell_1.
    \end{align*}
    The proof of the remaining equality is split in four cases.
    We use Equation~\eqref{eq:r-and-t-from-l-and-b} in the computations below.
    \begin{itemize}
\item If $(b_0,\ell_0)=(0,0)$, then
    \begin{align*}
        \langle d, t+\ell\rangle &= (t_2-t_1) + (\ell_2-\ell_1) 
                      = (\ell_1+b_0)-(\ell_2-n) + (\ell_2-\ell_1) 
                      = b_0 + n = n\\
        \langle d, r+b\rangle &= (r_2-r_1) + (b_2-b_1) =(b_1+\ell_0)-(b_2-n) + (b_2-b_1) = \ell_0 + n = n\\
        n\langle e, \ell-b\rangle &= n(\ell_0-b_0) = 0
    \end{align*}
\item If $(b_0,\ell_0)=(0,1)$, then $\ell_1=1$ and
    \begin{align*}
        \langle d, t+\ell\rangle &= (t_2-t_1) + (\ell_2-\ell_1)  
                      =(b_2+1)-(\ell_2-n) + (\ell_2-\ell_1)  
                      = b_2 + n\\
        \langle d, r+b\rangle &= (r_2-r_1) + (b_2-b_1) =(b_1+\ell_0)-(1) + (b_2-b_1) = b_2\\
        n\langle e, \ell-b\rangle &= n(\ell_0-b_0) = n
    \end{align*}
\item If $(b_0,\ell_0)=(1,0)$, then $b_1=1$ and
    \begin{align*}
        \langle d, t+\ell\rangle &= (t_2-t_1) + (\ell_2-\ell_1)  
                      =(\ell_1+b_0)-(1) + (\ell_2-\ell_1) = \ell_2\\
        \langle d, r+b\rangle &= (r_2-r_1) + (b_2-b_1) 
                      =(\ell_2+1)-(b_2-n) + (b_2-b_1)  
                      = \ell_2 + n\\
        n\langle e, \ell-b\rangle &= n(\ell_0-b_0) = -n
    \end{align*}
\item If $(b_0,\ell_0)=(1,1)$, then $b_1=\ell_1=1$ and
    \begin{align*}
        \langle d, t+\ell\rangle &= (t_2-t_1) + (\ell_2-\ell_1) 
                                  =(b_2+1)-(1) + (\ell_2-\ell_1) 
                                = b_2 + \ell_2 - \ell_1\\
        \langle d, r+b\rangle &= (r_2-r_1) + (b_2-b_1)
                               =(\ell_2+1)-(1)+ (b_2-b_1) 
                                  = \ell_2 + b_2 - b_1\\
        n\langle e, \ell-b\rangle &= n(\ell_0-b_0) = 0
    \end{align*}
\end{itemize}
    In all the four cases, we have
    $\langle d, t+\ell\rangle=\langle d, r+b\rangle+n\langle e, \ell-b\rangle$.
\end{proof}

The two sets in the statement of Theorem~\ref{thm:equations-satisfied-by-tiles}
are not equal. For instance
    $\raisebox{-11mm}{
    \begin{tikzpicture}[auto]
        \tile{white}{0}{0}{(1, 1, 3)}{(0, 0, 3)}{(1, 1, 5)}{(0, 0, 1)}
    \end{tikzpicture}}$
satisfy the equations when $n=4$, but it is not a tile in $\Ccal_n$.

Equation~\eqref{eq:wang-kari} behaves well with valid tiling of an
horizontal strip by Wang tiles associated with the same multiplication factor
$q\in\Q$. The same holds with tiles in $\Ccal_n$ which are related to some
addition of a certain value modulo 1.

\begin{figure}[h]
\begin{center}
    \includegraphics{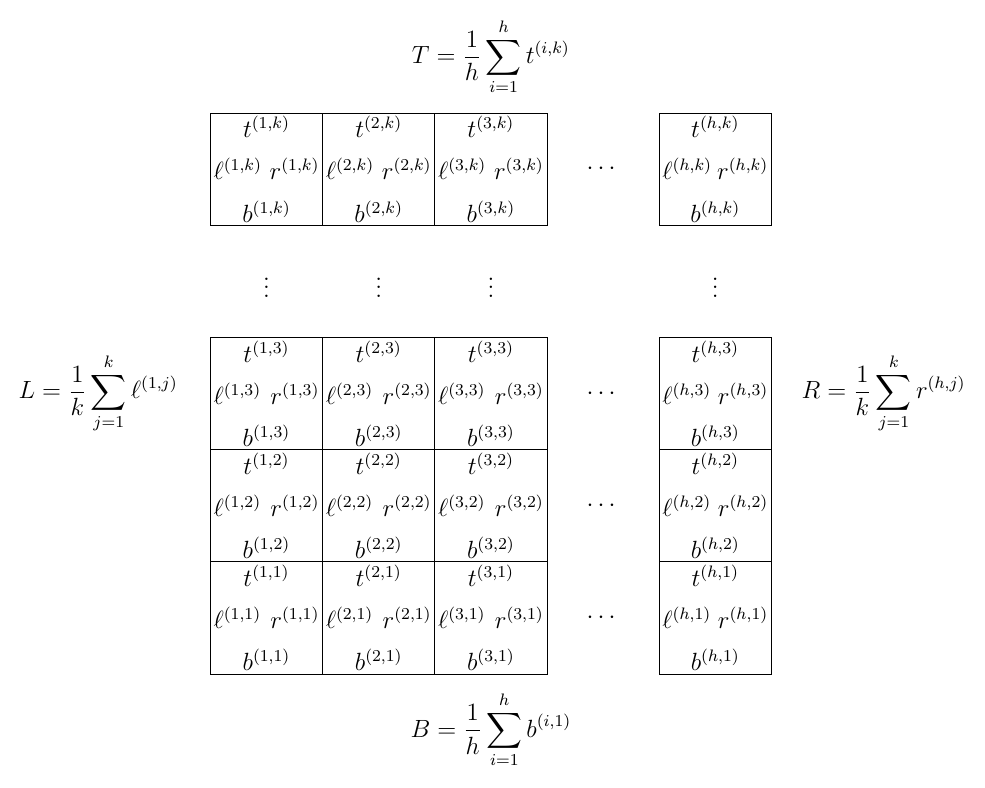}
\end{center}
    \caption{A $h\times k$ rectangular tiling of tiles from $\Ccal_n$.}
    \label{fig:rectangle-tiling}
\end{figure}

The equation satisfied by the tiles proved in
Theorem~\ref{thm:equations-satisfied-by-tiles}
extends to an equation for $h\times k$ rectangular valid tilings.

\begin{lemma}\label{lem:equations-for-a-rectangle}
    Let $n,h,k\geq1$ be integers and $d=(0,-1,1)$ and $e=(1,0,0)$. 
    Let 
    \[
        \{(r^{(i,j)},t^{(i,j)},\ell^{(i,j)},b^{(i,j)})\}_{1\leq i\leq h,1\leq j\leq k}
    \]
    be a family of tiles in $\Ccal_n$
    forming a valid tiling of a $h\times k$ rectangle, 
    see Figure~\ref{fig:rectangle-tiling}. 
    Let
    \[
        R=\textstyle\frac{1}{k}\sum_{j=1}^k r^{(h,j)},\quad
        T=\textstyle\frac{1}{h}\sum_{i=1}^h t^{(i,k)},\quad
        L=\textstyle\frac{1}{k}\sum_{j=1}^k\ell^{(1,j)}\quad\text{ and }\quad
        B=\textstyle\frac{1}{h}\sum_{i=1}^h b^{(i,1)}
    \]
    be the average of the right, top, left and bottom labels of the rectangular tiling.
    Then the following equation holds
        \begin{equation}\label{eq:equations-for-a-rectangle}
    \frac{1}{k}\left\langle \textstyle\frac{1}{n}d, T - B\right\rangle 
    -\langle e,L\rangle
    =
    \frac{1}{h} \left\langle \textstyle\frac{1}{n}d, R - L \right\rangle
    - \langle e, B \rangle.
        \end{equation}
\end{lemma}

\begin{proof}
    From Theorem~\ref{thm:equations-satisfied-by-tiles}, we have
        $\langle e, \ell^{(i,j)}\rangle=\langle e, r^{(i,j)}\rangle$,
        $\langle e, b^{(i,j)}\rangle=\langle e, t^{(i,j)}\rangle$ and
\[
    \langle \textstyle\frac{1}{n}d, t^{(i,j)} - b^{(i,j)}\rangle 
        - \langle e , \ell^{(i,j)} \rangle
        = \langle \textstyle\frac{1}{n}d, r^{(i,j)} - \ell^{(i,j)} \rangle 
        - \langle e ,  b^{(i,j)} \rangle,
\]
for every integers $i$ and $j$ such that $1\leq i\leq h$ and $1\leq j\leq k$.
We have
{\allowdisplaybreaks
\begin{align*}
    \frac{1}{k}\left\langle \textstyle\frac{1}{n}d, T - B\right\rangle 
    -\langle e,L\rangle
    &=
    \frac{1}{k}
    \left\langle \textstyle\frac{1}{n}d, \frac{1}{h}\sum_{i=1}^ht^{(i,k)} - 
                                         \frac{1}{h}\sum_{i=1}^hb^{(i,1)}\right\rangle
    -\langle e, \textstyle\frac{1}{k}\sum_{j=1}^k\ell^{(1,j)} \rangle\\
    &=
    \frac{1}{kh}\sum_{i=1}^h
    \left\langle \textstyle\frac{1}{n}d, t^{(i,k)} - b^{(i,1)}\right\rangle
    -\frac{1}{k}\sum_{j=1}^k\langle e, \ell^{(1,j)} \rangle\\
    &=
    \frac{1}{kh}\sum_{i=1}^h
    \left\langle \textstyle\frac{1}{n}d, 
    \sum_{j=1}^k t^{(i,j)} - 
    \sum_{j=1}^k b^{(i,j)}
    \right\rangle
    -\frac{1}{k}\sum_{j=1}^k\langle e, \textstyle\frac{1}{h}\sum_{i=1}^h\ell^{(i,j)} \rangle\\
    &=
    \frac{1}{kh}\sum_{i=1}^h\sum_{j=1}^k
    \left(
    \left\langle \textstyle\frac{1}{n}d, 
     t^{(i,j)} - b^{(i,j)}
    \right\rangle 
    -\langle e, \ell^{(i,j)} \rangle\right)\\
    &=
    \frac{1}{kh}\sum_{i=1}^h\sum_{j=1}^k
    \left(
    \left\langle \textstyle\frac{1}{n}d, 
     r^{(i,j)} - \ell^{(i,j)}
    \right\rangle
    -
    \left\langle e,
     b^{(i,j)}
    \right\rangle\right) \\
    &=
    \frac{1}{kh}\sum_{j=1}^k
    \left\langle \textstyle\frac{1}{n}d, 
     \sum_{i=1}^hr^{(i,j)} - \sum_{i=1}^h\ell^{(i,j)}
    \right\rangle
    -
    \frac{1}{h}\sum_{i=1}^h
    \left\langle e,
     \textstyle\frac{1}{k}\sum_{j=1}^kb^{(i,j)}
    \right\rangle \\
    &=
    \frac{1}{kh}\sum_{j=1}^k
    \left\langle \textstyle\frac{1}{n}d, 
     r^{(h,j)} - \ell^{(1,j)}
    \right\rangle
    -
    \frac{1}{h}\sum_{i=1}^h
    \left\langle e,
     b^{(i,1)}
    \right\rangle \\
    &=
    \frac{1}{h}
    \left\langle \textstyle\frac{1}{n}d, 
     \frac{1}{k}\sum_{j=1}^kr^{(h,j)} - \frac{1}{k}\sum_{j=1}^k\ell^{(1,j)}
    \right\rangle
    -
    \left\langle e,
    \textstyle\frac{1}{h}\sum_{i=1}^h
     b^{(i,1)}
    \right\rangle \\
    &=
    \frac{1}{h}
    \left\langle \textstyle\frac{1}{n}d, 
     R - L
    \right\rangle
    -
    \left\langle e, B \right\rangle.\qedhere
\end{align*}}
\end{proof}

Equation~\eqref{eq:equations-for-a-rectangle} is a simple consequence of the equations
satisfied by the tiles, but it has important implications.
If $L=R$, then
    $\left\langle \textstyle\frac{1}{n}d, 
     R - L
    \right\rangle=0$
    and $k\langle e,L\rangle$ is an integer.
    Thus,
    the average of the inner product with $\textstyle\frac{1}{n}d$ of the top labels is obtained 
    from the average of the inner product with $\textstyle\frac{1}{n}d$ of the bottom
    labels by $k$ rotations on the unit circle by a fixed angle:
    \begin{equation}\label{eq:equations-for-a-rectangle-with-equal-sides}
    \langle\textstyle\frac{1}{n} d,T\rangle 
    = \langle\textstyle\frac{1}{n} d,B\rangle 
    - k\langle e,B\rangle \pmod 1.
    \end{equation}

If $\Omega_n$ admits a periodic tiling, then there exists
a $h\times k$ rectangular tiling of tiles from $\Ccal_n$ such that
$L=R$ and $B=T$.
From Equation~\eqref{eq:equations-for-a-rectangle}, 
we get that 
$\left\langle e,L \right\rangle 
=\left\langle e,B \right\rangle$.
This equation means that the frequency of rows with no junction tiles
is equal to the frequency of columns with no junction tiles.
This holds if and only if
$h$ times the number of rows with no junction tile
is equal to 
$k$ times the number of columns with no junction tiles.
Copies of the $h\times k$ rectangular tiling can be used
to tile periodically a $hk\times hk$ square
respecting all matching rules
containing as many rows with no junction tile
as columns with no junction tile.
But this is not sufficient to prove that no periodic tiling exist.

Kari's \cite{MR1417578} and Culik's \cite{MR1417576} equations 
allow to show in a few lines that their sets of Wang tiles admit no periodic
tiling.
Proving the same for $\Omega_n$ directly from the equations remains an open question.

\section{Valid tilings obtained from floors of linear forms}
\label{sec:valid-tilngs-as-codings}

In this section, we present a method to construct valid tilings in $\Omega_n$.
It is based on the integer-floor value of three specific linear form over two
variables.

\begin{figure}[h]
\[
    \begin{array}{cccc}
         \lfloor y+\beta^*+1\rfloor
        &\lfloor \beta^{-1}x + y+\beta^*+1\rfloor
        &\lfloor \beta x + y+\beta^*+1\rfloor
        &\Lambda_n(x,y)\\[2mm]
         \includegraphics[height=4cm]{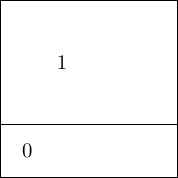}
        &\includegraphics[height=4cm]{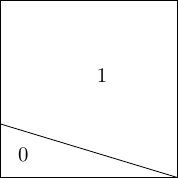}
        &\includegraphics[height=4cm]{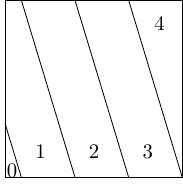}
        &\includegraphics[height=4cm]{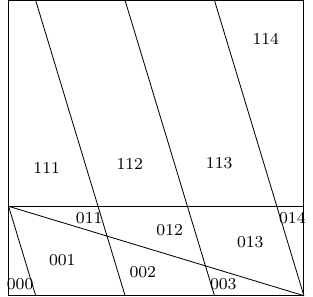}
    \end{array}
\]
    \caption{The preimage sets of the map $(x,y)\mapsto
        \Lambda_n(x,y)$ defines a partition of $[0,1)^2$ which
        is the refinement of the three partitions on the left.
        The above images are when $n=3$.}
        \label{fig:Lambda_n}
\end{figure}

Let $n\geq1$ be an integer and let $\beta$ be the positive root of $x^2- nx-1$.
We denote the negative root by $\beta^*$ which satisfies $\beta\beta^*=-1$ and $\beta+\beta^*=n$.
We consider the matrix
\[
    M_n=\left(\begin{array}{cc}
    0 & 1 \\
    \beta^{-1} & 1 \\
    \beta & 1
\end{array}\right)
\]
and the map $\lambda_n:\R^2\to\R^3$ defined by
\[
    \lambda_n(x,y)=M_n\cdot
    \left(\begin{array}{cc}
        \{x\} \\\{y\}
\end{array}\right)
+
\left(\begin{array}{cc}
    \beta^*+1 \\
    \beta^*+1 \\
    \beta^*+1
\end{array}\right)
\]
where $\{x\}=x-\lfloor x\rfloor$ is the fractional part of $x$.
Since $\lambda_n(x,y)=\lambda_n(x+1,y)=\lambda_n(x,y+1)$,
it is also well-defined on the torus $\lambda_n:\torus^2\to\R^3$.
Then, we define a coding function $\Lambda_n$ as the coordinate-wise floor of $\lambda_n$
when restricted to the domain $[0,1)^2$.
More precisely, we have
\[
\begin{array}{rcl}
    \Lambda_n:[0,1)^2 & \to & \Z^3\\
    (x,y) & \mapsto &
             \left(
             \begin{array}{r}
                 \lfloor y+\beta^*              +1\rfloor\\
                 \lfloor \beta^{-1}x + y+\beta^*+1\rfloor\\
                 \lfloor \beta x + y+\beta^*    +1\rfloor
             \end{array}
             \right),
\end{array}
\]
see Figure~\ref{fig:Lambda_n}.

Recall that, for every integer $n\geq 1$, we have
\[
    V_n = \{(v_0,v_1,v_2)\in\Z^3\colon 
                   0\leq v_0\leq v_1\leq v_2\leq n+1\text{ and } v_1\leq 1 \}.
\]

\begin{lemma}\label{lem:ceillambdan-nondecreasing}
    For every $(x,y)\in[0,1)^2$,
    $\Lambda_n(x,y)\in V_n$.
\end{lemma}

\begin{proof}
    Let $(x,y)\in[0,1)^2$.
    Since $\beta>1$, we have 
    \[
        0 <    \beta^*+1
          \leq y+\beta^*+1
          \leq \beta^{-1}x + y+\beta^*+1
          \leq \beta x + y+\beta^*+1
          < \beta + 1+\beta^*+1
          = n + 2.
    \]
    Thus, taking the floor function, we obtain
    \[
          0 \leq \lfloor \beta^* +1\rfloor
            \leq \lfloor y+\beta^*+1\rfloor
            \leq \lfloor \beta^{-1}x + y+\beta^*+1\rfloor
            \leq \lfloor \beta x + y+\beta^*+1\rfloor
            < n+2.
    \]
    Therefore, if $(v_0,v_1,v_2)=\Lambda_n(x,y)$, we have
    $0\leq v_0\leq v_1\leq v_2\leq n+1$.
    Also
    \[
         \beta^{-1}x + y+\beta^*+1
           <\beta^{-1} + 1+\beta^*+1
             =1+1= 2.
    \]
    Thus,
    \[
         v_1 = \lfloor \beta^{-1}x + y+\beta^*+1\rfloor
           \leq1.
    \]
    We conclude $\Lambda_n(x,y)=(v_0,v_1,v_2)\in V_n$.
\end{proof}

    The following lemma shows a relation between $\Lambda_n$ and the map 
    $\theta_n$ defined in Equation~\eqref{eq:defintion-map-theta}.

    \begin{lemma}\label{lem:satisfies-function-theta}
        If $x,y\in[0,1)$,
        then
        \[
        \Lambda_n(x,y)
        =\theta_n\big(
        \Lambda_n(\{x+\beta^*\},y),
        \Lambda_n(\{y+\beta^*\},x)
        \big).
        \]
    \end{lemma}

    \begin{proof}
        Let $x,y\in[0,1)$.
        We want to show that if $\ell_0, \ell_1, \ell_2,b_0,b_1,b_2\in\Z$
        are such that
        \[
            \Lambda_n(\{x+\beta^*\},y)
            =
             \left(
             \begin{array}{r}
                 \lfloor y+\beta^*+1\rfloor\\
                 \lfloor \beta^{-1}\{x+\beta^*\} + y+\beta^*+1\rfloor\\
                 \lfloor \beta \{x+\beta^*\} + y+\beta^*+1\rfloor
             \end{array}
             \right)
             =
             \left(
             \begin{array}{r}
                 \ell_0\\
                 \ell_1\\
                 \ell_2
             \end{array}
             \right)
        \]
        and
        \[
            \Lambda_n(\{y+\beta^*\},x)
            =
             \left(
             \begin{array}{r}
                 \lfloor x+\beta^*+1\rfloor\\
                 \lfloor \beta^{-1}\{y+\beta^*\} + x+\beta^*+1\rfloor\\
                 \lfloor \beta \{y+\beta^*\} + x+\beta^*+1\rfloor
             \end{array}
             \right)
             =
             \left(
             \begin{array}{r}
                 b_0\\
                 b_1\\
                 b_2
             \end{array}
             \right),
        \]
        then $\Lambda_n(x,y)=\theta_n\left((\ell_0,\ell_1,\ell_2),(b_0,b_1,b_2)\right)$.
        Let $r_0,r_1,r_2\in\Z$ be such that 
        \[
        \Lambda_n(x,y)=
             \left(
             \begin{array}{r}
                 \lfloor y+\beta^*+1\rfloor\\
                 \lfloor \beta^{-1}x + y+\beta^*+1\rfloor\\
                 \lfloor \beta x + y+\beta^*+1\rfloor
             \end{array}
             \right)
             =
             \left(
             \begin{array}{r}
                 r_0\\
                 r_1\\
                 r_2
             \end{array}
             \right).
        \]
        We want to show that the variables satisfy the definition
        of the function $\theta_n$ given in Equation~\eqref{eq:defintion-map-theta}.
        We have $r_0= \lfloor y+\beta^*+1\rfloor=\ell_0$.
        Therefore, the first equation defining the map $\theta_n$ is satisfied.

        Assume that $\ell_0=\lfloor y+\beta^*+1\rfloor=0$.
        Then $-\beta^{-1}=\beta^*\leq y+\beta^*<0$.
        Also $0\leq \beta^{-1}x <\beta^{-1}$.
        Thus, $-\beta^{-1}\leq\beta^{-1}x + y+\beta^*<\beta^{-1}$.
        We have
        \begin{align*}
            r_1 
            &=\lfloor \beta^{-1}x + y+\beta^*\rfloor+1\\
            &=\lfloor \beta(\beta^{-1}x + y+\beta^*)\rfloor +1
                            &(\text{because } {-}\beta^{-1}\leq\beta^{-1}x + y+\beta^*<\beta^{-1})\\
            &=\lfloor \beta (y+\beta^*) + x\rfloor +1      \\
            &=\lfloor \beta (y+\beta^*+1) + x+\beta^*\rfloor+1-n
                                       &(\text{because }\beta+\beta^*=n)\\
            &=\lfloor \beta \{y+\beta^*\} + x+\beta^*\rfloor+1-n\\
            &=b_2 -n
        \end{align*}

        Assume that $\ell_0=\lfloor y+\beta^*+1\rfloor=1$.
        Then $0\leq y+\beta^*< 1$.
        Also, we have $y<1$, so that $y+\beta^*< 1+\beta^*$.
        Moreover, 
        $0\leq \beta^{-1}x<\beta^{-1}$.
        Thus,
        $0< \beta^{-1}x + y+\beta^*< \beta^{-1} + 1 +\beta^* =1$.
        We have
        \[
            r_1 
            =\lfloor \beta^{-1}x + y+\beta^*\rfloor+1
            = 0 + 1 = \ell_0.
        \]
        Therefore, the second equation defining the map $\theta_n$ is satisfied.

        Assume that $b_0=\lfloor x+\beta^*+1\rfloor=0$.
        This implies that $-1\leq x+\beta^*<0$, which implies 
        $x<\beta^{-1}$.
        Thus, $0\leq \beta x< 1$.
        We need to consider the cases $\ell_0=0$ and $\ell_0=1$ separately.
        First, suppose that $\ell_0=\lfloor y+\beta^*+1\rfloor=0$.
        Then $-1\leq y+\beta^*< 0$.
        Thus, $-1\leq \beta x + y+\beta^*<1$.
        We have
        \begin{align*}
            r_2 
            &=\lfloor \beta x + y+\beta^*+1\rfloor\\
            &=\lfloor \beta^{-1}(\beta x + y+\beta^*)\rfloor +1
                            &(\text{because } -1\leq(\beta x + y+\beta^*)<1)\\
            &=\lfloor \beta^{-1}(\beta x + y+\beta^*) +\beta^{-1}+\beta^*\rfloor +1\\
            &=\lfloor \beta^{-1}(1+y+\beta^*) +x+\beta^*\rfloor +1\\
            &=\lfloor \beta^{-1}\{y+\beta^*\} +x+\beta^*\rfloor +1\\
            &= b_1 = b_1 + 0 = b_1+\ell_0.
        \end{align*}
        Secondly, suppose that $\ell_0=\lfloor y+\beta^*+1\rfloor=1$.
        Then $0\leq y+\beta^*< 1$,
        which implies $\{y+\beta^*\}=y+\beta^*$.
        Thus, $0\leq \beta x + y+\beta^*<2$.
        We have
        \begin{align*}
            r_2 
            &=\lfloor \beta x + y+\beta^*+1\rfloor\\
            &=\lfloor \beta x + y+\beta^*-1\rfloor+2\\
            &=\lfloor \beta^{-1}(\beta x + y+\beta^*-1) \rfloor+2
                            &(\text{because } -1\leq(\beta x + y+\beta^*-1)<1)\\
            &=\lfloor \beta^{-1}(y+\beta^*) +x +\beta^*\rfloor+2\\
            &=\lfloor \beta^{-1}\{y+\beta^*\} +x+\beta^*\rfloor+2\\
            &= b_1 + 1= b_1+\ell_0.
        \end{align*}

        Assume that $b_0=\lfloor x+\beta^*+1\rfloor=1$.
        This implies that $0\leq x+\beta^*<1$, 
        which implies $\{x+\beta^*\}=x+\beta^*$.
        We have
        \begin{align*}
            r_2 
            &=\lfloor \beta x + y+\beta^*+1\rfloor\\
            &=\lfloor \beta x +\beta\beta^*+1+ y+\beta^* +1\rfloor
                                       &(\text{because }\beta\beta^*=-1)\\
            &= \lfloor \beta (x+\beta^*) + y+\beta^*+1\rfloor + 1\\
            &= \lfloor \beta \{x+\beta^*\} + y+\beta^*+1\rfloor + 1\\
            &= \ell_2 + 1 = \ell_2+b_1.
        \end{align*}
        Therefore, the third equation defining the map $\theta_n$ is satisfied.
    \end{proof}

    For every $(x,y)\in\R^2$, let
    \[
        \sctile_n(x,y)=
        \left(
        \Lambda_n(\{x\},\{y\}),
        \Lambda_n(\{y\},\{x\}),
        \Lambda_n(\{x+\beta^*\},\{y\}),
        \Lambda_n(\{y+\beta^*\},\{x\})
        \right)
    \]
    which can be interpreted geometrically as a Wang tile:
    \[
        \sctile_n(x,y)=
    \raisebox{-11mm}{
    \begin{tikzpicture}[auto]
        \tile{white}{0}{0}{\Lambda_n(\{x\},\{y\})        }%
                          {\Lambda_n(\{y\},\{x\})        }%
                          {\Lambda_n(\{x+\beta^*\},\{y\})}%
                          {\Lambda_n(\{y+\beta^*\},\{x\})}%
    \end{tikzpicture}}
    \]

\begin{lemma}\label{lem:tilen-in-Tcal'n}
    If $(x,y)\in\R^2$,
    then
    \begin{itemize}
        \item $\widehat{\sctile_n}(x,y)=\sctile_n(y,x)$,
        \item $\sctile_n(x,y)\in(V_n)^4$,
        \item $\sctile_n(x,y)\in\Ccal_n$ is an instance of a $\theta_n$-chip tile.
    \end{itemize}
\end{lemma}

\begin{proof}
    We observe that
        $\sctile_n(x,y)$ is the image of $\sctile_n(y,x)$
        under the tile reflection $t\mapsto\widehat{t}$ by the positive slope diagonal.

From Lemma~\ref{lem:ceillambdan-nondecreasing},
for every $(x,y)\in[0,1)^2$, we have $\Lambda_n(x,y)\in V_n$.
Therefore, for every $(x,y)\in\R^2$,
\[
    \Lambda_n\left(\{x\},\{y\}\right),\quad
    \Lambda_n\left(\{y\},\{x\}\right),\quad
    \Lambda_n\left(\{x+\beta^*\},\{y\}\right),\quad
    \Lambda_n\left(\{y+\beta^*\},\{x\}\right)
    \in V_n.
\]
From Lemma~\ref{lem:satisfies-function-theta}, for every $(x,y)\in\R^2$, we have
\[
\Lambda_n(\{x\},\{y\})
=\theta_n\big(
\Lambda_n(\{x+\beta^*\},\{y\}),
\Lambda_n(\{y+\beta^*\},\{x\})
\big).
\]
Also
\[
    \Lambda_n(\{y\},\{x\})
=\theta_n\big(
\Lambda_n(\{y+\beta^*\},\{x\}),
\Lambda_n(\{x+\beta^*\},\{y\})
\big).
\]
Thus, $\sctile_n(x,y)\in\Ccal_n$.
\end{proof}

Here is another characterization of the set of Wang tiles $\Tcal_n$.

\begin{proposition}\label{prop:Tcaln_is_image_of_sctilexy}
    The following holds:
    \[
        \Tcal_n=\left\{\sctile_n(x,y) \colon (x,y)\in[0,1)^2\right\}.
    \]
\end{proposition}

\begin{proof}
    First, recall from Proposition~\ref{prop:Ccal-equals-Tcaln'} that
    \begin{equation}\label{eq:Tn'-sorted-is-Tn-plus-more-stuff}
        \Ccal_n
        = \Tcal_n'
        =
        \Tcal_n\cup
        \{j_n^{0,0,1,1},
          j_n^{1,1,0,0}\}
        \cup
        \left\{ a_n^i, \widehat{a_n^i} \mid 1\leq i \leq n \right\}
        \cup
        \left\{ b_n^n,\widehat{b_n^n} \right\}
    \end{equation}
    where
    \[
        \{j_n^{0,0,1,1},
          j_n^{1,1,0,0}\}
        = 
        \left\{
        \raisebox{-9mm}{
        \begin{tikzpicture}[auto]
            \tikzstyle{every node}=[font=\footnotesize]
            \tileJunctionOIIO{0}{0}{000}{011}{01\nbar}{00n}
            \tileJunctionIOOI{3}{0}{011}{000}{00n}{01\nbar}
        \end{tikzpicture}}
        \right\}.
    \]

    Let
    \[
        U_n=\left\{\sctile_n(x,y) \colon (x,y)\in[0,1)^2
        \right\}.
    \]
    First we show that $U_n\subseteq\Tcal_n$.
    It follows from Lemma~\ref{lem:tilen-in-Tcal'n} that
    $ U_n \subset \Ccal_n$.
    Thus, using Equation~\eqref{eq:Tn'-sorted-is-Tn-plus-more-stuff},
    the goal is to show that
    \begin{equation}\label{eq:U_n-cap-other=varnothing}
        U_n
        \cap\left(
        \{j_n^{0,0,1,1},
          j_n^{1,1,0,0}\}
        \cup
        \left\{ a_n^i, \widehat{a_n^i} \mid 1\leq i \leq n \right\}
        \cup
        \left\{ b_n^n,\widehat{b_n^n} \right\}
        \right)=\varnothing.
    \end{equation}

    Suppose that there exists
    $(x,y)\in[0,1)^2$
    such that $\sctile_n(x,y)=j_n^{0,0,1,1}$.
    Then $\Lambda_n(x,y)=000$ and $\Lambda_n(y,x)=011$.
    More precisely, we have
        \begin{align*}
        \Lambda_n(x,y)&=
             \left(
             \begin{array}{r}
                 \lfloor y+\beta^*+1\rfloor\\
                 \lfloor \beta^{-1}x + y+\beta^*+1\rfloor\\
                 \lfloor \beta x + y+\beta^*+1\rfloor
             \end{array}
             \right)
             = \left( \begin{array}{r} 0\\ 0\\ 0 \end{array} \right),\\
        \Lambda_n(y,x)&=
             \left(
             \begin{array}{r}
                 \lfloor x+\beta^*+1\rfloor\\
                 \lfloor \beta^{-1}y + x+\beta^*+1\rfloor\\
                 \lfloor \beta y + x+\beta^*+1\rfloor
             \end{array}
             \right)
             = \left( \begin{array}{r} 0\\ 1\\ 1 \end{array} \right).
        \end{align*}
    In particular,
    \[
     0 =    \lfloor \beta x + y+\beta^*+1\rfloor
       \geq \lfloor \beta^{-1}y + x+\beta^*+1\rfloor =1,
    \]
    which is a contradiction.
    The same contradiction is obtained if
    $\sctile_n(x,y)=j_n^{1,1,0,0}$.
    Therefore, these two junction tiles are not in $U_n$.

    Suppose that there exists
    $(x,y)\in[0,1)^2$
    such that $\sctile_n(x,y)=a_n^i$
    for some integer $i$ satisfying $1\leq i\leq n$.
    Then $\Lambda_n(x,y)=00\ibar$ and $\Lambda_n(y,x)=112$.
    More precisely, we have
        \begin{align*}
        \Lambda_n(x,y)&=
             \left(
             \begin{array}{r}
                 \lfloor y+\beta^*+1\rfloor\\
                 \lfloor \beta^{-1}x + y+\beta^*+1\rfloor\\
                 \lfloor \beta x + y+\beta^*+1\rfloor
             \end{array}
             \right)
             = \left( \begin{array}{c} 0\\ 0\\ i+1 \end{array} \right),\\
        \Lambda_n(y,x)&=
             \left(
             \begin{array}{r}
                 \lfloor x+\beta^*+1\rfloor\\
                 \lfloor \beta^{-1}y + x+\beta^*+1\rfloor\\
                 \lfloor \beta y + x+\beta^*+1\rfloor
             \end{array}
             \right)
             = \left( \begin{array}{c} 1\\ 1\\ 2 \end{array} \right).
        \end{align*}
    In particular, 
                 $\lfloor y+\beta^*+1\rfloor=0$
                 implies that $-\beta^{-1}\leq y+\beta^*<0$.
                 Also $0\leq\beta^{-1}x<\beta^{-1}$, so that
    $-\beta^{-1}\leq\beta^{-1}x+y+\beta^*<\beta^{-1}$.
    Therefore,
    \[
                 0=\lfloor \beta^{-1}x + y+\beta^*+1\rfloor
                  =\lfloor \beta(\beta^{-1}x + y+\beta^*)\rfloor+1
                  =\lfloor \beta y + x -1\rfloor+1
                  =\lfloor \beta y + x \rfloor.
    \]
    On the other hand,
    using $\lfloor a+b\rfloor\leq\lfloor a\rfloor+\lfloor b\rfloor+1$ for every $a,b\in\R$,
    we obtain
    \begin{align*}
                  2=\lfloor \beta y + x+\beta^*+1\rfloor
                   \leq\lfloor \beta y + x\rfloor+\lfloor\beta^*+1\rfloor+1
                   = 0 + 0 + 1 = 1,
    \end{align*}
    which is a contradiction.
    A similar contradiction is obtained if we suppose that
    such that $\sctile_n(x,y)=\widehat{a_n^i}$.
    Therefore, there is no antigreen tile in $U_n$.

    Suppose that there exists
    $(x,y)\in[0,1)^2$
    such that $\sctile_n(x,y)=b_n^n$.
    Then $\Lambda_n(x,y)=00\nbar$ and $\Lambda_n(y,x)=111$.
    More precisely, we have
        \begin{align*}
        \Lambda_n(x,y)&=
             \left(
             \begin{array}{r}
                 \lfloor y+\beta^*+1\rfloor\\
                 \lfloor \beta^{-1}x + y+\beta^*+1\rfloor\\
                 \lfloor \beta x + y+\beta^*+1\rfloor
             \end{array}
             \right)
             = \left( \begin{array}{c} 0\\ 0\\ n+1 \end{array} \right).
        \end{align*}
    In particular, using $\beta=n+\beta^{-1}$ and $x<1$, we obtain
    \begin{align*}
        n+1
        &=\lfloor \beta x + y+\beta^*+1\rfloor\\
        &=\lfloor (n+\beta^{-1}) x + y+\beta^*+1\rfloor\\
        &\leq\lfloor n+\beta^{-1} x + y+\beta^*+1\rfloor\\
        &=\lfloor \beta^{-1} x + y+\beta^*+1\rfloor+n
        = 0 + n = n,
    \end{align*}
    which is a contradiction.
    A similar contradiction is obtained if we suppose that
    such that $\sctile_n(x,y)=\widehat{b_n^n}$.
    Therefore, the blue tiles $b_n^n$ and $\widehat{b_n^n}$ are not in $U_n$.
    This shows that Equation~\eqref{eq:U_n-cap-other=varnothing} holds.
    Thus, $U_n\subseteq\Tcal_n$.

    Now we show that
    $\Tcal_n\subseteq U_n$.
    We have $J_n\subset U_n$ since
    \begin{align*}
    j_n^{0,0,0,0}        &= \sctile_n(0,0),\\
    j_n^{0,1,0,0}        &= \sctile_n\left(\beta^{-2},0\right),\\
    j_n^{0,0,0,1}        &= \sctile_n\left(0,\beta^{-2}\right),\\
    j_n^{0,1,0,1}        &= \sctile_n\left(\frac{1}{\beta(\beta+1)},\frac{1}{\beta(\beta+1)}\right),\\
    j_n^{1,1,0,1}    &= \sctile_n(x,y),
        \text{ where } (x,y) \text{ is on the segment from } (0,\beta^{-1}) \text{ to } 
        ((\beta+1)^{-1},(\beta+1)^{-1}),\\
    j_n^{0,1,1,1}    &= \sctile_n(x,y)
        \text{ where } (x,y) \text{ is on the segment from } (\beta^{-1},0) \text{ to } 
        ((\beta+1)^{-1},(\beta+1)^{-1}),\\
    j_n^{1,1,1,1}&= \sctile_n\left(\frac{1}{\beta+1},\frac{1}{\beta+1}\right).
    \end{align*}
    We have $B_n\subset U_n$ since
    \begin{align*}
        b_n^0 &= \sctile_n(\beta^{-1},0),\\
        b_n^i &= \sctile_n(\beta^{-2}+\beta^{-1}i,0) 
                 \text{ for every integer $i$ with $1\leq i\leq n-1$}. %
    \end{align*}
    We have $G_n\subset U_n$ since
    \begin{align*}
        g_n^0 &= \sctile_n(\beta^{-1},\beta^{-2}(\beta-1))\\
        g_n^i &= \sctile_n\left(\textstyle\frac{i}{n},\beta^{-1}(1-\frac{i}{n})\right)
                 \text{ for every integer $i$ with $1\leq i\leq n$}. %
    \end{align*}
    We have $Y_n\subset U_n$ since
    \begin{align*}
        y_n^1 &= \sctile_n(\beta^{-1}+\varepsilon,\beta^{-1}-\varepsilon\beta^{-1}) 
                     \text{ for some small }\varepsilon>0,\\
        y_n^i &= \sctile_n\left(\textstyle\frac{i-\beta^{-2}}{n},
                                \frac{\beta^{-1}}{n}(n-i+\beta^{-1}-\beta^{-2})\right)
                 \text{ for every integer $i$ with $2\leq i\leq n$}. %
    \end{align*}
    We have $W_n\subset U_n$ since
    \begin{align*}
        w_n^{1,1} &= \sctile_n(\beta^{-1},\beta^{-1}),\\
        w_n^{1,j} &= \sctile_n(\beta^{-1},j\beta^{-1}-\beta^{-2})
                 \text{ for every integer $j$ with $2\leq j\leq n$},\\
        w_n^{i,1} &= \sctile_n(i\beta^{-1}-\beta^{-2},\beta^{-1})
                 \text{ for every integer $i$ with $2\leq i\leq n$},\\
        w_n^{i,j} &= \sctile_n\left(
        \textstyle\beta^{-1}+\frac{1}{n}\left((i-1)-(j-1)\beta^{-1}\right), 
        \textstyle\beta^{-1}+\frac{1}{n}\left((j-1)-(i-1)\beta^{-1}\right)
        \right)\\
                 &\qquad\text{ for every integer $i,j$ with $2\leq i,j\leq n$}.
    \end{align*}
    Therefore, $J_n\cup B_n \cup G_n \cup Y_n\cup W_n \subseteq U_n$.
    Since $\widehat{U_n}=U_n$, we also have
    $\widehat{B_n} \cup \widehat{G_n} \cup \widehat{Y_n}\subseteq U_n$.
    We conclude that 
    $\Tcal_n \subseteq U_n$ and $\Tcal_n = U_n$.
\end{proof}

\begin{figure}[h]
\begin{center}
\begin{tikzpicture}[scale=3]
    \foreach \x in {1,2,3,4,5}
    \draw (\x+.5,1.4) -- (\x+.5,6-.4);
    \foreach \y in {1,2,3,4,5}
    \draw (1.4,\y+.5) -- (6-.4,\y+.5);
    \node at (2,5) {$\sctile_n\begin{pmatrix}x{-}\frac{1}{\beta}\\
                                            y{+}\frac{2}{\beta}\end{pmatrix}$};
    \node at (3,5) {$\sctile_n\begin{pmatrix}x\\y{+}\frac{2}{\beta}\end{pmatrix}$};
    \node at (4,5) {$\sctile_n\begin{pmatrix}x{+}\frac{1}{\beta}\\
                                            y{+}\frac{2}{\beta}\end{pmatrix}$};
    \node at (5,5) {$\sctile_n\begin{pmatrix}x{+}\frac{2}{\beta}\\
                                            y{+}\frac{2}{\beta}\end{pmatrix}$};
    \node at (2,4) {$\sctile_n\begin{pmatrix}x{-}\frac{1}{\beta}\\
                                            y{+}\frac{2}{\beta}\end{pmatrix}$};
    \node at (3,4) {$\sctile_n\begin{pmatrix}x\\y{+}\frac{1}{\beta}\end{pmatrix}$};
    \node at (4,4) {$\sctile_n\begin{pmatrix}x{+}\frac{1}{\beta}\\
                                            y{+}\frac{1}{\beta}\end{pmatrix}$};
    \node at (5,4) {$\sctile_n\begin{pmatrix}x{+}\frac{2}{\beta}\\
                                            y{+}\frac{1}{\beta}\end{pmatrix}$};
    \node at (2,3) {$\sctile_n\begin{pmatrix}x{-}\frac{1}{\beta}\\
                                            y\end{pmatrix}$};
    \node at (3,3) {$\sctile_n\begin{pmatrix}x\\y\end{pmatrix}$};
    \node at (4,3) {$\sctile_n\begin{pmatrix}x{+}\frac{1}{\beta}\\
                                            y\end{pmatrix}$};
    \node at (5,3) {$\sctile_n\begin{pmatrix}x{+}\frac{2}{\beta}\\
                                            y\end{pmatrix}$};
    \node at (2,2) {$\sctile_n\begin{pmatrix}x{-}\frac{1}{\beta}\\
                                            y{-}\frac{1}{\beta}\end{pmatrix}$};
    \node at (3,2) {$\sctile_n\begin{pmatrix}x\\y{-}\frac{1}{\beta}\end{pmatrix}$};
    \node at (4,2) {$\sctile_n\begin{pmatrix}x{+}\frac{1}{\beta}\\
                                            y{-}\frac{1}{\beta}\end{pmatrix}$};
    \node at (5,2) {$\sctile_n\begin{pmatrix}x{+}\frac{2}{\beta}\\
                                            y{-}\frac{1}{\beta}\end{pmatrix}$};
\end{tikzpicture}
\end{center}
    \caption{For every $(x,y)\in[0,1)^2$ the map
    $\Z^2\to\Tcal_n$ defined by
    $(i,j)\mapsto
    \sctile_n(x{+}\frac{i}{\beta},y{+}\frac{j}{\beta})$
    is a valid tiling of the plane by the set of Wang tiles $\Tcal_n$.}
    \label{fig:xy-to-tiling}
\end{figure}

This allows to construct valid configurations $\Z^2\to\Tcal_n$ from any
starting point $(x,y)$ on the torus. 
See Figure~\ref{fig:xy-to-tiling}.

\begin{THEOREMC}
    \MainTheoremC
\end{THEOREMC}

\begin{proof}
    Let $(x,y)\in[0,1)^2$ and $(i,j)\in\Z^2$. 
    We have $c_{(x,y)}(i,j)\in\Tcal_n$ from Proposition~\ref{prop:Tcaln_is_image_of_sctilexy}.
    Also the right color of the tile $c_{(x,y)}(i,j)$ is 
    $\Lambda_n(\{x+i\beta^{-1}\},\{y+j\beta^{-1}\})$
    which is equal to the left color of the tile $c_{(x,y)}(i+1,j)$.
    Finally, the top color of the tile $c_{(x,y)}(i,j)$ is 
    $\Lambda_n(\{y+j\beta^{-1}\},\{x+i\beta^{-1}\})$
    which is equal to the bottom color of the tile $c_{(x,y)}(i,j+1)$.
    Therefore, $c_{(x,y)}$ is a valid configuration of Wang tiles from the set $\Tcal_n$.
\end{proof}

The set $\{c_{(x,y)}\colon (x,y)\in[0,1)^2\}$ is not a subshift because it is
not topologically closed. Indeed, if $(x_0,y_0)$ lies on the boundary of the partition,
there is more than one configuration associated with it. 
The configuration $c_{(x_0,y_0)}$ is one of them, but
$\lim_{(x,y)\to(x_0,y_0)}c_{(x,y)}$ might be a different configuration if the limit is taken
coming from another direction.
The same issue happens with the representation of numbers in base 10.
For example, the number 1 has two base-10 representations, one being
$1.000000\dots$ and the other $0.999999\dots$.

This implies that the set
$\{c_{(x,y)}\colon (x,y)\in[0,1)^2\}$ 
is not the set of all valid configurations of $\Tcal_n$.
In other terms, $c:(x,y)\mapsto c_{(x,y)}$ is not surjective in
the set $\Omega_n$ of all valid configurations of $\Tcal_n$.
One way to solve this issue is to take the topological closure
\[
    C = \overline{\left\{c_{(x,y)}\colon (x,y)\in[0,1)^2\right\}}
\]
which is a nonempty subshift satisfying $C\subseteq\Omega_n$.
Since $\Omega_n$ is minimal \cite{labbe_metallic_I_2025}, we conclude 
the equality $C=\Omega_n$ must hold.

A standard approach is to create the subshift $C$ as the symbolic
extension of a dynamical system defined on the 2-torus $\torus^2$. 
This is what we do in the next two sections.

\section{An explicit factor map}\label{sec:explicit-factor-map}

The goal of this section is to introduce a factor map
$\Omega_n\to\torus^2$ explicitly defined from the average of
inner products of the labels of the Wang tiles in a configuration,
see Equation~\eqref{eq:Phi}.
Then, we prove Theorem~\ref{thm:factor-map} using this explicit factor map.

First, it is convenient to make some observation
on the inner product with the vector $d=(0,-1,1)$ of the tile labels.
In the statement below, we use the indicator function
$\I_{A}\colon\R\to\{0,1\}$ of a subset $A\subset\R$ defined as
\[
    \I_A(x)=
    \begin{cases}
        1 & \text{ if } x\in A,\\
        0 & \text{ if } x\notin A.
    \end{cases}
\]

    \begin{lemma}\label{lem:inner-product-with-d-of-Lambda_n}
        Let $n\geq1$ be an integer and
        $d=(0,-1,1)$.
        If $x,y\in[0,1)$, then
        \[
            \langle d, \Lambda_n(x,y)\rangle
        =
            \lfloor nx\rfloor + \I_{[1-\{nx\},1)}(\{\delta_{x} + y\})
        \]
        where $\delta_{x}=1-\beta^{-1}(1-x)$.
    \end{lemma}

    \begin{proof}
        Let $x,y\in[0,1)$.
        Observe that $\delta_{x}=1-\beta^{-1}(1-x)=\beta^{-1}x +\beta^*+1$.
        We have
        \begin{align*}
            \langle d, \Lambda_n(x,y)\rangle
            &=   \lfloor \beta x + y+\beta^*+1\rfloor
               - \lfloor \beta^{-1}x + y+\beta^*+1\rfloor\\
            &=   \lfloor (n+\beta^{-1}) x + y+\beta^*+1\rfloor
               - \lfloor \beta^{-1}x + y+\beta^*+1\rfloor\\
            &=   \lfloor nx+\delta_{x}+y\rfloor
               - \lfloor    \delta_{x}+y\rfloor\\
            &=   \left( \lfloor nx\rfloor
                       +\lfloor \delta_{x}+y\rfloor
                       +\lfloor \{nx\}+\{\delta_{x}+y\}\rfloor\right)
               - \lfloor    \delta_{x}+y\rfloor\\
            &=   \lfloor nx\rfloor
                 +\lfloor \{nx\}+\{\delta_{x}+y\}\rfloor\\
            &= \lfloor nx\rfloor +
            \begin{cases} 
                0  & \text{ if } \{nx\}+\{\delta_{x}+y\} < 1,\\
                1  & \text{ if } \{nx\}+\{\delta_{x}+y\} \geq 1.
            \end{cases}
        \end{align*}
        The conclusion follows.
    \end{proof}

As illustrated in Figure~\ref{fig:T3-tiling-10x5-with-values}
for a finite rectangular pattern,
the average of the values of $\langle\frac{1}{n}d,v\rangle$ 
for labels $v$ appearing along an horizontal line
can be considered for valid configurations $w:\Z^2\to\Tcal_n$.
For some reason
(in order to have the equality $\phi_n(c_{(x,y)}) = y$ in
Proposition~\ref{prop:phi-is-a-factor-map}), 
it is convenient to consider the average of the top label of
the tiles on the horizontal row passing through the origin.
Assuming that the limit exists for every configuration, 
this leads to a map from the Wang shift to the
interval $[0,1]$ defined as follows:
\begin{equation}\label{eq:phi_n}
\begin{array}{rcl}
    \phi_n:\Omega_n & \to & [0,1]\\
w & \mapsto &
\displaystyle
    \lim_{k\to\infty}\frac{1}{2k+1}\sum_{i=-k}^k
    \langle \textstyle\frac{1}{n} d,\sctop(w_{i,0})\rangle
\end{array}
\end{equation}
where $\sctop(t)$ denotes the top label of the Wang tile $t$.

We show in the next proposition that $\phi_n$ is well-defined and that
it recovers the parameter $y$ of a configuration $c_{(x,y)}$.

\begin{proposition}\label{prop:phi-is-a-factor-map}
    For every integer $n\geq1$, the following holds:
    \begin{enumerate}[\rm (i)]
        \item for every $(x,y)\in[0,1)^2$, $\phi_n(c_{(x,y)}) = y$,
        \item $\phi_n:\Omega_n\to[0,1]$ is continuous,
        \item $\phi_n:\Omega_n\to[0,1]$ is onto,
        \item if $\beta$ denotes the positive root of the polynomial $x^2-nx-1$, then
            \begin{align*}
            \phi_n(\sigma^{\be_1}w) &= \phi_n(w),\\
            \phi_n(\sigma^{\be_2}w) &= \phi_n(w) +\beta^{-1} \pmod 1.
            \end{align*}
    \end{enumerate}
\end{proposition}

\begin{proof}
    (i)
    Let $R_\alpha(x)=\{x+\alpha\}$ be the rotation by angle $\alpha$ on the interval $[0,1)$.
    If $\alpha$ is irrational, then for every $x\in[0,1)$ the sequence $(R^i_\alpha(x))_{i\in\Z}$ is
    uniformly distributed modulo 1
    \cite[Exercise~2.5]{zbMATH03440485}.
    Therefore, using Weyl's equidistribution theorem for Riemann-integrable functions
    \cite[Corollary~1.1]{zbMATH03440485},
    for every $(x,y)\in[0,1)^2$, we have
    {\allowdisplaybreaks
    \begin{align*}
        \phi_n(c_{(x,y)})
        &=
        \lim_{k\to\infty}
        \frac{1}{2k+1}
        \sum_{i=-k}^k
        \langle 
        \textstyle \frac{1}{n} d, 
        \sctop(c_{(x,y)}(i,0)) \rangle\\
        &=
        \lim_{k\to\infty}
        \frac{1}{2k+1}
        \sum_{i=-k}^k
        \langle 
        \textstyle \frac{1}{n} d, 
        \sctop(\sctile_n(x+i\beta^{-1}, y))
        \rangle\\
        &=
        \lim_{k\to\infty}
        \frac{1}{2k+1}
        \sum_{i=-k}^k
        \langle 
        \textstyle \frac{1}{n} d, 
        \Lambda_n(y, \{x+i\beta^{-1}\})
        \rangle\\
        &=
        \frac{1}{n}
        \lim_{k\to\infty}
        \frac{1}{2k+1}
        \sum_{i=-k}^k
            \left(\lfloor n y\rfloor 
                  + \I_{[1-\{ny\},1)}(\{\delta_{y} + \{x+i\beta^{-1}\}\})\right) 
                  \text{\qquad(Lemma~\ref{lem:inner-product-with-d-of-Lambda_n})} \\
        &=
        \frac{1}{n}
        \left(
        \lfloor ny\rfloor 
        +
        \lim_{k\to\infty}
        \frac{1}{2k+1}
        \sum_{i=-k}^k
             \I_{[1-\{ny\},1)}(R_{\beta^{-1}}^i(\delta_{y} + x))\right) \\
        &=
        \frac{1}{n}
        \left(
        \lfloor ny\rfloor 
        + \int_0^1 \I_{[1-\{ny\},1)}(t) dt \right) \text{\qquad(Weyl's equidistribution theorem)}\\
        &=
        \frac{1}{n}
        \left( \lfloor ny\rfloor + \{ny\} \right)
        = \frac{1}{n}(ny) = y.
    \end{align*}
    }

    (ii) 
    Now we want to show that the rule $\phi_n$ defines a continuous map $\Omega_n\to\torus$.
    Since $\Omega_n$ is minimal \cite{labbe_metallic_I_2025},
    we have that the orbit 
    $\shiftclosure{\{c_{(0,0)}\}}
    =\{\sigma^k c_{(0,0)} \mid k\in\Z^2\}
    =\{c_{\beta^{-1}k\pmod{\Z^2}} \mid k\in\Z^2\}$ 
    is a dense
    subset of $\Omega_n$.
    Therefore, $\{c_{(x,y)} \mid x,y\in[0,1)\}$ is dense in $\Omega_n$.
    Let $w\in\Omega_n$.
    There exists a sequence $(x^{(\ell)},y^{(\ell)})_{l\in\N}$
    with $x^{(\ell)},y^{(\ell)}\in[0,1)$
    such that 
    $w=\lim_{\ell\to\infty} c_{(x^{(\ell)},y^{(\ell)})}$.

    Notice that the limit
    $(x^{(\infty)},y^{(\infty)})=\lim_{\ell\to\infty} (x^{(\ell)},y^{(\ell)})\in[0,1]^2$
    exist.
    This essentially follows from \cite[Lemma 3.4]{MR4213162} allowing to define 
    another factor map, see Equation~\eqref{eq:factor-map-fn}.
    Indeed, suppose on the contrary that
    the sequence
    $(x^{(\ell)},y^{(\ell)})_{l\in\N}$
    has two distinct accumulation points $(p_1,q_1)$ and $(p_2,q_2)$.
    Recall that $\{\Int{\sctile_n^{-1}(t)}\}_{t\in\Tcal_n}$ is a topologicial
    partition of $\torus^2$.
    Since the orbits under the $\Z^2$-action $R_n$ are dense,
    there exists $(i,j)\in\Z^2$ such that
    $R_n^{(i,j)}(p_1,q_1)\in\Int{\sctile_n^{-1}(t_1)}$ and
    $R_n^{(i,j)}(p_2,q_2)\in\Int{\sctile_n^{-1}(t_2)}$
    where $t_1$ and $t_2$ are two distinct tiles in $\Tcal_n$.
    Therefore, for sufficiently large $\ell\in\N$, we have
    \begin{align*}
        w(i,j)&= c_{(x^{(\ell)},y^{(\ell)})}(i,j)= \sctile_n(R_n^{(i,j)}(p_1,q_1)) = t_1,\\
        w(i,j)&= c_{(x^{(\ell)},y^{(\ell)})}(i,j)= \sctile_n(R_n^{(i,j)}(p_2,q_2)) = t_2,
    \end{align*}
    which is a contradiction.

    We split the proof according to the behavior
    of $\lim_{\ell\to\infty}ny^{(\ell)}$,
    and more precisely if it converges to an integer and if so from above or
    from below (the fact that it converges from above or from below when it
    converges to an integer follows from the existence of the configuration
    $w$ because the boundary of the topological partition 
    $\{\Int{\sctile_n^{-1}(t)}\}_{t\in\Tcal_n}$
    contains the vertical and horizontal lines passing through integers points). 
    We proceed as above using Weyl equidistribution theorem.
    We have
    {\allowdisplaybreaks
    \begin{align*}
        \phi_n(w)
        &=
        \phi_n\left(\lim_{\ell\to\infty} c_{(x^{(\ell)},y^{(\ell)})}\right)\\
        &=
        \lim_{k\to\infty}
        \frac{1}{2k+1}
        \sum_{i=-k}^k
        \lim_{\ell\to\infty}
        \langle 
        \textstyle \frac{1}{n} d, 
        \sctop(c_{(x^{(\ell)},y^{(\ell)})}(i,0)) \rangle\\
        &=
        \frac{1}{n}
        \lim_{k\to\infty}
        \frac{1}{2k+1}
        \sum_{i=-k}^k
        \lim_{\ell\to\infty}
            \left(\lfloor n y^{(\ell)}\rfloor 
                  + \I_{[1-\{ny^{(\ell)}\},1)}(\{\delta_{y^{(\ell)}} + \{x^{(\ell)}+i\beta^{-1}\}\})\right) \\
        &=
        \frac{1}{n}
        \lim_{k\to\infty}
        \frac{1}{2k+1}
        \sum_{i=-k}^k
        \lim_{\ell\to\infty}
        \left(
        \lfloor ny^{(\ell)}\rfloor 
        +
             \I_{[1-\{ny^{(\ell)}\},1)}(R_{\beta^{-1}}^i(\delta_{y^{(\ell)}} + x^{(\ell)}))\right) \\
        &=\begin{cases}
        \frac{1}{n}
        \lim_{k\to\infty}
        \frac{1}{2k+1}
        \sum_{i=-k}^k
        \left(
        \lfloor ny^{(\infty)}\rfloor 
        +
             \I_{\varnothing}(R_{\beta^{-1}}^i(\delta_{y^{(\infty)}} + x^{(\infty)}))\right) 
             & \text{ if } \{ny^{(\ell)}\}\to 0, \\
        \frac{1}{n}
        \lim_{k\to\infty}
        \frac{1}{2k+1}
        \sum_{i=-k}^k
        \left(
        \lfloor ny^{(\infty)}\rfloor - 1
        +
             \I_{(0,1)}(R_{\beta^{-1}}^i(\delta_{y^{(\infty)}} + x^{(\infty)}))\right) 
             & \text{ if } \{ny^{(\ell)}\}\to 1, \\
        \frac{1}{n}
        \lim_{k\to\infty}
        \frac{1}{2k+1}
        \sum_{i=-k}^k
        \left(
        \lfloor ny^{(\infty)}\rfloor 
        +
             \I_{[1-\{ny^{(\infty)}\},1)}(R_{\beta^{-1}}^i(\delta_{y^{(\infty)}} + x^{(\infty)}))\right) 
             & \text{ if } \{ny^{(\ell)}\}\not\to 0,1, \\
        \end{cases}\\
        &=\begin{cases}
        \frac{1}{n}
        \left(
        \lfloor ny^{(\infty)}\rfloor 
            + \int_0^1 \I_{\varnothing}(t) dt\right)
             & \text{ if } \{ny^{(\ell)}\}\to 0, \\
        \frac{1}{n}
        \left(
        \lfloor ny^{(\infty)}\rfloor - 1
            + \int_0^1 \I_{(0,1)}(t) dt\right)
             & \text{ if } \{ny^{(\ell)}\}\to 1, \\
        \frac{1}{n}
        \left(
        \lfloor ny^{(\infty)}\rfloor 
            + \int_0^1 \I_{[1-\{ny^{(\infty)}\},1)}(t) dt
             \right) 
             & \text{ if } \{ny^{(\ell)}\}\not\to 0,1, \\
        \end{cases}\\
        &=\begin{cases}
        \frac{1}{n}
        \lfloor ny^{(\infty)}\rfloor + 0
             & \text{ if } \{ny^{(\ell)}\}\to 0, \\
        \frac{1}{n}
        \lfloor ny^{(\infty)}\rfloor -1 + 1
             & \text{ if } \{ny^{(\ell)}\}\to 1, \\
        \frac{1}{n}
        \left(
        \lfloor ny^{(\infty)}\rfloor 
        + \{ny^{(\infty)}\} \right) 
             & \text{ if } \{ny^{(\ell)}\}\not\to 0,1, \\
        \end{cases}\\
        &=y^{(\infty)}
        = \lim_{\ell\to\infty} y^{(\ell)}
        = \lim_{\ell\to\infty} \phi_n(c_{(x^{(\ell)},y^{(\ell)})}).
    \end{align*}
    }
    This shows that the rule $\phi_n$ defines a map $\Omega_n\to[0,1]$
    and that this map is continuous.

    (iii)
    If $y\in[0,1)$, then $y= \phi_n(c_{(0,y)})$.
    If $y=1$, then $y= \phi_n(\lim_{y\to1^-}c_{(0,y)})$.
    Thus, the map $\phi_n$ is onto.

    (iv)
    Since the map $\phi_n$ is continuous, we only need to show the equalities
    for a dense subset of $\Omega_n$. 
    Let $(x,y)\in[0,1)^2$.
    We have
    \[
        \phi_n(\sigma^{\be_1} c_{(x,y)}) 
        = \phi_n(c_{(\{x+\beta^{-1}\},y)}) 
        = y
        = \phi_n(c_{(x,y)}).
    \]
    Moreover, we have
    \[
        \phi_n(\sigma^{\be_2}c_{(x,y)}) 
        = \phi_n(c_{(x,\{y+\beta^{-1}\})})
        = \{y+\beta^{-1}\}
        = \phi_n(c_{(x,y)}) +\beta^{-1} \pmod 1.\qedhere
    \]
\end{proof}

Since $\phi_n(\sigma^{\be_1}w) = \phi_n(w)$ for every configuration $w\in\Omega_n$,
the factor map $\phi_n$ is far from being injective. We may improve this as follows.
We use the symmetry of the tiles in $\Tcal_n$ to define an involution on $\Omega_n$.
If $w\in\Omega_n$ is a configuration, then its image under a reflection by the
positive diagonal is the configuration $\widehat{w}\in\Omega_n$ defined as
\[
\begin{array}{rccl}
\widehat{w}:&\Z^2 & \to & \Tcal_n\\
    &(i,j) & \mapsto & \widehat{w_{j,i}}.
\end{array}
\]
This allows to define a map from the Wang shift to the $2$-dimensional torus
\begin{equation}\label{eq:Phi}
\begin{array}{rccl}
    \Phi_n:&\Omega_n & \to & \torus^2\\
    &w & \mapsto & (\phi_n(\widehat{w}), \phi_n(w)).
\end{array}
\end{equation}
The first coordinate $\phi_n(\widehat{w})$ computes the average of the inner product with $d$
of the right-hand labels of the Wang tiles in the column containing the origin
of the configuration $w$.  We show in the next theorem that $\Phi_n$
is a factor map.

\begin{THEOREMD}
    \MainTheoremD
\end{THEOREMD}

\begin{proof}
    From Proposition~\ref{prop:phi-is-a-factor-map},
    $\phi_n$ is continuous. Thus, $\Phi_n$ is also continuous.

    Let $(x,y)\in[0,1)^2$. 
    Using Lemma~\ref{lem:tilen-in-Tcal'n}, 
    for every $(i,j)\in\Z^2$, we have
    \[
        \widehat{c_{(x,y)}}(i,j)
        =\widehat{\sctile_n}(x{+}j\beta^{-1}, y{+}i\beta^{-1})
        =\sctile_n(y{+}i\beta^{-1}, x{+}j\beta^{-1})
        =c_{(y,x)}(i,j).
    \]
    Thus, the identity $\widehat{c_{(x,y)}}=c_{(y,x)}$ holds.
    We obtain
    \[
        (x,y) 
        = 
        (
        \phi_n(c_{(y,x)}),
        \phi_n(c_{(x,y)})
        )
        = 
        (
        \phi_n(\widehat{c_{(x,y)}}),
        \phi_n(c_{(x,y)})
        )
        =\Phi_n(c_{(x,y)}).
    \]
    Therefore, $\Phi_n$ is onto.

    Let $w\in\Omega_n$ be a configuration.
    Let $k=(k_1,k_2)\in\Z^2$.
    Using Proposition~\ref{prop:phi-is-a-factor-map},
    we have
    \begin{align*}
        \Phi_n\circ\sigma^k (w)
        &= \left(\phi_n(\widehat{\sigma^k w}), \phi_n(\sigma^k w)\right)\\
        &= \left(\phi_n(\sigma^{(k_2,k_1)} \widehat{w}), \phi_n(\sigma^{(k_1,k_2)} w)\right)\\
        &= \left(\phi_n(\widehat{w})+\beta^{-1}k_1, \phi_n(w)+\beta^{-1}k_2\right) \pmod{\Z^2}\\
        &= (\phi_n(\widehat{w}), \phi(w)) + \beta^{-1}(k_1,k_2) \pmod{\Z^2}\\
        &= \Phi_n(w) + \beta^{-1}k \pmod{\Z^2}\\
        &= R_n^k\circ \Phi_n (w).\qedhere
    \end{align*}
\end{proof}

\begin{corollary}\label{cor:Omegan_is_aperiodic}
    For every $n\geq 1$, $\Omega_n$ is aperiodic.
\end{corollary}

\begin{proof}
    By contradiction, suppose that $\Omega_n$ contains a periodic configuration $w$
    such that $\sigma^k(w)=w$ for some $k\in\Z^2\setminus\{(0,0)\}$.
    The image $\Phi_n(w)\in\torus^2$ must be a periodic point for the
    $\Z^2$-action $R_n$ because, using Theorem~\ref{thm:factor-map}, we have
    \[
        \Phi_n(w)
        =\Phi_n(\sigma^k(w))
        = R_n^k(\Phi_n(w))
        = R_n^k(\Phi_n(w)).
    \]
    The $\Z^2$-action $R_n$ has no periodic point,
    since the metallic mean $\beta$ is an irrational number.
    Thus, we must have $k=0$, which is a contradiction.
    The subshift $\Omega_n$ is nonempty.
    Thus, $\Omega_n$ is aperiodic.
\end{proof}

\begin{remark}
    Note that Corollary~\ref{cor:Omegan_is_aperiodic} can not be considered
    as a totally independent proof of aperiodicity of $\Omega_n$.
    Recall that aperiodicity of $\Omega_n$ was proved in \cite{labbe_metallic_I_2025}
    from the self-similarity of $\Omega_n$.
    Indeed, Corollary~\ref{cor:Omegan_is_aperiodic} 
    uses Theorem~\ref{thm:factor-map} which
    depends on Proposition~\ref{prop:phi-is-a-factor-map}.
    In the proof of Proposition~\ref{prop:phi-is-a-factor-map},
    we use the minimality of $\Omega_n$ which was proved in
    \cite{labbe_metallic_I_2025}
    and deduced from its self-similarity.
\end{remark}

    In other words, the following question remains open.

\begin{question}
    Can the aperiodicity of $\Omega_n$ be proved independently of its self-similarity?
\end{question}

\section{The factor map is an isomorphism (mod 0)}
\label{sec:isomorphism}

The goal of this section is to show more properties of the factor map
$\Phi_n:\Omega_n\to\torus^2$ introduced in the previous section.
Based on the approach presented in \cite{MR4213162},
we prove
Theorem~\ref{thm:Markov-partition} and
Theorem~\ref{thm:maximal-equicontinuous-factor}.

Let $n\geq1$ be an integer.
We consider
the continuous $\Z^2$-action $R_n$ defined
on $\torus^2=\R^2/\Z^2$ by
\[
\begin{array}{rccl}
    R_n&:\Z^2\times\torus^2 & \to & \torus^2\\
    &(\bn,\bx) & \mapsto & R_n^\bn(\bx):=\bx + \beta\bn
\end{array}
\]
where $\beta=\frac{n+\sqrt{n^2+4}}{2}$
is the positive root of the polynomial $x^2-nx-1$.
We say that $R_n$ is a \defn{toral $\Z^2$-rotation}
and it defines a dynamical system that we denote
$\dynsys{\Z^2}{R_n}{\torus^2}$.
In this section, we encode this dynamical system symbolically
using a partition associated with the Wang tiles $\Tcal_n$.

\begin{figure}
\begin{center}
\begin{tabular}{ccc}
        $\north_1$ & $\east_1$ & $\north_1\wedge\east_1$\\
 \includegraphics[scale=.80]{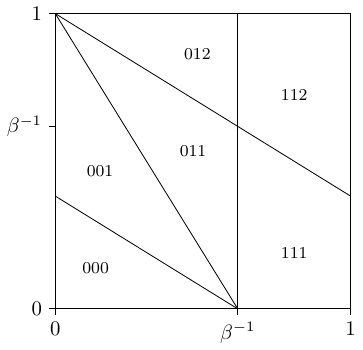}
&\includegraphics[scale=.80]{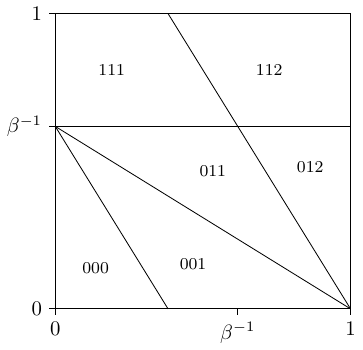}
&\includegraphics[scale=.80]{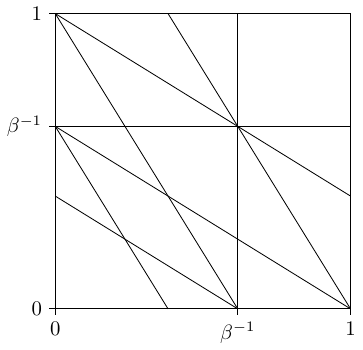}\\
        $\south_1$ & $\west_1$ & $\south_1\wedge\west_1$\\
 \includegraphics[scale=.80]{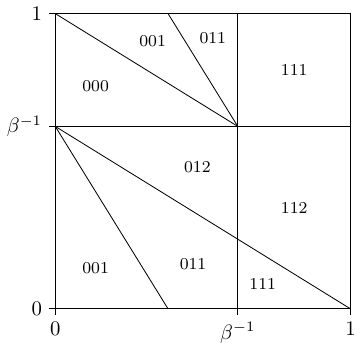}
&\includegraphics[scale=.80]{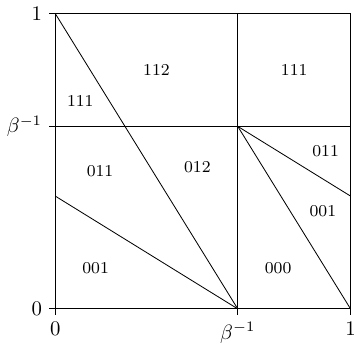}
&\includegraphics[scale=.80]{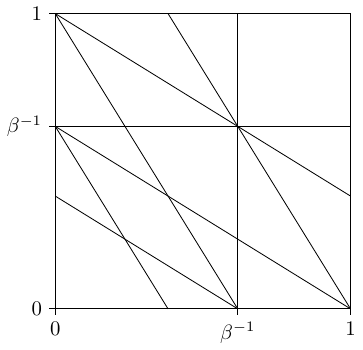}\\
    \end{tabular}
\end{center}
\caption{The partitions $\north_1$, $\east_1$, $\south_1$ and $\west_1$.}
\label{fig:partition-n1}
\end{figure}

\begin{figure}
\begin{center}
    \begin{tabular}{ccc}
        $\north_2$ & $\east_2$ & $\north_2\wedge\east_2$\\
 \includegraphics[scale=.80]{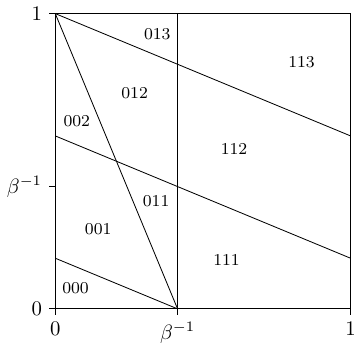}
&\includegraphics[scale=.80]{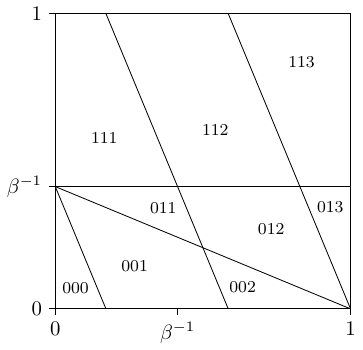}
&\includegraphics[scale=.80]{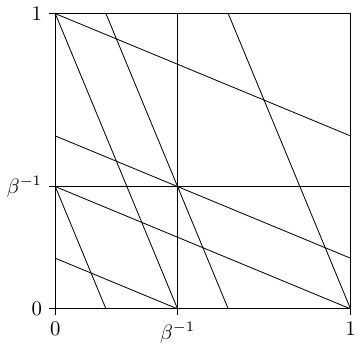}\\
        $\south_2$ & $\west_2$ & $\south_2\wedge\west_2$\\
 \includegraphics[scale=.80]{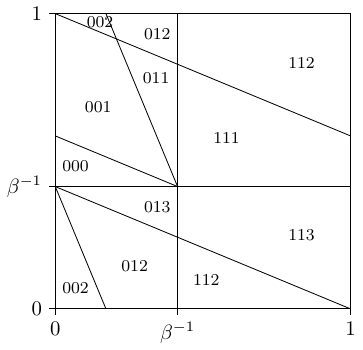}
&\includegraphics[scale=.80]{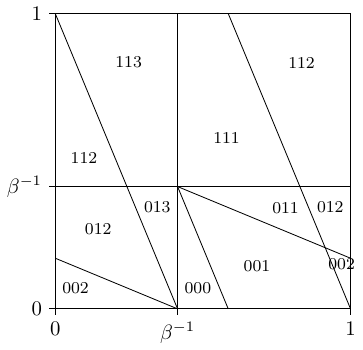}
&\includegraphics[scale=.80]{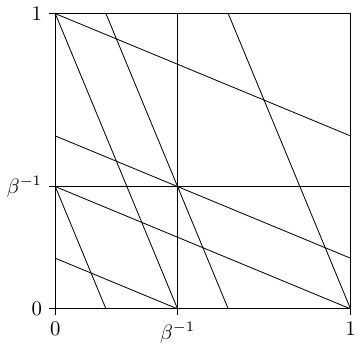}\\
    \end{tabular}
\end{center}
\caption{The partitions $\north_2$, $\east_2$, $\south_2$ and $\west_2$.}
\label{fig:partition-n2}
\end{figure}

\begin{figure}[h]
\begin{center}
    \begin{tabular}{ccc}
        $\north_3$ & $\east_3$ & $\north_3\wedge\east_3$\\
 \includegraphics[scale=.74]{SAGEOUTPUT/W3_partition_top_with_axis_labels.pdf}
&\includegraphics[scale=.74]{SAGEOUTPUT/W3_partition_right_with_axis_labels.pdf}
&\includegraphics[scale=.74]{SAGEOUTPUT/W3_partition_top_right_with_axis_labels.pdf}\\
        $\south_3$ & $\west_3$ & $\south_3\wedge\west_3$\\
 \includegraphics[scale=.74]{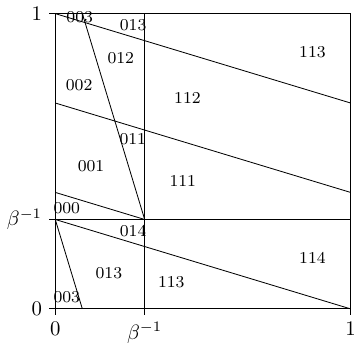}
&\includegraphics[scale=.74]{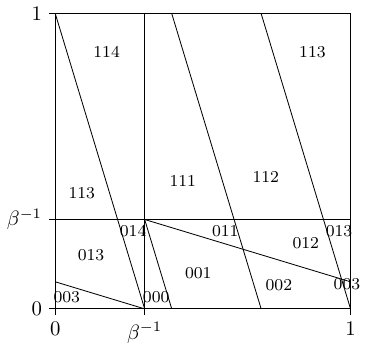}
&\includegraphics[scale=.74]{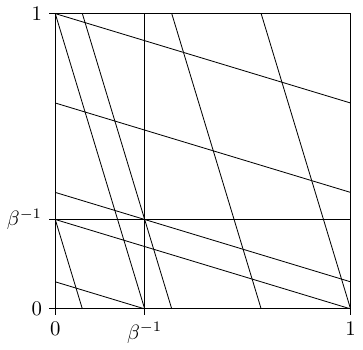}\\
    \end{tabular}
\caption{The partitions $\north_3$, $\east_3$, $\south_3$ and $\west_3$.}
\label{fig:partition-n3}
\end{center}
\end{figure}

\begin{figure}[h]
\begin{center}
    \begin{tabular}{ccc}
        $\north_4$ & $\east_4$ & $\north_4\wedge\east_4$\\
 \includegraphics[scale=.74]{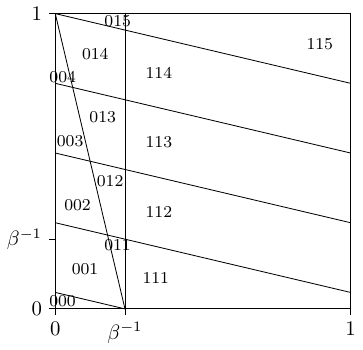}
&\includegraphics[scale=.74]{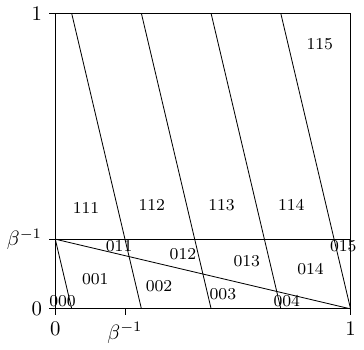}
&\includegraphics[scale=.74]{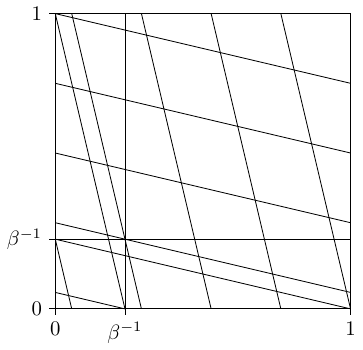}\\
        $\south_4$ & $\west_4$ & $\south_4\wedge\west_4$\\
 \includegraphics[scale=.74]{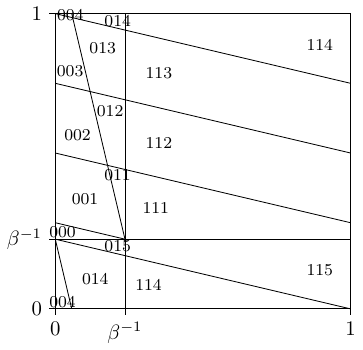}
&\includegraphics[scale=.74]{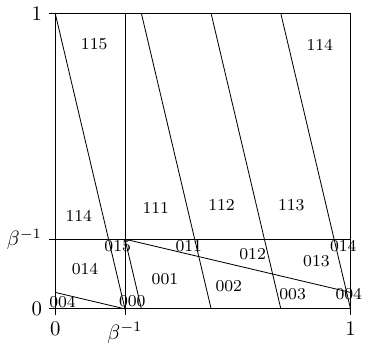}
&\includegraphics[scale=.74]{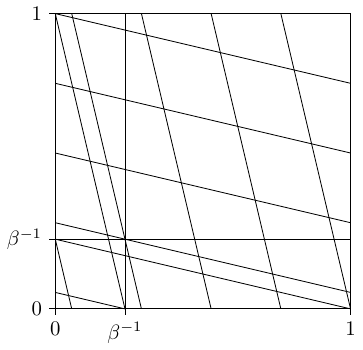}\\
    \end{tabular}
\caption{The partitions $\north_4$, $\east_4$, $\south_4$ and $\west_4$.}
\label{fig:partition-n4}
\end{center}
\end{figure}

Recall that
\[
\begin{array}{rcl}
    \Lambda_n:[0,1)^2 & \to & \Z^3\\
    (x,y) & \mapsto &
             \left(
             \begin{array}{r}
                 \lfloor y+\beta^*              +1\rfloor\\
                 \lfloor \beta^{-1}x + y+\beta^*+1\rfloor\\
                 \lfloor \beta x + y+\beta^*    +1\rfloor
             \end{array}
             \right).
\end{array}
\]
From Lemma~\ref{lem:ceillambdan-nondecreasing}, we have in fact that
$\Lambda_n$ is a map $[0,1)^2 \to V_n$.
Therefore, 
\[
     \east_n = \{\Lambda_n^{-1}(v)\colon v\in V_n\}
\]
is a partition of $[0,1)^2$. Its symmetric image is
\[
    \north_n = \left\{\eta\circ\Lambda_n^{-1}(v)\colon v\in V_n\right\}
\]
which is another partition of $[0,1)^2$, where $\eta:(x,y)\mapsto(y,x)$.
Also, we let 
\begin{align*}
    \west_n  &= R_n^{\be_1}(\east_n),\\
    \south_n &= R_n^{\be_2}(\north_n)
\end{align*}
where $\be_1=(1,0)$ and $\be_2=(0,1)$.
These partitions are illustrated for $n=1,2,3,4$ in
Figure~\ref{fig:partition-n1},
Figure~\ref{fig:partition-n2},
Figure~\ref{fig:partition-n3} and
Figure~\ref{fig:partition-n4}.
We may observe in these figures a nice property of the partitions:
$\east_n\wedge\north_n$
is the same partition (with different indices) as
$\west_n\wedge\south_n$
(this is related to the fact that the set of Wang tiles $\Tcal_n$ is both
NE-deterministic and SW-deterministic, see
Theorem~\ref{thm:NE-SW-deterministic}).

We now want to construct the refined partition
     $\east_n \wedge
    \north_n \wedge
    \west_n  \wedge
    \south_n$ whose atoms are defined as follows.
For each $(v_1,v_2,v_3,v_4)\in (V_n)^4$,
we define the interior of the intersection
\[
    P_{(v_1,v_2,v_3,v_4)} = 
                \Int{
                  \Lambda_n^{-1}(v_1)
                  \cap \eta\circ\Lambda_n^{-1}(v_2)
                  \cap R^{\be_1}(\Lambda_n^{-1}(v_3))
                  \cap R^{\be_2}(\eta\circ\Lambda_n^{-1}(v_4))}.
\]
It follows from Proposition~\ref{prop:Tcaln_is_image_of_sctilexy}
that the quadruples $\tau$ for which $P_\tau$ has nonempty interior
define a set which is equal to the set of Wang tiles $\Tcal_n$:
\[
    \Tcal_n = \left\{\tau\in  (V_n)^4
                  \mid P_\tau \neq\varnothing\right\}.
\]
Recall that, for some finite set $A$,
a \defn{topological partition} of a compact metric space $M$ is a finite
collection $\{P_a\}_{a\in A}$ of disjoint open sets $P_a\subset M$
such that $M = \bigcup_{a\in A} \overline{P_a}$.
Naturally, the set $\Tcal_n$ defines a topological partition
\[
    \Pcal_n = \{P_\tau\}_{\tau\in\Tcal_n}
\]
of $\R^2/\Z^2$ which is the refinement of the four partitions
$\east_n$ (the right color), 
$\north_n$ (the top color),
$\west_n$ (the left color) and 
$\south_n$ (the bottom color).

\subsection{Symbolic dynamical system $\Xcal_{\Pcal_n,R_n}$}

We now define the symbolic dynamical system associated with the toral $\Z^2$-rotation
$R_n$ generated by the partition $\Pcal_n$.
We adapt \cite{MR1369092} to the 2-dimensional setting as it was done in
\cite{MR3525488} and \cite{MR4213162}.

If $S\subset\Z^2$ is a finite set,
we say that a pattern $w\in\Acal^S$
is \defn{allowed} for $\Pcal_n,R_n$ if
\begin{equation}\label{eq:allowed-if-nonempty}
    \bigcap_{\bk\in S} R_n^{-\bk}(P_{w_\bk}) \neq \varnothing.
\end{equation}
Let $\Lcal_{\Pcal_n,R_n}$ be the collection of all allowed patterns for $\Pcal_n,R_n$.
The set $\Lcal_{\Pcal_n,R_n}$ is the language of a subshift 
$\Xcal_{\Pcal_n,R_n}\subseteq\Acal^{\Z^2}$ defined as follows,
see \cite[Prop.~9.2.4]{MR3525488},
\[
    \Xcal_{\Pcal_n,R_n} = 
    \{x\in\Acal^{\Z^2} \mid \pi_S\circ\sigma^\bn(x)\in\Lcal_{\Pcal_n,R_n}
    \text{ for every } \bn\in\Z^2 \text{ and finite subset } S\subset\Z^2\}.
\]
We say that $\Xcal_{\Pcal_n,R_n}$ is the \defn{symbolic dynamical
system} corresponding to $\Pcal_n,R_n$.

For each $w\in\Xcal_{\Pcal_n,R_n}\subset\Acal^{\Z^2}$ and $m\geq 0$ there is a corresponding nonempty open set
\[
    D_m(w) = \bigcap_{\Vert\bk\Vert\leq m} R_n^{-\bk}(P_{w_\bk}) \subseteq \torus^2.
\]
The closures $\overline{D}_m(w)$ of these sets are compact
and decrease with $m$, so that
$\overline{D}_0(w)\supseteq
\overline{D}_1(w)\supseteq
\overline{D}_2(w)\supseteq
\dots$.
It follows that $\cap_{m=0}^{\infty}\overline{D}_m(w)\neq\varnothing$.
In order for points in
$\Xcal_{\Pcal_n,R_n}$
to correspond to points in $\torus^2$, this intersection should contain only one point.
This leads to the following definition.
A topological partition $\Pcal_n$ of $\torus^2$ \defn{gives a symbolic representation}
of 
$\dynsys{\Z^2}{R_n}{\torus^2}$
    if for every $w\in\Xcal_{\Pcal_n,R_n}$ the intersection
$\cap_{m=0}^{\infty}\overline{D}_m(w)$ consists of exactly one
point $\bx\in\torus^2$.
We call $w$ a \defn{symbolic representation of $\bx$}.

Markov partitions were originally defined for one-dimensional dynamical
systems $\dynsys{\Z}{T}{\torus^2}$
and were extended to $\Z^d$-actions by automorphisms of
compact Abelian group in \cite{MR1632169}.
Following \cite{MR4213162,MR4347332},
we use the same terminology 
and extend the definition proposed in \cite[\S 6.5]{MR1369092}
for dynamical systems defined by higher-dimensional actions by rotations.

\begin{definition}\label{def:Markov}
A topological partition $\Pcal$ of $\torus^2$ is a \defn{Markov partition} for
$\dynsys{\Z^2}{R}{\torus^2}$  if
\begin{itemize}
    \item $\Pcal$ gives a symbolic representation of $\dynsys{\Z^2}{R}{\torus^2}$ and 
    \item $\Xcal_{\Pcal,R}$ is a shift of finite type (SFT).
\end{itemize}
\end{definition}

\subsection{Proofs of main results}

First, we have the following result.

\begin{lemma}\label{lem:minimal-aperiodic}
    The dynamical system $\dynsys{\Z^2}{\sigma}{\Xcal_{\Pcal_n,R_n}}$ is minimal
    and $\Xcal_{\Pcal_n,R_n}$ is aperiodic.
\end{lemma}

\begin{proof}
    Since $R_n^{\be_1}$ and $R_n^{\be_2}$ are linearly independent irrational
    rotations on $\R^2/\Z^2$, we have
    that $R_n$ is a free $\Z^2$-action.
    Thus, from \cite[Lemma 5.2]{MR4213162},
    $\Xcal_{\Pcal_n,R_n}$ is minimal and aperiodic.
\end{proof}

Each atom of the partition $\Pcal_n$ is invariant only under the trivial translation.
Therefore, from \cite[Lemma 3.4]{MR4213162},
$\Pcal_n$ gives a symbolic representation of the dynamical system 
$\dynsys{\Z^2}{R_n}{\torus^2}$.
Thus, we can define the following function:
\begin{equation}\label{eq:factor-map-fn}
f_n:\Xcal_{\Pcal_n,R_n}\to\torus^2
\end{equation}
be such that $f_n(w)$ is the unique point in the intersection
$\cap_{m=0}^{\infty}\overline{D}_m(w)$.

\begin{proposition}\label{prop:f_n-is-factor-map}
    Let $n\geq1$ be an integer.
    The map $f_n:\Xcal_{\Pcal_n,R_n}\to\torus^2$ is a factor map
    satisfying
    \[
        f_n\circ\sigma^k
        = R_n^k\circ f_n 
    \]
    for every $k\in\Z^2$.
\end{proposition}

\begin{proof}
    The result is an application of Proposition 5.1 from \cite{MR4213162}.
\end{proof}

From the minimality of the Wang shift $\Omega_n$ proved separately in
\cite{labbe_metallic_I_2025},
we may now prove Theorem~\ref{thm:Markov-partition}
using the same method as in \cite{MR4213162}.

\begin{THEOREMF}
    \MainTheoremE
\end{THEOREMF}

\begin{proof}
    From Proposition 8.1 in \cite{MR4213162},
    we have that
    $\Xcal_{\Pcal_n,R_n}\subseteq\Omega_n$ for every integer $n\geq 1$.
    It was proved in \cite{labbe_metallic_I_2025}
    that the Wang shift $\Omega_n$ is minimal
    for every integer $n\geq1$.
    Thus, $\Xcal_{\Pcal_n,R_n}=\Omega_n$.

    Each atom of the partition $\Pcal_n$ is invariant only under the trivial translation.
    Therefore, from \cite[Lemma 3.4]{MR4213162},
    $\Pcal_n$ gives a symbolic representation of $\dynsys{\Z^2}{R_n}{\torus^2}$.
    Since $\Xcal_{\Pcal_n,R_n}=\Omega_n$ is a shift of finite type, 
    we conclude that the partition $\Pcal_n$
    is a Markov partition for the dynamical system
    $\dynsys{\Z^2}{R_n}{\torus^2}$.
\end{proof}

In fact, we can show that the factor map $f_n$ is equal to the map $\Phi_n$
explicitly defined in Section~\ref{sec:explicit-factor-map}
from the average of the labels of Wang tiles on the row and column
containing the origin.
It follows from the next lemma.

\begin{lemma}\label{lem:f-cyx=xy}
    For every $(x,y)\in[0,1)^2$, we have $f_n(c_{(x,y)})=(x,y)$.
\end{lemma}

\begin{proof}
    Let $v_1,v_2,v_3,v_4\in V_n$. Observe that
    \begin{align*}
        \sctile_n^{-1}(v_1,v_2,v_3,v_4)
        &\subseteq
        \Lambda_n^{-1}(v_1)
                  \cap \eta\circ\Lambda_n^{-1}(v_2)
                  \cap R^{\be_1}(\Lambda_n^{-1}(v_3))
                  \cap R^{\be_2}(\eta\circ\Lambda_n^{-1}(v_4))\\
        &\subset
        \overline{\Lambda_n^{-1}(v_1)
                  \cap \eta\circ\Lambda_n^{-1}(v_2)
                  \cap R^{\be_1}(\Lambda_n^{-1}(v_3))
                  \cap R^{\be_2}(\eta\circ\Lambda_n^{-1}(v_4))}\\
        &=
        \overline{P_{(v_1,v_2,v_3,v_4)}}.
    \end{align*}
    For every $k\in\Z^2$, we have
    \[
        c_{(x,y)}(k) = \sctile_n\circ R_n^{k}(x,y),
    \]
    so that
    \[
        (x,y)
        \in R_n^{-k}\circ\sctile_n^{-1}(c_{(x,y)}(k))
        \subset
        R_n^{-k}(\overline{P_{c_{(x,y)}(k)}}).
    \]
    Therefore, for every $m\in\N$, we have
    \[
        (x,y)
        \in 
        \bigcap_{\Vert k\Vert\leq m} R_n^{-k}(\overline{P_{c_{(x,y)}(k)}}) 
        =
        \overline{D_m}(c_{(x,y)}).
    \]
    Since $\Pcal_n$ gives a symbolic representation of the dynamical system
    $\dynsys{\Z^2}{R_n}{\torus^2}$, we have that
    $\cap_{m=0}^{\infty}\overline{D}_m(c_{(x,y)})$
    is a singleton
    and 
    \[
        \cap_{m=0}^{\infty}\overline{D}_m(c_{(x,y)})=\{(x,y)\}.
    \]
    Therefore, $f(c_{(x,y)})=(x,y)$.
\end{proof}

\begin{proposition}\label{prop:fn=Phin}
    The factor map $f_n:\Omega_n\to\torus^2$ is equal to the factor map
    $\Phi_n:\Omega_n\to\torus^2$ explicitly defined in Equation~\eqref{eq:Phi}:
    \[
        f_n=\Phi_n.
    \]
\end{proposition}

\begin{proof}
    From Lemma~\ref{lem:f-cyx=xy}, we have $ f_n(c_{(0,0)}) =(0,0)$.
    Also, observe that the configuration $c_{(0,0)}$ is symmetric: $\widehat{c_{(0,0)}}=c_{(0,0)}$.
    Thus, we have 
    \[
        \Phi_n(c_{(0,0)})
             =(\phi_n(\widehat{c_{(0,0)}}), \phi_n(c_{(0,0)}))
             =(\phi_n(c_{(0,0)}), \phi_n(c_{(0,0)}))=(0,0).
    \]
    Let $w\in\Omega_n$ be any configuration.
    Since $\Omega_n$ is minimal \cite{labbe_metallic_I_2025},
    there exists a sequence $(k_\ell)_{\ell\in\N}$
    such that $k_\ell\in\Z^2$ such that
    $w=\lim_{\ell\to\infty}\sigma^{k_\ell}(c_{(0,0)})$.
    From Proposition~\ref{prop:f_n-is-factor-map}
    and Theorem~\ref{thm:factor-map},
    $f_n$ and $\Phi_n$ are factor maps
    commuting the shift map with the $\Z^2$-action $R_n$ on the torus $\torus^2$.
    Thus, we obtain
    \begin{align*}
        \Phi_n(w)
        &= \Phi_n\left(\lim_{\ell\to\infty}\sigma^{k_\ell}(c_{(0,0)})\right)\\
        &= \lim_{\ell\to\infty}\Phi_n\circ\sigma^{k_\ell}(c_{(0,0)})\\
        &= \lim_{\ell\to\infty}R_n^{k_\ell}\circ\Phi_n(c_{(0,0)})\\
        &= \lim_{\ell\to\infty}R_n^{k_\ell}\left((0,0)\right)\\
        &= \lim_{\ell\to\infty}R_n^{k_\ell}\circ f_n(c_{(0,0)})\\
        &= \lim_{\ell\to\infty}f_n\circ\sigma^{k_\ell}(c_{(0,0)})\\
        &= f_n\left(\lim_{\ell\to\infty}\sigma^{k_\ell}(c_{(0,0)})\right)
         = f_n(w).\qedhere
    \end{align*}
\end{proof}

The factor map $\Phi_n$ between the dynamical system
$\dynsys{\Z^2}{\sigma}{\Omega_n}$ and the $\Z^2$-action $R_n$ on the torus
$\torus^2$ satisfies additional properties.
In particular, $\Phi_n$ is an isomorphism of measure-preserving dynamical systems.
Their proofs follow the structure of similar results proved in \cite{MR4213162}
for Jeandel--Rao tilings.

\begin{THEOREME}
    \MainTheoremF
\end{THEOREME}

\begin{proof}
    From Theorem~\ref{thm:Markov-partition}, we have $\Xcal_{\Pcal_n,R_n}=\Omega_n$.

    (i)
    From Proposition~\ref{prop:f_n-is-factor-map},
    the factor map $f_n:\Xcal_{\Pcal_n,R_n}\to\torus^2$
    commutes the actions
    $\dynsys{\Z^2}{\sigma}{\Xcal_{\Pcal_n,R_n}}$
    and $\dynsys{\Z^2}{R_n}{\torus^2}$.
    From \cite[Proposition 5.1]{MR4213162}, 
    $f_n$ is one-to-one on 
    $f_n^{-1}(\torus^2\setminus\Delta_{\Pcal_n,R_n})$ where
    \begin{equation*}\label{eq:boundaries}
        \Delta_{\Pcal_n,R_n}:=\bigcup_{\bk\in\Z^2}R_n^\bk
            \left(\bigcup_{\tau\in\Tcal_n}\partial P_\tau\right)
            \subset \torus^2
    \end{equation*}
    is the set of points whose orbit under the $\Z^2$-action $R_n$ intersect
    the boundary of the topological partition $\Pcal_n=\{P_\tau\}_{\tau\in\Tcal_n}$.
    From \cite[Corollary 5.3]{MR4213162} (which is a consequence of \cite[Lemma 3.11]{MR3381481}),
    $\dynsys{\Z^2}{R_n}{\torus^2}$ is the maximal equicontinuous factor of
    $\dynsys{\Z^2}{\sigma}{\Xcal_{\Pcal_n,R_n}}$.

    (ii)
    We have that $\{y\in
    \torus^2:\card(f_n^{-1}(y))=1\}=\torus^2\setminus\Delta_{\Pcal_n,R_n}$ is a
    countable intersection of open sets and is dense in $\torus^2$.
    Thus, it is a $G_\delta$-dense set in $\torus^2$.
    Therefore, the factor map $f_n:\Xcal_{\Pcal_n,R_n}\to\torus^2$ is almost one-to-one.
    From Proposition~\ref{prop:fn=Phin}, we have $f_n=\Phi_n$.

Suppose that $\bx\in\Delta_{\Pcal_n,R_n}$.
We have $\card(f_n^{-1}(\bx))\geq2$.
If $\card(f_n^{-1}(\bx))>2$, then we may show that
    there exists $\bn\in\Z^2$ such that $\bx=R_n^\bn(\zero)$.
    If $\bx=R_n^{\bn}(\zero)$ for some $\bn\in\Z^2$, then
the set $f_n^{-1}(\bx)$ contains 8 different configurations
    of the form $\lim_{\varepsilon\to0}c_{\varepsilon\bv}$ for some $\bv\in\R^2\setminus\Theta^{\Pcal_n}$
    where $\Theta^{\Pcal_n}=\R\cdot\{(1,0),(0,1),(1,-\beta),(1,\beta^*)\}$.
If $\bx\in\Delta_{\Pcal_n,R_n}$ but not
    in the orbit of $\zero$ under $R_n$,
    then $\card(f_n^{-1}(\bx))=2$.
    We conclude that
    $\{\card(f_n^{-1}(\bx))\mid\bx\in\torus^2\}=\{1,2,8\}$.

    (iii)
    The dynamical system $\dynsys{\Z^2}{R_n}{\torus^2}$ is minimal.
    We have that $\lambda(\partial P)=0$ for each atom $P\in\Pcal_n$
    where $\lambda$ is the Haar measure on $\torus^2$.
    The partition $\Pcal_n$ gives a symbolic representation of the dynamical system
    $\dynsys{\Z^2}{R_n}{\torus^2}$.
    Thus, from \cite[Proposition 6.1]{MR4213162},
    the dynamical system
    $\dynsys{\Z^2}{\sigma}{\Xcal_{\Pcal_n,R_n}}$
    is uniquely ergodic.

    (iv)
    Since the dynamical system
    $\dynsys{\Z^2}{\sigma}{\Xcal_{\Pcal_n,R_n}}$ is uniquely ergodic,
    it admits a unique
    shift-invariant probability measure $\nu$ on $\Omega_n$.
    From \cite[Proposition 6.1]{MR4213162},
    the measure-preserving dynamical system $(\Omega_n,\Z^2,\sigma,\nu)$
    is isomorphic to $(\torus^2,\Z^2,R_n,\lambda)$ where 
    $\lambda$ is the Haar measure on $\torus^2$.
\end{proof}

\section{Renormalization and Rauzy induction of $\Z^2$-rotations}
\label{sec:renormalization-rauzy}

Another consequence of Theorem~\ref{thm:Markov-partition} is that the symbolic
dynamical system $\Xcal_{\Pcal_n,R_n}$ is self-similar because this was proved
in \cite{labbe_metallic_I_2025} for the Wang shift $\Omega_n$.
The Rauzy induction of polygonal partitions and of toral $\Z^2$-rotations
defined in \cite{MR4347332} can be used to compute the self-similarity of the
symbolic dynamical system $\Xcal_{\Pcal_n,R_n}$. We illustrate below how this
can be done for a fixed value of an integer $n\geq1$.

For some postive integer $n\geq1$,
we define the positive root $\beta$ of the polynomial $x^2-nx-1$. Computations
will be done in the number field generated by this root.
We perform the computations below with $n=3$, but it works with other integers.
For instance, the computation of the self-similarity for $n=7$ from the Rauzy induction is
done in about 200 seconds on a recent laptop.
\begin{sagecommandline}
sage: n = 3  # try with another integer
sage: x = polygen(QQ, "x")
sage: K.<beta> = NumberField(x^2 - n*x - 1, embedding=RR(n))
sage: beta.n()
3.30277563773199
\end{sagecommandline}
We define a function that computes the atoms $\Lambda_n^{-1}(v)$ for every $v\in V_n$.
Note that in SageMath, an entry equal to \texttt{[-1,7,3,4]} represents the inequality
$7x_1+3x_2+4x_3\geq1$.
\begin{sagecommandline}
sage: unit_square_ieqs = [[0, 1, 0], [0, 0, 1], [1, -1, 0], [1, 0, -1]]
sage: def Lambda_inv(a,b,c):
....:     ieqs = list(unit_square_ieqs)
....:     ieqs.extend([[-1/beta+1-a, 0, 1], [a+1/beta, 0, -1]])
....:     ieqs.extend([[-1/beta+1-b, 1/beta, 1], [b+1/beta, -1/beta, -1]])
....:     ieqs.extend([[-1/beta+1-c, beta, 1], [c+1/beta, -beta, -1]])
....:     return Polyhedron(ieqs=ieqs)
\end{sagecommandline}
We define the set $V_n$ and we check that the sum of the area of the polygons
$\{\Lambda_n^{-1}(v)\}_{v\in V_n}$ is 1.
\begin{sagecommandline}
sage: Vn = [(a,b,c) for a in range(2) for b in range(2) for c in range(n+2) if a<=b<=c]
sage: Vn
[(0, 0, 0), (0, 0, 1), ..., (1, 1, ...)]
sage: assert sum(Lambda_inv(*v).volume() for v in Vn) == 1
sage: Lambda_inv(0,0,n+1).volume()      # one of the atom has empty interior
0
\end{sagecommandline}
For readability reason, we define a map which concatenates the entries of a
vector into a string.
\begin{sagecommandline}
sage: def vector_to_str(v):
....:     return "".join(str(a) for a in v)
sage: vector_to_str((0,1,4))        # for example
'014'
\end{sagecommandline}
We define the $\Z^2$-action $R_n$ on $\R^2/\Z^2$ as two polyhedron exchange
transformations on the unit square.
\begin{sagecommandline}
sage: lattice_base = identity_matrix(2)
sage: from slabbe import PolyhedronExchangeTransformation as PET
sage: Re1 = PET.toral_translation(lattice_base, vector((1/beta,0)))
sage: Re2 = PET.toral_translation(lattice_base, vector((0,1/beta)))
\end{sagecommandline}
We construct the $\east_n$ partition (ignoring the atom with empty interior)
and the three other partitions from it.
\begin{sagecommandline}
sage: from slabbe import PolyhedronPartition
sage: EAST = PolyhedronPartition({vector_to_str(v):Lambda_inv(*v) for v in Vn 
....:                                if Lambda_inv(*v).volume() > 0})
sage: M = matrix(K, 2, (0,1,1,0))
sage: NORTH = EAST.apply_linear_map(M)
sage: WEST = Re1(EAST)
sage: SOUTH = Re2(NORTH)
sage: G = graphics_array([EAST.plot(),NORTH.plot(), SOUTH.plot(),WEST.plot()])
sage: G.show(figsize=10)
\end{sagecommandline}
\begin{center}
\sageplot[width=.95\linewidth][pdf]{G,figsize=10}
\end{center}
We compute the refinement 
of the $\east_n$ and $\north_n$ partitions
and of the $\west_n$ and $\south_n$ partitions.  
\begin{sagecommandline}
sage: PEN,dEN = EAST.refinement(NORTH, certificate=True)
sage: PWS,dWS = WEST.refinement(SOUTH, certificate=True)
sage: G = graphics_array([PEN.plot(),PWS.plot()])
sage: G.show(figsize=5)
\end{sagecommandline}
\begin{center}
\sageplot[width=.50\linewidth][pdf]{G,figsize=7}
\end{center}
In general, we would need to compute the refinement of the two partitions.
But here, since they are equal up to relabeling,
we may take one as the refinement and compute the
bijection of the labels between them.
\begin{sagecommandline}
sage: PWS.is_equal_up_to_relabeling(PEN)
True
sage: P = PEN           # faster than P = PEN.refinement(PWS)
sage: bijection = P.keys_permutation(PWS)
sage: bijection[9]     # for example
16
\end{sagecommandline}
We compute the set of Wang tiles defined by the refinement of the four partitions
$\east_n$, $\north_n$, $\west_n$ and $\south_n$:
\begin{sagecommandline}
sage: from slabbe import WangTileSet
sage: tiles = [dEN[i]+dWS[bijection[i]] for i in sorted(dEN)]
sage: T3 = WangTileSet(tiles)
sage: t = T3.tikz(ncolumns=10, scale=1.2)
\end{sagecommandline}
\begin{center}
\sageplot[width=.6\linewidth][pdf]{t}
\end{center}
We perform the Rauzy induction on the square window $[0,\beta^{-1}]\times[0,\beta^{-1}]$
using the algorithms 
\texttt{induced\_partition}
and
\texttt{induced\_transformation}
defined in \cite{MR4347332}.
First, we perform the induction on the domain restricted to the inequality $x\leq\beta^{-1}$.
\begin{sagecommandline}
sage: x_le_beta_inv = [1/beta,-1,0]
sage: P1,s1 = Re1.induced_partition(x_le_beta_inv, P, substitution_type="row")
sage: R1e1,_ = Re1.induced_transformation(x_le_beta_inv)
sage: R1e2,_ = Re2.induced_transformation(x_le_beta_inv)
\end{sagecommandline}
Secondly, we perform the induction on the domain restricted to the inequality $y\leq\beta^{-1}$.
\begin{sagecommandline}
sage: y_le_beta_inv = [1/beta,0,-1]
sage: P2,s2 = Re2.induced_partition(y_le_beta_inv, P1, substitution_type="column")
sage: R2e1,_ = R1e1.induced_transformation(y_le_beta_inv)
sage: R2e2,_ = R1e2.induced_transformation(y_le_beta_inv)
\end{sagecommandline}
We rescale the induced partition by the factor $-\beta$ and translate it back
to the unit square in the positive quadrant.
Then we apply each generator of the $\Z^2$-action once on the rescaled induced partition.
\begin{sagecommandline}
sage: P2_scaled = (-beta * P2).translate((1,1))
sage: P3 = Re2(Re1(P2_scaled))
sage: G = graphics_array([P2_scaled.plot(), P3.plot()])
sage: G.show(figsize=5)
\end{sagecommandline}
\begin{center}
\sageplot[width=.5\linewidth][pdf]{G,figsize=7}
\end{center}
We check that the resulting partition is equal to the initial partition.
We check that the induced action is equal to the initial action.
\begin{sagecommandline}
sage: P.is_equal_up_to_relabeling(P3)
True
sage: Re1 == (beta * R2e1).inverse()
True
sage: Re2 == (beta * R2e2).inverse()
True
\end{sagecommandline}
The self-similarity computed by this Rauzy induction is the product of the above
2-dimensional substitutions by the bijection of the labels.
\begin{sagecommandline}
sage: from slabbe import Substitution2d
sage: s3 = Substitution2d.from_permutation(P.keys_permutation(P3))
sage: s123 = s1*s2*s3
\end{sagecommandline}

The computed self-similarity \texttt{s123} is:
\begin{center}
\resizebox{.95\textwidth}{!}{
$\sagestr{s123._latex_(ncolumns=6)}$
}
\end{center}

The above self-similarity can be illustrated with the Wang tiles computed above
as follows:
\begin{sagecommandline}
sage: s123_tikz = s123.wang_tikz(domain_tiles=T3, codomain_tiles=T3, ncolumns=6, scale=1.2, label_shift=.15)
\end{sagecommandline}

\begin{center}
    \sageplot[width=.95\linewidth][pdf]{s123_tikz}
\end{center}

We may observe that the self-similarity computed here from the Rauzy induction on
polygonal partition on $\Pcal_3$ and toral $\Z^2$-action $R_3$ is the same as
the self-similarity proved for the Wang shift $\Omega_3$ in
\cite{labbe_metallic_I_2025}.

\section{Open questions}
\label{sec:open-questions}

\begin{figure}[h]
\begin{center}
    \includegraphics{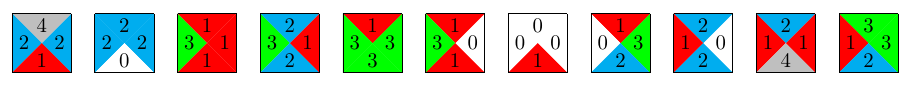}
\end{center}
    \caption{The Jeandel--Rao aperiodic set of 11 Wang tiles.}
    \label{fig:JR-11-tiles}
\end{figure}

For almost twenty years, the Kari and Culik sets of Wang tiles were the
smallest known aperiodic sets of Wang tiles.
In 2015, Jeandel and Rao performed an exhaustive search on all sets of Wang tiles of
cardinality up to 11 \cite{zbMATH07421483}
and proved that sets of Wang tiles of cardinality at most 10 either do not
tile the plane or tile the plane and one of the valid tilings is periodic.
Moreover, they provided a list of 36 sets of 11 Wang tiles considered to be
candidates for being aperiodic. One of candidates was intriguing because
Fibonacci numbers appeared in the structure of the transducers involved
in the computation of valid tilings. Jeandel and Rao focused on the intriguing candidate,
shown in Figure~\ref{fig:JR-11-tiles}, and they proved it to be aperiodic.
The set of valid configurations over these 11 tiles forms a subshift that we
call the Jeandel--Rao Wang shift.

In \cite{MR4213162}, it was proved 
that a minimal subshift within the Jeandel--Rao Wang shift 
is the coding of a dynamical system
defined by the following $\Z^2$-action $R_0$
on the $2$-dimensional torus $\R^2/\Gamma_0$,
where $\Gamma_0=\left(\begin{smallmatrix} \varphi & 1\\ 0 & \varphi+3\end{smallmatrix}\right)\Z^2$
is a lattice in $\R^2$ involving the golden ratio
$\varphi=\frac{1+\sqrt{5}}{2}$:
\[
\begin{array}{rccl}
    R_0:&\Z^2\times\R^2/\Gamma_0 & \to & \R^2/\Gamma_0\\
    &(\bk,\bx) & \mapsto &\bx+\bk.
\end{array}
\]
The symbolic coding is obtained through a polygonal partition $\Pcal_0$ of a fundamental domain
of $\R^2/\Gamma_0$.
The partition was proved to be a Markov partition for $R_0$
after comparing the substitutive structure computed from the Rauzy induction of
$R_0$ and $\Pcal_0$ \cite{MR4347332} with the substitutive structure of the
associated Wang shift \cite{MR3978536,MR4226493}.

Intuitively, this means that the Jeandel--Rao Wang tiles shown in
Figure~\ref{fig:JR-11-tiles} correspond to computing the orbit of points in the plane
$\R^2$ under the translations by $+1$ horizontally and $+1$
vertically modulo the lattice $\Gamma_0$.
How come this is possible is still a mystery.
The link between the 11 Jeandel--Rao Wang tiles themselves
and the golden ratio or toral rotation $R_0$ remains unclear.
Unlike the Kari example, the values 0, 1, 2, 3, 4 of the labels of the
Jeandel--Rao Wang tiles are five distinct symbols rather than arithmetic values.
They do not satisfy a known equation. 

In general, the following questions can be raised.

\begin{mainquestion}\label{question:finding-JR-values}
    Let $\Tcal$ be a set of Wang tiles such that the Wang shift $\Omega_\Tcal$
    is aperiodic.
    \begin{itemize}
        \item Is it multiplicative (Kari-Culik-like)? More precisely,
    can we replace the labels of the tiles in $\Tcal$
    by arithmetic values in such a way that an 
    equation similar to \eqref{eq:wang-kari} is satisfied?
        \item Is it additive (metallic mean-like)? More precisely,
        can we replace the labels of the tiles in $\Tcal$ by integer vectors 
        computed from floors of linear forms as in
        Proposition~\ref{prop:Tcaln_is_image_of_sctilexy}
        and satisfying additive equations as in
            Theorem~\ref{thm:equations-satisfied-by-tiles}?
    \end{itemize}
    Does there exists an aperiodic set of Wang tiles which is neither multiplicative 
    nor additive?
\end{mainquestion}

Solving Question~\ref{question:finding-JR-values} for Jeandel--Rao Wang tiles
would improve our understanding of the Jeandel--Rao Wang shift. 
Hopefully it would allow to generate more examples maybe not related to the golden ratio
and that are not self-similar. Remember that the computations made by Jeandel
and Rao took one year using 100 cpus to explore exhaustively the sets of 11
Wang tiles \cite{zbMATH07421483}. Finding new examples by exploring all sets
of 12, 13 or 14 Wang tiles becomes soon out of reach. We need to understand
what is happening in order to find other examples and characterize them.

\begin{mainquestion}\label{question:arithmetical-aperiodic-proof}
    If an aperiodic set of Wang tiles is additive (metallic mean-like) with
    labels given by integer vectors satisfying equations,
    can we use the equations to directly prove that the Wang shift
    $\Omega_\Tcal$ is aperiodic following the short arithmetical argument for
    the nonperiodicity of Kari's tile set?
\end{mainquestion}

Finding an answer to Question~\ref{question:arithmetical-aperiodic-proof}
for the Ammann set of 16 Wang tiles was the original motivation
of the author which led to the discovery of the family of metallic mean Wang
tiles. As we discussed in Section~\ref{sec:equations},
Question~\ref{question:arithmetical-aperiodic-proof}
remains open even for the Ammann 16 Wang tiles and the family of metallic mean
Wang tiles.

In general, we may ask the following question.

\begin{mainquestion}\label{question:find-the-wang-tiles}
    For which invertible matrix $M\in\mathrm{GL}_2(\R)$
    does there exist a set of Wang tiles $\Tcal$ such
    that the Wang shift $\Omega_\Tcal$ is isomorphic,
    as a measure-preserving dynamical system,
    to the toral $\Z^2$-rotation
    $R:\Z^2\times\torus^2\to\torus^2$
    defined by $R^\bk(\bx)=\bx+M\bk$
    on the 2-dimensional torus $\torus^2=(\R/\Z)^2$?
\end{mainquestion}

The Markov partition associated with Jeandel--Rao tiles and action $R_0$ on $\R^2/\Gamma_0$
is related to the golden ratio \cite{MR4213162}. In this contribution,
we describe a family of $\Z^2$-actions related to the metallic-mean quadratic integers.
Can we find examples related to other numbers?

\begin{mainquestion}
For which $\Z^2$-actions defined by rotations on a $2$-dimensional torus does
    there exist a Markov Partition? When is this partition smooth/polygonal?
\end{mainquestion}

As for toral hyperbolic automorphisms, we can expect that smooth Markov partitions
are associated with algebraic integers of degree 2 and that the partition is
piecewise linear in this case \cite{MR1145614}. Markov partitions for typical
toral hyperbolic automorphisms have fractal boundaries \cite{MR474415}.

The relation with toral hyperbolic automorphisms does not come out of nowhere.
Indeed, the self-similarity of $\Omega_n$ proved in
\cite{labbe_metallic_I_2025} has an incidence matrix of size $(n+3)^2\times(n+3)^2$.
Its eigenvalues are all quadratic integers, 0 or $\pm 1$.
This incidence matrix acts hyperbolically as a toral automorphism on a subspace
of $\R^{(n+3)^2}$ thus admits a Markov partition with piecewise linear boundaries.
A link between this Markov partition and the partition $\Pcal_n$ can be expected,
because this is what happens for $1$-dimensional sequences. Indeed, 
the Markov partition associated with the toral automorphism 
$\left(\begin{smallmatrix}1&1&1\\1&0&0\\0&1&0 \end{smallmatrix}\right)$
is a suspension of the Rauzy fractal \cite{MR667748} as nicely illustrated
in a talk by Timo Jolivet \cite{jolivet_toral_2012}.

\begin{mainquestion}
    What is the relation between the Markov partition for the hyperbolic toral automorphism
    defined from the incidence matrix of the self-similarity of $\Omega_n$ and
    the Markov partition $\Pcal_n$ associated with
    $\dynsys{\Z^2}{\sigma}{\Omega_n}$?
\end{mainquestion}

The symmetric properties of $\Omega_n$ and of the partition $\Pcal_n$ make them
a good object of study to tackle these questions in more generality.

\subsection*{Acknowledgments}
The author would like to thank the referees for the extensive and helpful
comments that have helped to improve the presentation.
The author is also thankful to Hugo Parlier for his comments on an earlier draft
of the introduction and to Vincent Delecroix for making the author realize
that it is not the Birkhoff ergodic theorem which is needed in the proof of 
Proposition~\ref{prop:phi-is-a-factor-map}
but rather simply the Weyl's equidistribution theorem.

\subsection*{Competing interest}

The authors have no competing interests to declare.

\subsection*{Financial Support}
This work was partly funded from France's Agence Nationale de la Recherche
(ANR) projects CODYS (ANR-18-CE40-0007) and IZES (ANR-22-CE40-0011).
It was also supported by grants from the 
\emph{Symbolic Dynamics and Arithmetic Expansions}
(SymDynAr) Project, co-funded by ANR 
(ANR-23-CE40-0024)
and FWF (\href{https://dx.doi.org/10.55776/I6750}{I 6750}),
the Austrian Science Fund.

\subsection*{Reproducibility statement}

All results proved in this article are proved by hand.
Computations performed in Section~\ref{sec:renormalization-rauzy}
are based on the open-source
mathematical software SageMath \cite{sagemathv10.6} and the
optional package \texttt{slabbe} \cite{labbe_slabbe_0_8_0_2025}.
All SageMath input/output blocks in this article were created using
the \texttt{sageexample} environment with 
SageTeX version \texttt{2021/10/16 v3.6}
and with the following software versions:
\begin{sagecommandline}
sage: version()
'SageMath version ..., Release Date: ...'
sage: import importlib.metadata
sage: importlib.metadata.version("slabbe")
'...'
\end{sagecommandline}
The fact that these software are open-source means that anyone
is free to use, reproduce, verify, adapt for their own needs all of the
computations performed therein according to the
GNU General Public License (version 2, 1991, \url{http://www.gnu.org/licenses/gpl.html}).

The contents of all of the \texttt{sageexample} environments from the tex
source are gathered in the file
\texttt{demos/arXiv\_2403\_03197\_doctest.sage}
autogenerated by SageTeX when running \texttt{pdflatex}.
This file is included in the \texttt{slabbe} package and available at
\url{https://gitlab.com/seblabbe/slabbe/}.
It allows to make sure that future releases of the package do not break the
code included in this article.
It is possible to reproduce all computations present in this article and 
check that all outputs are correct, by \emph{doctesting} this file, that is,
by running the command
\texttt{sage -t demos/arXiv\_2403\_03197\_doctest.sage}.
It should output \texttt{All tests passed!} and 
\texttt{[58 tests, 11.75s wall]} (most probably with a different timing).

\bibliographystyle{alpha-first-name-initials-with-doi} %
\bibliography{biblio}

\end{document}